# INTERVAL SEMIRINGS

W. B. Vasantha Kandasamy
Florentin Smarandache

**2011**

# INTERVAL SEMIRINGS

**W. B. Vasantha Kandasamy**
**Florentin Smarandache**

**2011**

# CONTENTS

# PREFACE

In this book the notion of interval semirings are introduced. The authors study and analyse semirings algebraically. Methods are given for the construction of non-associative semirings using loops and interval semirings or interval loops and semirings. Another type of non-associative semirings are introduced using groupoids and interval semirings or interval groupoids and semirings. Examples using integers and modulo integers are given. Also infinite semirings which are semifields are given using interval semigroups and semirings or semigroups and interval semirings or using groups and interval semirings.

Interval groups are introduced to construct interval group interval semirings, and properties related with them are analysed. Interval matrix semirings are introduced.

This book has seven chapters. In chapter one we give the basics needed to make this book a self contained one. Chapter two introduces the notion of interval semigroups and interval semifields and are algebraically analysed. Chapter three introduces special types of interval semirings like matrix interval semirings and interval polynomial semirings. Chapter four for the first time introduces the notion of group interval semirings, semigroup interval semirings, loop interval semirings and groupoid interval semirings and these structures are studied. Interval neutrosophic semirings are introduced in chapter five. Applications of these structures are given in chapter six. The final chapter suggests around 120 problems for the reader. We thank Dr. K.Kandasamy for proof reading.

<div style="text-align: right;">
W.B.VASANTHA KANDASAMY  
FLORENTIN SMARANDACHE
</div>

**Chapter One**

# BASIC CONCEPTS

In this book we introduce several types of semirings using intervals, mainly of the form [0, a]; a real. Study of such type of interval semirings is carried out for the first time in this book.

A semiring is a nonempty set S endowed with two binary operations "+" and "." such that (S, +) is an additive abelian semigroup and (S, .) is a semigroup and both right and the left distributive laws are satisfied.

*Example 1.1:* The set of positive integers with zero is a semiring.

*Example 1.2:* Set of positive reals with zero is a semiring. A semiring is said to be commutative if (S, .) is a commutative semigroup.

The semirings given in examples 1.1 and 1.2 are commutative semirings.

A semiring S is said to be a strict semiring if for a and b in S a + b = 0 implies a = 0 or b = 0.

Semirings mentioned in the above examples are strict semirings. A semiring S is said to be a semifield if S is a strict commutative semiring with the multiplicative identity 1. Semirings given in examples 1.1 and 1.2 are semifields. A semiring is said to be of infinite cardinality if it has infinite number of elements in them. A finite semiring is one which has finite number of elements in them. All distributive lattices are semirings. If the distributive lattices are of finite order we get finite semirings. All Boolean algebras are semirings. All finite Boolean algebras of order greater than two are only commutative semirings and are not strict semirings, hence they are not semifields. All chain lattices are semifields.

For more about semirings and semifields please refer [13].

$I(Z^+ \cup \{0\})$ denotes the collection all intervals [0, a] with a in $Z^+ \cup \{0\}$. $I(Q^+ \cup \{0\}) = \{[0, a] \mid a \in Q^+ \cup \{0\} \}$ and $I(R^+ \cup \{0\})$ denotes the collection of all [0, a] with a in $R^+ \cup \{0\}$.

I denotes the indeterminancy and $I^2 = I$. All intervals of the form [0, a+ bI] denotes the neutrosophic interval with a and b reals.

For more about neutrosophic concepts please refer [5-7].

**Chapter Two**

# INTERVAL SEMIRINGS AND INTERVAL SEMIFIELDS

In this chapter for the first time we define the new notions of interval semirings and interval semifields we describe these structures with examples. These structures will be used in the later chapters to construct neutrosophic interval semirings and other types of interval semialgebras. We will be using the notations given in chapter one.

**DEFINITION 2.1:** *Let $S = \{[0, a], +, . \}$ be a collection of intervals from $Z_n$ (or $Z^+ \cup \{0\}$ or $R^+ \cup \{0\}$ or $Q^+ \cup \{0\}$). S is a semigroup under '+' which is commutative with the additive identity and $(S, .)$ is a semigroup under multiplication. For $[0, a], [0, b], [0, c]$ in S if we have*

$$\begin{aligned}[0, a] . ([0, b] + [0, c]) &= [0, a] . [0, b] + [0, a] . [0, c] \\ &= [0, ab] + [0, ac] \\ &= [0, ab + ac] \\ &= [0, a(b + c)]\end{aligned}$$

*for all elements in S, then S is defined to be a interval semiring.*

We will give examples of interval semirings.

*Example 2.1:* Let $S = \{[0, a] \mid a \in Z_5\}$; $(S, +, .)$ is an interval semiring of finite order.

*Example 2.2:* Let $S = \{[0, a] \mid a \in Z^+ \cup \{0\}\}$; $(S, +, .)$ is an interval semiring of infinite order.

*Example 2.3:* Let $P = \{[0, a] \mid a \in Q^+ \cup \{0\}\}$ is an interval semiring of infinite order.

We can as in case of usual semirings define interval subsemirings.

**DEFINITION 2.2:** *Let $V = \{[0, a] / a \in Z_n$ or $Z^+ \cup \{0\}$ or $Q^+ \cup \{0\}$ or $R^+ \cup \{0\}\}$ be an interval semiring under addition and multiplication. Let $P \subseteq V$, if $(P, +, .)$ is an interval semiring ($P$ a proper subset of $V$) then we define $P$ to be a interval subsemiring of $P$.*

If in an interval semiring $(S, +, .)$, the interval semigroup $(S, .)$ is commutative then we define $S$ to be a commutative interval semiring. If $(S, +, .)$ is such that $(S, .)$ is not a commutative interval semigroup, then we define $(S, +, .)$ to be a non commutative interval semigroup.

We will give examples of them.

*Example 2.4:* Let $(S, +, .) = \{[0, a] \mid a \in Z_{15}, +, .\}$ be an interval semiring. Consider $(P, +, .) = \{[0, a] \mid a \in \{0, 3, 6, 9, 12\} \subseteq Z_{15}, +, .\} \subseteq (S, +, .)$; $P$ is an interval subsemiring of $S$. Both $S$ and $P$ are of finite order.

*Example 2.5:* Let $(S, +, .) = \{[0, a] \mid a \in Z^+ \cup \{0\}, +, .\}$ be an interval semiring. $(P, +, .) = \{[0, a] \mid a \in 3Z^+ \cup \{0\}, +, .\} \subseteq (S, +, .)$ is an interval subsemiring of $(S, +, .)$.

*Example 2.6:* Let $(S, +, .) = \{[0, a] \mid a \in Q^+ \cup \{0\}, +, .\}$ be an interval semiring. Take $(P, +, .) = \{[0, a] \mid a \in Z^+ \cup \{0\}, +, .\}$
$(P, +, .)$ is an interval subsemiring.

***Example 2.7:*** Let $(S, +, .) = \{[0, a] \mid a \in R^+ \cup \{0\}, +, .\}$ be an interval semiring. Consider $(W, +, .) = \{[0, a] \mid a \in R^+ \cup \{0\}, +, .\} \subseteq (S, +, .)$, W is an interval subsemiring of S.

All the examples of interval semirings given are interval commutative semirings. If an interval semiring S has an interval [0, 1] such that [0, 1] acts as identity that is [0, 1] . [0, a] = [0, a] . [0, 1] = [0, a] for all [0, a] $\in$ S. then we call S an interval semiring with identity.

All interval semirings in general need not be an interval semiring with identity. To this effect we will give some examples.

***Example 2.8:*** Let $(S, +, .) = \{[0, a] \mid a \in 2Z^+ \cup \{0\}, +, .\}$ be an interval semiring. $(S, +, .)$ is not an interval semiring with identity.

***Example 2.9:*** Let $(S, +, .) = \{[0, a] \mid a \in Q^+ \cup \{0\}, +, .\}$ be an interval semiring . S is an interval semiring with identity. We call an interval semiring S to be of characteristic 0 if n[0, a] = 0 implies n = 0, where a $\neq$ 0 and a $\in Z^+ \cup \{0\}$ or $Q^+ \cup \{0\}$ or $R^+ \cup \{0\}$. If in case n [0, a] = 0 still n $\neq$ 0 and a $\neq$ 0 then we say the interval semiring S is of characteristic n, with m[0, a] = 0 for no m < n.

We will illustrate this situation before we proceed on to define different types of interval semirings.

***Example 2.10:*** Let $S = \{[0, a], +, . \mid a \in R^+ \cup \{0\}\}$ be an interval semiring of characteristic zero.

***Example 2.11:*** Let $V = \{[0, a], +, ., Z_{12}\}$ be an interval semiring of characteristic twelve.

***Example 2.12:*** Let $P = \{[0, a], +, ., Z_7\}$ be an interval semiring of characteristic seven.

***Example 2.13:*** Let $V = \{[0, a], +, ., 3Z^+ \cup \{0\}\}$ be an interval semiring of characteristic zero.

We say an interval semiring $S = \{[0, a], +, .\}$ has zero divisors if $[0, a] . [0, b] = [0, ab] = [0, 0] = 0$ where $a \neq 0$ and $b \neq 0$.

We will give some examples of interval semirings with zero divisors.

*Example 2.14:* Let $S = \{[0, a], +, . \mid a \in Z_{18}\}$ be an interval semiring. Consider $[0, 6], [0, 3]$ in S; $[0, 6] [0, 3] = [0, 18] = [0, 0]$ is a zero divisor in S. Also for $[0, 2], [0, 9] \in S$; we have $[0, 2].[0, 9] = [0, 18] = [0, 0] = 0$ is again a zero divisor. We have several other zero divisors in S.

*Example 2.15:* Let $S = \{[0, a], +, . \mid a \in Z^+ \cup \{0\}$ be an interval semiring. S has no zero divisors.

*Example 2.16:* Let $S = \{[0, a], +, . \mid a \in Z_{41}\}$ be an interval semiring S has no zero divisors.

**THEOREM 2.1:** *All interval semirings $S = \{[0, a], +, . \mid a \in Z_p,$ p a prime\} have no zero divisors.*

The proof is straight forward and hence is left as an exercise for the reader.

**THEOREM 2.2:** *All interval semirings $S = \{[0, a], +, . \mid a \in Z_n;$ n a composite number\} have zero divisors.*

This proof is also left as an exercise for the reader to prove.

**THEOREM 2.3:** *Let $S = [[0, a], +, . \mid a \in Z^+ \cup \{0\}$ or $Q^+ \cup \{0\}$ or $R^+ \cup \{0\}$ or $nZ^+ \cup \{0\}$; n >0. n = 1, 2, ..., $\infty$\} be an interval semiring . S has no zero divisors.*

The proof is left as an exercise for the reader. We say an interval semiring S has idempotents if $[0, a] [0, a] = [0, a]$; $[0, a] \in S$; $a \neq 0$ or $[0, 1]$.

We will illustrate this situation by some examples.

***Example 2.17:*** Let $S = \{[0, a], +, . \mid a \in Z_{12}\}$ be an interval semiring $[0, 4] \in S$ is such that $[0, 4] . [0, 4] = [0, 16(\text{mod } 12)] = [0, 4]$; so $[0, 4]$ is an idempotent of S. $[0, 9] \in S$ is also an idempotent for $[0, 9].[0, 9] = [0, 81(\text{mod } 12)] = [0, 9]$.

Thus S has nontrivial idempotents.

***Example 2.18:*** Let $S = \{[0, a], +, . \mid a \in Z_p\}$ where p is a prime be interval semiring S has no idempotents and zero divisors. However S has units.

In view of this we have the following theorem.

**THEOREM 2.4:** *Let $S = \{[0, a] \mid a \in Z_p, p$ a prime; $+, . \}$ be an interval semiring, S has no idempotents and zero divisors.*

The proof is straightforward and is left as an exercise for the reader.

***Example 2.19:*** Let $S = \{[0, a] \mid a \in Z_{23}, +, .\}$ be an interval semiring. S has non trivial units. For take $[0, 8] [0, 3] = [0, 8.3 (\text{mod } 23)] = [0, 1]$ is a unit of S. $[0, 12] [0, 2] = [0, 12.2 (\text{mod } 23)] = [0, 1]$ is a unit of S. Consider $[0, 6] [0, 4] = [0, 6.4(\text{mod } 23)] = [0, 1]$ is a unit of S. Now $[0, 22] [0, 22] = [0, 48.4 (\text{mod} 23)]$ is a unit in S. Thus S has several units but has no idempotents or zero divisors.

**THEOREM 2.5:** *Let $S = \{[0, a] / a \in Z_p, p$ a prime, $+, .\}$ be an interval semiring. S has units.*

The proof is also left as an exercise to the reader.

Now we can have interval semirings which are commutative or non commutative. Those examples of interval semirings built using $Z_n$ or $Z^+ \cup \{0\}$ or $Q^+ \cup \{0\}$ or $R^+ \cup \{0\}$ are commutative. To construct interval semirings we have to proceed differently which will be carried out in chapter three of this book.

Now we proceed on to define interval semifields.

**DEFINITION 2.3:** *Let $S = \{[0, a] \mid a \in Z^+ \cup \{0\}$ or $R^+ \cup \{0\}\}$ with operation + and . so that S is an interval commutative semiring. If (a) $[0, a] + [0, b] = [0\ 0] = 0$ if and only if $[0, a] = [0, 0]$ and $[0, b] = [0, 0]$. That is the S is a strict interval semiring and (b) in S, $[0, a] . [0, b] = [0, 0]$ then $[0, a] = [0, 0]$ or $[0, b] = [0, 0]$ then we define the interval semiring to be an interval semifield. If in the interval semifield $0 . [0, x] = [0, 0]$ and for no integer n, $n [0, x] = 0$ or equivalently if $[0, x] \in S \setminus \{0\} = [0, 0]$, $n [0, x] = [0, x] + [0, x] + [0, x] + \ldots + [0, x]$, n times equal to is zero is impossible for any positive integer n; then we say the semifield S is of characteristic 0.*

We will illustrate this by some examples.

***Example 2.20:*** Let $S = \{[0, a], +, . \mid a \in Z^+ \cup \{0\}\}$ be a interval semifield of characteristic zero.

***Example 2.21:*** Let $S = \{[0, a], a \in R^+ \cup \{0\}; +, .\}$ be an interval semifield of characteristic zero.

***Example 2.22:*** Let $S = \{[0, a]; a \in Q^+ \cup \{0\}, +, .\}$ be an interval semifield of characteristic zero.

Now having seen examples of interval semifields we now proceed on to show that all interval commutative semirings are not interval semifields.

***Example 2.23:*** Let $S = \{[0, a], a \in Z_n, n < \infty, +, .\}$ be an interval semiring. Take n = 6, clearly S is not an interval semifield. For consider
$$[0, 2] + [0, 2] + [0, 2] = 3[0, 2]$$
$$= [0, 0].$$
Also $[0, 3] + [0, 3] = [0, 0]$; thus S is not an interval semifield. Further $[0, 2] . [0, 3] = [0, 2.3 \pmod 6] = [0, 0]$. So conditions (a) and (b) of the definition 2.3 is not true.

***Example 2.24:*** Let $S = \{[0, a], +, .; a \in Z_{11}\}$ be an interval semiring. We see $[0, 6] + [0, 5] = [0, (5+6) \bmod 11] = [0, 0]$, Clearly $[0, 6] \neq [0, 0]$ and $[0, 5] \neq [0, 0]$. So S is not a strict interval semiring. So S is not an interval semifield. However S has no zero divisors. For in S we do not have $[0, a] \cdot [0, b] = [0, 0]$ if $a, b \in Z_{11} \setminus \{0\}$.

*Example 2.25:* Let $S = \{[0, a] \mid a \in Z_{15}, +, \cdot\}$ be an interval semiring. We see S is not an interval semifield. For we have $[0, 11] + [0, 4] = [0, (11+4) \mod (15)] = [0, 0]$ so condition (a) of the definition 2.3 is not true.

Further $[0, 6] \cdot [0, 5] = [0, 6.5 \pmod{15}] = [0, 0]$ is a nontrivial zero divisor in S. Infact we have several zero divisors in S. Thus S is not an interval semifield.

*Example 2.26:* Let $S = \{([0, a], [0, b], [0, c]) \mid a, b, c \in Z^+ \cup \{0\}, +, \cdot\}$ be an interval semiring under component wise addition and component wise multiplication. That is if $x = ([0, a], [0, b], [0, c])$ and $y = ([0, t], [0, u], [0, w])$ where $x, y \in S$ and $a, b, c, u, t, w \in Z^+ \cup \{0\}$.

We see

$x + y$ = $([0, a], [0, b], [0, c]) + ([0, t], [0, u], [0, w])$
= $([0, a + t], [0, b + u], [0, c + w])$.

$x \cdot y$ = $([0, a], [0, b], [0, c]) ([0, t], [0, u], [0, w]$
= $([0, a] [0, t], [0, b] [0, u], [0, c] [0, w])$
= $([0, at], [0, bu], [0, cw])$.

We see S has zero divisors for take $p = ([0, 0], [0, b], [0, 0])$ and $q = ([0, a] [0, 0], [0, c])$ ; We have

p.q = $([0, 0] [0, b], [0, 0] \cdot ([0, a], [0, 0], [0, c])$
= $([0, 0] [0, a], [0, b] [0, 0], [0, 0] [0, c])$
= $([0, 0.a], [0, b.0], [0, 0.c])$
= $([0, 0], [0, 0], [0, 0])$.

Thus S is not an interval semifield only an interval semiring. Thus we can construct by this method infinitely many interval semirings which are not interval semifields using $Z^+ \cup \{0\}$ or $Z_n$ or $Q^+ \cup \{0\}$ or $R^+ \cup \{0\}$.

**THEOREM 2.6:** *Let $S = \{([0, a_1], [0, a_2], ..., [0, a_n]) \mid a_i \in Z^+ \cup \{0\}$ or $a_i \in R^+ \cup \{0\}$ or $a_i \in Q^+ \cup \{0\}; 1 \leq i \leq n; +, .\}$ be an interval semiring. S is not an interval semifield.*

This proof is left as an exercise for the reader.
Now we can as in case of other algebraic structures define the notion of interval subsemifield, this task is also left to the reader. However we give examples of them.

***Example 2.27:*** Let $S = \{[0, a], +, \text{'.'}$ where $a \in Q^+ \cup \{0\}\}$ be an interval semifield. Take $T = \{[0, a], +, .$ where $a \in Z^+ \cup \{0\}\} \subseteq S$; T is an interval subsemifield of S.

***Example 2.28:*** Let $P = \{[0, a], +, \text{'o'}$ where $a \in R^+ \cup \{0\}\}$ be an interval semifield. Take $W = \{[0, a], +, \text{'o'}$ where $a \in Z^+ \cup \{0\}\} \subseteq P$; W is an interval subsemifield of P. Infact P has many interval subsemifields.

Now having seen examples of interval subsemifields we now proceed on to define the notion of direct product of interval semirings.

**DEFINITION 2.4:** *Let $S_1$ and $S_2$ be two interval semirings. The direct product $S_1 \times S_2 = \{([0, a], [0, b]) \mid [0, a] \in S_1$ and $[0, b] \in S_2\}$ is also an interval semiring with component wise operation. Similarly if $S_1, S_2, ..., S_t$ are t interval semirings. The direct product of these interval semirings denoted by $S_1 \times S_2 \times ... \times S_t = \{([0, a_1], [0, a_2], ..., [0, a_t]) \mid [0, a_i] \in S_i; 1 \leq i \leq t\}$ is an interval semiring.*

We will illustrate this situation by some examples.

***Examples 2.29:*** Let $S_1 = \{[0, a] \mid a \in Z_{12}; +, .\}$, $S_2 = \{[0, b] \mid b \in Z^+ \cup \{0\}, +, .\}$ and $S_3 = \{[0, c] \mid c \in Z_{19}, +, .\}$ be three interval semirings. The direct product of these interval semirings $S_1 \times S_2 \times S_3 = \{([0, a], [0, b], [0, c]) \mid [0, a] \in S_1, [0, b] \in S_2$ and $[0, c] \in S_3\}$ is an interval semiring.

We see $S = \{([0, a_1], \ldots, [0, a_n]) \mid a \in Z_m, +, .\}$ is an interval semiring; This interval semiring can be visualized as the direct product of n interval semirings. $P = \{[0, a] \mid a \in Z_m, +, .\}$ That is $\underbrace{P \times P \times \ldots \times P}_{n-\text{times}} = S$.

We can define ideals in an interval semiring S.

**DEFINITION 2.5:** *Let $S = \{[0, a] \mid a \in Z_n \text{ or } Z^+ \cup \{0\}, +, .\}$ be an interval semiring. Choose $I = \{[0, t] \mid t \in Z^+ \cup \{0\} \text{ or } Z_n\} \subseteq S$; I is an ideal of S if*
1. *I is an interval subsemiring and*
2. *For all $[0, t] \in I$ and $[0, s] \in S$ $[0, t] . [0, s]$ and $[0, s] . [0, t]$ are in I.*

We can have the notion of left ideals and right ideals of an interval semiring. When both right ideal and left ideal coincide for a non empty subset I of S we call it as an ideal of S.

We will illustrate this situation by some examples.

*Example 2.30:* Let $S = \{[0, a] \mid a \in Z^+ \cup \{0\}, +, \text{`·'}\}$ be an interval semiring. Let $I = \{[0, a] \mid a \in 3Z^+ \cup \{0\}, +, \text{`·'}\} \subseteq S$; I is an interval ideal of the interval semiring S.

*Example 2.31:* Let $S = \{[0, a] \mid a \in Q^+ \cup \{0\}, +, \text{`o'}\}$ be an interval semirings. It is easily verified S has no interval ideals, but S has nontrivial interval subsemirings. For $T = \{[0, a] \mid a \in nZ^+ \cup \{0\}, +, \text{`o'}\} \subseteq S$ for $n = 1, 2, \ldots, t$, $t < \infty$ are interval subsemirings of S which are not ideals.

*Example 2.32:* Let $S = \{[0, a] \mid a \in R^+ \cup \{0\}, +, \text{`o'}\}$ be an interval semiring. S has no interval ideals but has infinitely many interval subsemirings. For example take $P = \{[0, a] \mid a \in Q^+ \cup \{0\}, +, \text{`o'}\}$, S is an interval subsemiring and is not an ideal of S. Take $V_n = \{[0, a] \mid a \in nZ^+ \cup \{0\}; +, \text{`o'}\} \subseteq S; n < \infty$ but take all positive integer values; $V_n$ is an interval subsemiring but is not an ideal of S.

***Example 2.33:*** Let $W = \{[0, a] / a \in Z_{19}, +, \cdot\}$ be an interval semiring. Clearly W has no nontrivial ideals or interval subsemiring. It is interesting to note that W is not a strict semiring for we see

$[0, 10] + [0, 9] = [0, 10 + 19 \pmod{19}] = [0, 0]$.
$[0, 12] + [0, 7] \equiv 0 \pmod{19}$
$[0, 5] + [0, 14] \equiv 0 \pmod{19}$

and so on. However W has no zero divisors. Infact W is a finite ring. All rings are semirings.

***Example 2.34:*** Let $W = \{[0, a] | a \in Z_{25}, +, \cdot\}$ be an interval semiring. W has interval ideals as well as interval subsemirings. For consider $T = \{[0, a] | a \in \{0, 5, 10, 15, 20\} \subseteq Z_{25}, +, \cdot\} \subseteq W$ is an interval ideal so is also an interval subsemiring. W has no other ideals or subsemirings.

***Example 2.35:*** Let $M = \{[0, a] / a \in Z_{16}; +, \cdot\}$ be an interval semiring. $N = \{[0, a] / a \in \{0, 2, 4, 6, 8, 10, 12, 14\} \subseteq Z_{16}, +, \cdot\} \subseteq M$ is an interval ideal of M.

***Example 2.36:*** Let $S = \{[0, a] / a \in Z_{24}, +, \cdot\}$ be an interval semiring. Take $T = \{[0, a] / a \in \{0, 2, \ldots, 22\} \subseteq Z_{24}, +, \cdot\} \subseteq S$, T is an interval subsemiring as well as an interval ideal of S. $V = \{[0, a] / a \in \{0, 6, 12, 18\} \subseteq Z_{24}, +, \cdot\} \subseteq S$ is an interval subsemiring, as well as an ideal of S. $R = \{[0, a] / a \in \{0, 8, 16\} \subseteq Z_{24}, +, \cdot\} \subseteq S$, R is an interval subsemiring and an interval ideal of S.

Now having seen ideals and subsemirings of interval semirings we now proceed on to define direct product of interval semirings using interval rings and interval semirings so that the resultant is an interval semiring. We know all interval semirings built using $Z_n$, $n < \infty$ are only interval semirings.

***Example 2.37:*** Let $W = W_1 \times W_2 \times W_3 \times W_4 = \{[0, a] / a \in Z_5, +, \cdot\} \times \{[0, a] | a \in 3Z^+ \cup \{0\}, +, \cdot\} \times \{[0, a] / a \in Z_7, +, .\} \times \{[0, a] | a \in 7Z^+ \cup \{0\}\}$; W is an interval semiring and is not a ring. W has both ideals and subsemirings.

***Example 2.38:*** Let $V = V_1 \times V_2 \times V_3 \times V_4 \times V_5 \times V_6 = \{[0, a] / a \in Q^+ \cup \{0\}, +, \cdot\} \times \{[0, a] / a \in Z_{12}, +, \cdot\} \times \{[0, a] | a \in Z_{19}, +, \cdot\} \times \{[0, a] | a \in Z_3, +, \cdot\} \times \{[0, a] | a \in Z_{120}, +, \cdot\} \times \{[0, a] | a \in Z_{43}, +, \cdot\}$ be an interval semiring, V has both nontrivial interval ideals as well as nontrivial subsemigroups. It is clear none of these are interval semifields.

Now we can as in case of other algebraic structures define homomorphisms.

**DEFINITION 2.6:** *Let S and S' be any two interval semirings. An interval mapping $\phi : S \to S'$ is called the interval semiring homomorphism if $\phi([0, a]+[0, b]) = \phi([0, a]) + \phi([0, b])$ and $\phi([0, a] \cdot [0, b]) = \phi([0, a]) * \phi([0, b])$ where * is an operation on S'; this is true for all [0, a] and [0, b] $\in$ S.*

We will illustrate this situation by some examples.

***Example 2.39:*** Let $S = \{[0, a] | a \in 3Z^+ \cup \{0\}, +, \cdot\}$ and $S' = \{[0, a] | a \in Q^+ \cup \{0\}, +, \cdot\}$ be any two interval semirings. Define $\phi : S \to S'$ by $\phi([0, a]) = [0, a]\ \forall\ [0, a]\ 3Z^+ \cup \{0\}$. It is easily verified $\phi$ is an interval semiring homomorphism, that is $\phi$ is an embedding of S in to S' as clearly $S \subseteq S'$. Infact we see $\phi$ is an embedding of S in to S' as $S \subseteq S_1$. We see as in case of other homomorphisms, we can in case of interval semiring homomorphisms $\eta$ also define kernel of $\eta$.

Let $\eta: S \to S'$ where S and S' are interval semirings and $\eta$ an interval semiring homomorphism. Ker $\eta = \{[0, a] \in S / \eta [0, a] = [0, 0] = 0\}$.

We will call an interval semiring S to be a Smarandache interval semiring S-interval semiring if S has a proper subset $P \subseteq S$ such that P is an interval semifield.

We will illustrate this by some example.

***Example 2.40:*** Let S = {[0, a] | a ∈ $Q^+ \cup \{0\}$} be an interval semiring. Consider P = {[0, a] | a ∈ $Z^+ \cup \{0\}$} ⊆ S; P is an interval semifield, hence S is a S-interval semiring.

It is interesting to note that all interval semirings need not in general be S- interval semirings.

***Example 2.41:*** Let S = {[0, a] | a ∈ $Z^+ \cup \{0\}$} be an interval semiring. Clearly S is not a S- interval semiring as S does not contain any proper interval subset which is an interval semifield.

We can as in case of S- semirings define S- interval subsemirings. We call a non empty subset T of an interval semiring A to be a S-interval subsemiring if T is an interval subsemiring of A and T has a proper subset B, (B ⊆ T) such that B is an interval semifield.

In view of this we have the following interesting result; the proof of which is left to the reader.

**THEOREM 2.7:** *Let A be an interval semiring. Suppose A has a proper subset B ⊆ A such that B is a S- interval subsemiring then A is a S- interval semiring.*

We will illustrate this situation by an example.

***Example 2.42:*** Let V = {[0, a] | a ∈ $R^+ \cup \{0\}$} be an interval semiring. Take P = {[0, a] | a ∈ $Q^+ \cup \{0\}$} ⊆ V; P is an S-interval subsemiring of V as P contains a proper subset T = {[0, a] | a ∈ $Z^+ \cup \{0\}$} ⊆ P such that T is an interval semifield.

Now V is a S-interval semiring as T ⊆ P ⊆ V and T is an interval semifield of V which is a proper subset of V.

Hence the claim.

It is interesting to note that every interval subsemiring of an S- interval semiring need not in general be an S- interval subsemiring.

We will illustrate this situation by some examples.

***Example 2.43:*** Let $S = \{[0, a] \mid a \in Q^+ \cup \{0\}\}$ be a S- interval semiring. Take $W = \{[0, a] \mid a \in 3Z^+ \cup \{0\}\} \subseteq S$, W is only an interval subsemiring and is not an S- interval subsemiring.

***Example 2.44:*** Let $T = \{[0, a] \mid a \in R^+ \cup \{0\}\}$ be an S- interval semiring. Take $V = \{[0, a] \mid a \in Z^+ \cup \{0\}\} \subseteq T$; V is only an interval subsemiring which is not an S- interval subsemiring.

We can define S- interval ideals of an interval semiring in the usual way [13].

***Example 2.45:*** Let $V = \{[0, a] \mid a \in R^+ \cup \{0\}\}$ be an interval semiring. Take $P = \{[0, a] \mid a \in Q^+ \cup \{0\}\} \subseteq V$ be a S-interval subsemiring. Consider $A = \{[0, a] \mid a \in Z^+ \cup \{0\}\} \subseteq P$. A is an interval semifield of P. We see for every $q \in P$ and $a \in A$ aq and qa are in A. Hence P is a S-interval ideal of V.

We wish to state here that it is an interesting problem to study whether the interval semiring constructed using $Z^+ \cup \{0\}$ or $Q^+ \cup \{0\}$ or $R^+ \cup \{0\}$ have no S-ideals. However if we take direct product we have S-interval ideals.

***Example 2.46:*** Let $V = \{[0, a] \times [0, b] \mid a, b \in Q^+ \cup \{0\}\}$ be an interval semiring.

Take $P = \{[0, a] \times [0\ b] \mid a, b \in Z^+ \cup \{0\}\} \subseteq V$ be an S- interval subsemiring. For $A = \{[0, a] \times \{0\} \mid a \in Z^+ \cup \{0\}\} \subseteq P$ is an interval semifield contained in P. We see P is an S- interval ideal of V.

***Example 2.47:*** Let $V = \{([0, a], [0, b], [0, c], [0, d]) \mid a, b, c, d$ belongs to $Q^+ \cup \{0\}\}$ be an interval semiring. Clearly V is a S- interval semiring. Consider $P = \{([0, a], [0, b], \{0\}, \{0\}) \mid a, b \in Q^+ \cup \{0\}\} \subseteq V$ be a S- interval subsemiring of V. P is an S-interval ideal as $A = \{([0, a], 0, 0, 0) \mid a \in Q^+ \cup \{0\}\} \subseteq P$ is such that for every $p \in P$ and $a \in A$ pa and ap are in P.

We can define the notion of Smarandache pseudo interval subsemiring (S- pseudo interval subsemiring) of an interval semiring in an analogous way [13].

We will only illustrate this situation by an example.

***Example 2.48:*** Let S = {[0, a] / a ∈ $Q^+ \cup \{0\}$} be an interval semiring. Consider A = {[0, a] / a ∈ $3Z^+ \cup \{0\}$} ⊆ S be a subset of S. Consider P = { [0, a] / a ∈ $Z^+ \cup \{0\}$} ⊆ S, we see A ⊆ P and P is a S-interval subsemiring of S, so A is a Smarandache pseudo interval subsemiring (S- pseudo interval subsemiring) of S. Infact S has infinitely many S- pseudo interval subsemirings.

The interesting factor is we have interval semirings, V which has no S- pseudo subsemiring.

This is illustrated by the following example.

***Example 2.49:*** Let V = {[0, a] / a ∈ $Z^+ \cup \{0\}$} be an interval semiring. V has no S- pseudo interval subsemirings.

This can be easily verified by the reader.

We can define Smarandache pseudo right (left) ideal (S-pseudo (left) right ideal) of the interval semiring as incase of semirings [13].

We will only illustrate this situation by some examples.

***Example 2.50:*** Let V = {[0, a] / a ∈ $Q^+ \cup \{0\}$} be an interval semiring. P = {[0, a] | a ∈ $3Z^+ \cup \{0\}$} is a S- pseudo interval subsemiring of V. Consider A = {[0, a] | a ∈ $Z^+ \cup \{0\}$} ⊆ V, clearly P ⊆ A and pa and ap are in P for every p ⊆ P and a ⊆ A. Thus P is a Smarandache pseudo interval ideal of V.

Likewise we can define Smarandache dual interval ideal and Smarandache pseudo dual ideal of an interval semiring.

The reader is left with the task of defining them and giving examples of them. Now we can define Smarandache special elements of a Smarandache interval semiring.

**DEFINITION 2.7:** *Let S be any interval semiring. We say a, b ∈ S is a Smarandache interval zero divisor (S-interval zero divisor) if a.b = 0 and there exists x, y ∈ S \ {a, b, 0}, x ≠ y with*
  1. *ax = 0 or xa = 0*
  2. *by = 0 or yb = 0*
  3. *xy ≠ 0 or yx ≠ 0.*

We will illustrate this situation by an example.

***Example 2.51:*** Let $S = \{([0, a_1], [0, a_2], [0, a_3], [0, a_4], [0, a_5], [0, a_6]) \mid a_i \in Q^+ \cup \{0\}; 1 \leq i \leq 6\}$ be an interval semiring. Take $a = ([0, a_1], [0, a_2], 0, 0, 0, 0)$ and $b = (0, 0, [0, b_1] [0, b_2], 0, 0)$ in S. Clearly $a.b = (0, 0, 0, 0, 0, 0)$.

Now consider $x = (0, 0, 0, 0, [0, x_1], [0, x_2])$ and $y = ([0, y_1] [0, y_2], 0, 0, 0 [0, y_3]) \in S$, we see $ax = 0$, $by = 0$ and $xy \neq 0$. Thus a, b in S is a Smarandache interval zero divisor of an interval semiring S.

**DEFINITION 2.8:** *Let S be an interval semiring. An element $x \in S$ is said to be a Smarandache anti – interval zero divisor (S-anti interval zero divisor) if we can find a y such that $xy \neq 0$ and $a, b \in S \setminus \{0, x, y\}$ such that*
  i.   *$ax \neq 0$ or $xa \neq 0$*
  ii.  *$by \neq 0$ or $yb \neq 0$*
  iii. *$ab = 0$ or $ba = 0$.*

We will illustrate this situation by an example.

***Example 2.52:*** Let $S = \{[0, a_1], [0, a_2] [0, a_3], \ldots, [0, a_9]) / a_i \in R^+ \cup \{0\}; 1 \leq i \leq 9\}$ be an interval semiring. Let $x = ([0, a_1], [0, a_2] [0, a_3], 0, 0, 0, 0, 0, 0)$ and $y = (0, [0, a] [0, b] [0, c], 0, 0, 0, 0, 0)$ in S.

We see $x.y \neq (0, 0, 0, \ldots, 0)$ that is $x.y = (0, [0, a_2 a], [0, a_3 b], 0, 0, 0, 0, 0, 0) \neq (0, 0, \ldots, 0)$. Consider $a = ([0, x_1], 0, 0, \ldots, 0)$ and $b = (0, 0, 0, [0, y_1], 0, \ldots, 0)$ in S. We see

$$ax = ([0, a, x_1], 0, 0, \ldots, 0) \neq (0, 0, \ldots, 0)$$
and
$$by = (0, 0, 0, [0, cy_1], 0, \ldots, 0) \neq (0, 0, 0, \ldots, 0);$$
but
$$ab = (0, 0, \ldots, 0).$$

Thus $x, y \in S$ is a S- interval anti zero divisor of the interval semiring S.

However it is important and interesting to note that an interval semifield will not have S- interval zero divisors or S-interval anti zero divisors.

Further if x ∈ S, S an interval semiring and x is an S anti interval zero divisor then x need not be an interval zero divisor.

We now proceed on to define Smarandache interval idempotent of an interval semiring.

**DEFINITION 2.9:** *Let S be an interval semiring. An element $0 \neq [0, a] \in S$ is a Smarandache idempotent (S- idempotent) of S if*

*(i)* $a^2 = a$ *that is* $([0, a])^2 = [0, a^2] = [0, a]$
*(ii)* *There exists $[0, b] \in S \setminus [0, a]$ such that*
   *a.* $([0, b])^2 = [0, a]$
   *b.* $[0, a] [0, b] = [0, b]$ ($[0, b] [0, a] = [0, b]$) or $[0, b] [0, a] = [0, a]$ or ($[0, a] [0, b] = [0, a]$)*

*'or' used in condition (ii) is mutually exclusive [15-17].*

We see all the interval semirings built using $Z^+ \cup \{0\}$ or $Q^+ \cup \{0\}$ or $R^+ \cup \{0\}$ have no S- idempotents. On similar lines we can define Smarandache units in an interval semiring.

**DEFINITION 2.10:** *Let S be an interval semiring with identity $[0, 1]$. We say $x \in S \setminus \{[0, 1]\}$ to be a Smarandache unit (S-unit) if there exists a $[0, y] \in S$ such that*
*(i)* $[0, x] \cdot [0, y] = [0, 1]$
*(ii)* *There exists $[0, a], [0, b] \in S \setminus \{[0, x], [0, y], [0, 1]\}$ such that*
   *a.* $[0, x] [0, a] = [0, y]$ *or* $[0, a] [0, x] = [0, y]$
   *b.* $[0, y] [0, b] = [0, x]$ *or* $[0, b] [0, y] = [0, x]$ *and* $[0, a] [0, b] = [0, 1]$.*

For more refer [13, 15].

The reader is requested to construct examples of S-units in interval semirings.

**Chapter Three**

# NEW CLASSES OF INTERVAL SEMIRINGS

In this chapter we introduce new classes of interval semirings and describe them. These new classes of interval semirings satisfy several interesting properties which are satisfied by the interval semirings introduced in Chapter two of this book. This chapter has two sections. Section one introduces matrix interval semirings. Interval polynomial semirings are introduced in section two.

## 3.1 Matrix Interval Semirings

In this section we introduce the new classes of matrix interval semirings and study some of the properties enjoyed by them.

**DEFINITION 3.1.1:** *$S = \{([0, a_1], …, [0, a_n]) / a_i \in Q^+ \cup \{0\}$ or $R^+ \cup \{0\}$ or $Z^+ \cup \{0\}\}$ be a row interval matrix. Define usual componentwise addition and componentwise multiplication on S. S is an interval semiring known as the row matrix interval semiring.*

We illustrate this situation by some examples.

*Example 3.1.1:* Let $V = \{([0, a_1], [0, a_2], [0, a_3], [0, a_4]) \mid a_i \in Z^+ \cup \{0\}; 1 \leq i \leq 4\}$; V is a row interval matrix semiring.

*Example 3.1.2:* Let $W = \{([0, a_1], [0, a_2], \ldots, [0, a_7]) \mid a_i \in Q^+ \cup \{0\}; 1 \leq i \leq 7\}$ be a row interval matrix semiring.

*Example 3.1.3:* Let $T = \{([0, a_1], [0, a_2], [0, a_3]) \mid a_i \in R^+ \cup \{0\}; 1 \leq i \leq 3\}$ be a row interval matrix semiring.

We can define row interval matrix subsemiring and row interval matrix ideal.
We will illustrate these substructures by examples.

*Example 3.1.4:* Let $V = \{([0, a_1], [0, a_2], \ldots, [0, a_{12}]) \mid a_i \in Q^+ \cup \{0\}; 1 \leq i \leq 12\}$ be an interval row matrix semiring. $W = \{([0, a_1], [0, a_2], \ldots, [0, a_{12}]) / a_i \in 3Z^+ \cup \{0\}; 1 \leq i \leq 12\} \subseteq V$; W is easily verified to be an interval row matrix subsemiring. Clearly W is not a interval row matrix ideal of V.

*Example 3.1.5:* Let $V = \{([0, a_1], [0, a_2], \ldots, [0, a_6]) / a_i \in Z^+ \cup \{0\}; 1 \leq i \leq 6\}$ be an interval row matrix semiring. Choose $T = \{([0, a_1], [0, a_2], \ldots, [0, a_6]) \mid a_i \in 3Z^+ \cup \{0\}; 1 \leq i \leq 6\} \subseteq T$; T is a interval row matrix subsemiring as well as interval row matrix ideal of V.

*Example 3.1.6:* Let $V = \{([0, a_1], [0, a_2], \ldots, [0, a_{10}]) / a_i \in R^+ \cup \{0\}; 1 \leq i \leq 10\}$ be an interval matrix semiring. Choose $W = \{([0, a_1], [0, a_2], \ldots, [0, a_{10}]) \mid a_i \in Z^+ \cup \{0\}; 1 \leq i \leq 10\} \subseteq V$; W is an interval row matrix subsemiring which is not an interval matrix row ideal of V.

*Example 3.1.7:* Let $V = \{([0, a_1], [0, a_2], \ldots, [0, a_7]) \mid a_i \in Q^+ \cup \{0\}; 1 \leq i \leq 7\}$ be an interval row matrix semiring. Choose $W = \{([0, a_1], [0, a_2], 0, 0, \ldots, 0, [0, a_7]) \mid a_1, a_2 \text{ and } a_7 \in Q^+ \cup \{0\}\} \subseteq V$; W is an interval row matrix semiring as well as W is an interval row matrix ideal of V.

***Example 3.1.8:*** Let $V = \{([0, a_1], [0, a_2], \ldots, [0, a_{25}]) \mid a_i \in R^+ \cup \{0\}; 1 \leq i \leq 25\}$ be an interval row matrix semiring. $W = \{([0, a_1], 0, 0, \ldots, 0) \mid a_1 \in Z^+ \cup \{0\}\} \subseteq V$; W is an interval row matrix subsemigroup and W is not an interval row matrix ideal of V.

***Example 3.1.9:*** Let $V = \{([0, a_1], [0, a_2], \ldots, [0, a_{19}]) \mid a_i \in R^+ \cup \{0\}; 1 \leq i \leq 19\}$ be an interval row matrix semiring. Choose $W = \{([0, a_1], \ldots, 0, 0, 0, [0, a_9]) \mid a_9, a_1 \in Q^+ \cup \{0\}\} \subseteq V$; W is an interval row matrix subsemiring of V. However W is not an interval row matrix ideal of V.

***Example 3.1.10:*** Let $V = \{([0, a_1], [0, a_2], [0, a_3]) \mid a_i \in Z^+ \cup \{0\}; 1 \leq i \leq 3\}$ be an interval row matrix semiring. Choose $W = \{([0, a_1], [0, a_2], [0, a_3]) \mid a_i \in 5Z^+ \cup \{0\}; 1 \leq i \leq 3\}$ a proper subset of V. W is an interval row matrix subsemiring of V as well as an interval row matrix ideal of V.

Infact V has infinitely many interval ideals.

We see all these row matrix interval semirings are of infinite cardinality further all them are interval row matrix commutative semirings.

Now we call an interval row matrix semiring S to be a Smarandache row matrix interval semiring of S if S has a proper subset $P \subseteq S$ such that P is a row matrix interval semifield of S.

We will illustrate this situation by some examples.

***Example 3.1.11:*** Let $T = \{([0, a_1], [0, a_2], \ldots, [0, a_{12}]) \mid a_i \in Q^+ \cup \{0\}; 1 \leq i \leq 12\}$ be a row matrix interval semiring. Consider $G = \{([0, a_1], 0, \ldots, 0) / a_1 \in Q^+ \cup \{0\} \subseteq T$; G is a row matrix interval semifield contained in T. Clearly T is a S-row matrix interval semiring.

It is interesting and important to note that all row matrix interval semirings in general need not be a S-row matrix interval semirings.

We will illustrate this situation by some examples.

***Example 3.1.12:*** Let $V = \{([0, a_1], [0, a_2], \ldots, [0, a_{32}]) \mid a_i \in 5Z^+ \cup \{0\}\}$ be a row matrix interval semiring. Clearly V does not contain a proper subset A such that A is a row matrix interval semifield. So V is not a S-row matrix interval semiring.

Infact we have an infinite class of row matrix interval semirings which are not S-row matrix interval semirings. We see it is not possible to have a $1 \times n$ row interval matrix semifield where $n > 1$; with interval entries different from zero. The reader is expected to prove this.

We can define Smarandache pseudo row matrix interval subsemiring, S-ideal, S-subsemiring, S-pseudo ideal and S-pseudo dual ideal in case of row matrix interval semirings in an analogous way.

We reserve this work to the reader, however give examples of these concepts.

***Example 3.1.13:*** Let $V = \{([0, a_1], [0, a_2], [0, a_3], [0, a_4]) / a_i \in Q^+ \cup \{0\}; 1 \leq i \leq 4\}$ be a row interval matrix semiring. Take $P = \{([0, a_1], 0, 0, 0) \mid a_1 \in Q^+ \cup \{0\}\} \subseteq V$; P is a row interval matrix semifield; hence V is a S-row interval matrix semiring. Take $W = \{([0, a_1], 0, [0, a_2], 0) \mid a_1, a_2 \in Q^+ \cup \{0\}\} \subseteq V$; W is a row interval matrix subsemiring of V. Take $G = \{(0, 0, [0, a], 0) / a \in Q^+ \cup \{0\}\} \subseteq W$; G is a row interval matrix semifield of W. Hence W is a S-row matrix interval semiring.

In view of this we have the following theorem the proof of which is simple.

**THEOREM 3.1.1:** *Let $V = \{([0, a_1], [0, a_2], \ldots, [0, a_n]) / a_i \in Q^+ \cup \{0\}$ or $Z^+ \cup \{0\}$ or $R^+ \cup \{0\}; 1 \leq i \leq n\}$ be a row matrix interval semiring. If V has a nontrivial S-row matrix interval subsemiring, then V is a S-row matrix interval semiring. However if V is a S-row matrix interval semiring then every row matrix interval subsemiring need not be a S-row matrix interval subsemiring.*

We will substantiate this by an example.

***Example 3.1.14:*** Let $V = \{([0, a_1], [0, a_2], [0, a_3], [0, a_4], [0, a_5], [0, a_6]) \mid a_i \in Z^+ \cup \{0\}; 1 \le i \le 6\}$ be a S-row matrix interval semiring. Consider $W = \{([0, a_1], [0, a_2], \ldots, [0, a_6]) \mid a_i \in 15Z^+ \cup \{0\}; 1 \le i \le 6\} \subseteq V$; Clearly W is only a row matrix interval subsemiring of V and W is not a S-row matrix interval subsemiring of V.

We will now give some examples of S-ideal of a row matrix interval semiring.

***Example 3.1.15:*** Let $V = \{([0, a_1], [0, a_2], [0, a_3], [0, a_4], [0, a_5]) \mid a_i \in Q^+ \cup \{0\}; 1 \le i \le 5\}$ be a row matrix interval semiring. Consider $W = \{(0, [0, a], [0, b], 0, 0) / a, b \in Q^+ \cup \{0\}\} \subseteq V$, be a S-row matrix interval subsemiring of V. Clealry $A = \{(0, [0, a], 0, 0, 0) / a \in Q^+ \cup \{0\}\} \subseteq W$ is a row matrix interval semifield such that for every $a \in A$ and $p \in W$, ap and pa are in A. Thus W is a S-row matrix interval ideal (row matrix interval S-ideal) of V.

***Example 3.1.16:*** Let $S = \{([0, a_1], [0, a_2], [0, a_3], \ldots, [0, a_9]) \mid a_i \in R^+ \cup \{0\}; 1 \le i \le 9\}$ be a row matrix interval semiring. Consider $P = \{([0, a_1], [0, a_2], [0, a_3], 0, 0, 0, 0, 0, 0) / a_i \in Q^+ \cup \{0\}; 1 \le i \le 3\}$ be a S-row matrix interval subsemiring of S. Take $A = \{(0, 0, [0, a_3], 0, 0, 0, 0, 0, 0) / a_3 \in Q^+ \cup \{0\}\} \subseteq P$; A is a row matrix interval semifield of S such that for all $a \in A$ and $p \in P$, ap and pa are in A. Thus P is a S-ideal of S.

Now we proceed onto give an example of a S-pseudo subsemiring.

***Example 3.1.17:*** Let $V = \{([0, a_1], [0, a_2], \ldots, [0, a_6]) \mid a_i \in Z^+ \cup \{0\}; 1 \le i \le 6\}$ be a row matrix interval semiring. Let $A = \{([0, a], 0, 0, 0, 0, 0) \mid a \in 5Z^+ \cup \{0\}\} \subseteq V$ be a non empty subset of V. Choose $P = \{([0, a], 0, [0, b], 0, 0, 0) / a, b \in Z^+ \cup \{0\}\} \subseteq V$, P is a S-row matrix interval subsemiring and $B = \{([0, a], 0, 0, 0, 0, 0) \mid a \in Z^+ \cup \{0\}\} \subseteq P$ is a row matrix interval semifield. Thus $A \subset P$ so A is a Smarandache pseudo row matrix interval subsemiring of V.

Interested reader can give examples of this structure.

Now we proceed onto give example of Smarandache pseudo ideals of a row matrix interval semiring.

***Example 3.1.18:*** Let $V = \{([0, a_1], [0, a_2], \ldots, [0, a_{12}]) / a_i \in Q^+ \cup \{0\}; 1 \leq i \leq 12\}$ be row matrix interval semiring.

Consider $P = \{([0, a_1], [0, a_2], [0, a_3], 0, 0, \ldots, 0) / a_1, a_2, a_3 \in Z^+ \cup \{0\}\}$ be a S-row matrix interval subsemiring of V. $A = \{(0, 0, [0, a], 0, 0, \ldots, 0) / a \in Z^+ \cup \{0\}\} \subseteq P$, is a row matrix interval semifield in P.

For every $p \in P$ and for every $a \in A$ we have pa, ap $\in P$.

Thus P is a S-pseudo row matrix interval subsemiring of V.

Now we proceed onto give examples of Smarandache special elements in a row matrix interval semirings.

***Example 3.1.19:*** Let $V = \{([0, a_1], [0, a_2], \ldots, [0, a_{10}]) / a_i \in Q^+ \cup \{0\}; 1 \leq i \leq 10\}$ be a row matrix interval semiring. Let $A = \{([0, a_1], [0, a_2], 0, \ldots, 0)$ and $b = \{(0, 0, 0, [0, a_1], [0, a_2], [0, a_3], , 0, 0, 0, 0)$ be in V. We see $a \cdot b = (0, 0, 0, \ldots, 0)$.

Now take $x = (0, 0, [0, a], 0, \ldots, 0)$ and $y = (0, [0, a], [0, b], 0, 0, \ldots, 0)$ in V. We see $ax = (0, 0, \ldots, 0)$ and $by = (0, 0, \ldots, 0)$, further $xy = (0, 0, [0, a], 0, \ldots, 0) \times (0, [0, a], [0, b], 0, \ldots, 0)$
$= (0, 0, [0, ab], 0, \ldots, 0)$
$\neq (0, 0, 0, \ldots, 0)$.

Thus a, b $\in$ V is a S-zero divisor of the row matrix interval semiring. Interested reader can construct more examples of S-zero divisors.

Now we proceed onto give an example of a S-antizero divisors in row matrix interval semirings.

***Example 3.1.20:*** Let $V = \{([0, a_1], [0, a_2], \ldots, [0, a_{15}]) \mid a_i \in Z^+ \cup \{0\}; 1 \leq i \leq 15\}$ be a row matrix interval semiring. Let $x = \{([0, a_1], [0, a_2], [0, a_3], [0, a_4], \ldots, [0, a_{15}])$ and $y = ([0, a_1], [0, a_2], 0, 0, [0, a_3], \ldots, [0, a'_{15}])$ be in V. Clearly $xy = ([0, a_1.a_2], 0, 0, [0, a_4 . a_3], 0, \ldots, 0, [0, a_{15} . a'_{15}]) \neq (0, 0, 0, \ldots, 0)$ in S.

Let $a = ([0, a_1^1], 0, \ldots, 0)$ and $b = (0, 0, 0, [0, a], 0, \ldots, 0)$ in S. We see $ax = ([0, a_1 \cdot a_1^1], 0, \ldots, 0) \neq (0, 0, \ldots, 0)$ and $by = (0, 0, 0, [0, a_3 a], 0, 0, \ldots, 0) \neq (0, 0, \ldots, 0)$.

But $ab = (0, 0, \ldots, 0)$.

Thus $x, y \in V$ is a S-anti zero divisor of the row matrix interval semiring V.

We see in row matrix interval semiring a S-anti zero divisor in general is not a zero divisor. However if a row interval matrix semiring S has S-anti zero divisors then S has zero divisors.

The reader is given the task of defining S-units in a row matrix interval semiring and illustrate it with examples.

Now having studied some properties about row matrix interval semirings, we now proceed onto define a square matrix interval semirings. Clearly we cannot define $m \times n$ ($m \neq n$) matrix interval semirings where $m \neq 1$.

**DEFINITION 3.1.2:** *Let V = {Set of all $n \times n$ interval matrices with intervals of the form [0, a] where $a \in Z^+ \cup \{0\}$ or $Q^+ \cup \{0\}$ or $R^+ \cup \{0\}$}; V is a semiring under interval matrix addition and interval matrix multiplication. We define V to be a $n \times n$ square matrix interval semiring.*

We will illustrate this situation by some examples.

*Example 3.1.21:* Let
$$V = \left\{ \begin{bmatrix} [0, a_1] & [0, a_2] \\ [0, a_3] & [0, a_4] \end{bmatrix} \middle| a_i \in Q^+ \cup \{0\}; 1 \leq i \leq 4 \right\}$$
be a $2 \times 2$ square matrix interval semiring.

*Example 3.1.22:* Let S = {all $5 \times 5$ interval matrices with intervals of the form $[0, a_i]$ where $a_i \in Z^+ \cup \{0\}$} be a $5 \times 5$ square matrix interval semiring.

*Example 3.1.23:* Let T = {All $10 \times 10$ upper triangular interval matrices with intervals of the form [0, a] with $a \in R^+ \cup \{0\}$} be a $10 \times 10$ square matrix interval semiring.

**Example 3.1.24:** Let G = {All 8 × 8 lower triangular interval matrices with intervals of the form [0, $x_i$] where $x_i \in Q^+ \cup \{0\}$} be a 8 × 8 square matrix interval semiring.

Now as in case of semirings we can define the basic substructure in these semirings which is left as an exercise to the reader.

However we will illustrate this situation by some examples.

**Example 3.1.25:** Let P = {All 3 × 3 interval matrices with intervals of the form [0, a] where $a \in Z^+ \cup \{0\}$} be a 3 × 3 square matrix interval semiring. Choose G = {All 3 × 3 interval matrices with intervals of the form [0, a] where $a \in 5 Z^+ \cup \{0\}$} $\subseteq$ P; G is a square matrix interval subsemiring of P.

**Example 3.1.26:** Let

$$V = \left\{ \begin{pmatrix} [0,a] & [0,b] \\ [0,c] & [0,d] \end{pmatrix} \,\middle|\, a,b,c,d \in Z^+ \cup \{0\} \right\}$$

be a 2 × 2 interval matrix semiring.

$$T = \left\{ \begin{pmatrix} [0,a] & [0,a] \\ [0,a] & [0,a] \end{pmatrix} \,\middle|\, a \in Z^+ \cup \{0\} \right\} \subseteq V,$$

T is a 2 × 2 interval matrix subsemiring of V.

**Example 3.1.27:** Let V = {collection of all 6 × 6 interval matrices with intervals of the form [0, a], $a \in Q^+ \cup \{0\}$} be a square interval matrix semiring. Consider P = {Collection of all 6 × 6 upper triangular interval matrices with intervals of the form [0, a] where $a \in Q^+ \cup \{0\}$} $\subseteq$ V, P is an interval square matrix subsemiring of V.

*Example 3.1.28:* Let

$$T = \left\{ \begin{pmatrix} [0,a] & [0,a] \\ [0,a] & [0,a] \end{pmatrix} \middle| a,b,c,d \in Z^+ \cup \{0\} \right\}$$

be a square interval semigroup.

$$P = \left\{ \begin{pmatrix} [0,a] & [0,a] \\ [0,a] & [0,a] \end{pmatrix} \middle| a \in 15Z^+ \cup \{0\} \right\} \subseteq T;$$

P is a square interval subsemiring of T.

We now proceed onto give examples of square interval matrix ideal of the semiring.

*Example 3.1.29:* Let V = {all 10 × 10 square interval matrices with intervals of the form [0, a] with a ∈ $Z^+ \cup \{0\}$} be a matrix interval semiring. Take I = {all 10 × 10 square interval matrices with intervals of the form [0, a] where a ∈ $5Z^+ \cup \{0\} \subseteq$ V; I is an ideal of V.

*Example 3.1.30:* Let V = {all 7 × 7 square interval matrices with intervals of the form [0, a]; a ∈ $Q^+ \cup \{0\}$} be a square matrix interval semiring. W = {all 7 × 7 diagonal interval matrices with intervals of the form [0, a] where a ∈ $Q^+ \cup \{0\}$} $\subseteq$ V; W is a square matrix interval ideal of V.

*Example 3.1.31:* Let V = {5 × 5 interval matrices with intervals of the form [0, a] where a ∈ $Q^+ \cup \{0\}$} be an interval matrix semiring. Take S = {5 × 5 interval matries with intervals of the form a ∈ $Z^+ \cup \{0\}$} $\subseteq$ V is an interval matrix subsemiring. Clearly S is not an interval matrix ideal of V.

*Example 3.1.32:* Let V = {all 3 × 3 interval matrices with intervals of the form [0, a]; a ∈ $5Z^+ \cup \{0\}$} be a square matrix interval semiring. Consider W = {all 3 × 3 upper triangular interval matrices with intervals of the form [0, a] where a ∈ $15Z^+ \cup \{0\}$} ⊆ V, it is easily verified W is only an interval square matrix subsemiring and is not an ideal of V.

Now having seen ideals and subsemirings of square intervals matrix semirings now we proceed onto define Smarandache square matrix interval semiring.

**DEFINITION 3.1.3**: *Let V = {set of all n × n square interval matrices with intervals of the form [0, a], a ∈ $Q^+ \cup \{0\}$ or $R^+ \cup \{0\}$} be a square interval matrix semiring.*

*Let P ⊆ V, if P is a square interval matrix semifield, we define V to be a Smarandache square matrix interval semiring (S-square matrix interval semiring).*

*Example 3.1.33:* Let V = {2 × 2 interval matrices with intervals of the form [0, a]; a ∈ $Q^+ \cup \{0\}$} be a square interval matrix semiring. Take

$$W = \left\{ \begin{bmatrix} a & 0 \\ 0 & 0 \end{bmatrix} \middle| a \in Q^+ \cup \{0\} \right\} \subseteq V,$$

W is a square matrix interval semifield, hence W is a S-square matrix interval semifield.

*Example 3.1.34:* Let V = {All 10 × 10 triangular interval matrices with intervals of the form [0, a], a ∈ $Q^+ \cup \{0\}$} be a square matrix interval semiring.

Consider W = {All 10 × 10 upper diagonal interval matrices with elements of the form [0, a], a ∈ $Q^+ \cup \{0\}$; every element in the diagonal is the same} ⊆ V; W is a square matrix interval semifield, thus V is a S-interval matrix semiring.

*Example 3.1.35:* Let V = {all 4 × 4 interval matrices with intervals of the form [0, a] | a ∈ $Q^+ \cup \{0\}$} be a square matrix interval semiring.

Consider

$$W = \left\{ \begin{bmatrix} [0,a] & 0 & 0 & 0 \\ 0 & [0,a] & 0 & 0 \\ 0 & 0 & [0,a] & 0 \\ 0 & 0 & 0 & [0,a] \end{bmatrix} \middle| a \in Q^+ \cup \{0\} \right\} \subseteq V,$$

W is a square matrix interval semifield contained in V. Thus V is a S-square matrix interval semiring.

It is interesting to note that all square interval matrix semirings are not in general Smarandache square interval matrix semirings. We will illustrate this situation by some examples.

*Example 3.1.36:* Let V = {all 2 × 2 upper triangular interval matrices with intervals of the form [0, a]; a ∈ $3Z^+ \cup \{0\}$} be a matrix interval semiring. We see V has no proper subset T ⊆ V, such that T is a matrix interval semifield. Thus V is not a matrix interval semiring.

*Example 3.1.37:* Let

$$V = \left\{ \begin{bmatrix} [0,a] & 0 \\ 0 & 0 \end{bmatrix} \middle| a,b,c,d \in Z^+ \cup \{0\} \right\}$$

be a matrix interval semiring.

$$W = \left\{ \begin{bmatrix} [0,a] & 0 \\ 0 & 0 \end{bmatrix} \middle| a \in Z^+ \cup \{0\} \right\} \subseteq V;$$

W is a matrix interval semifield, thus V is a S-matrix interval semiring.

Now we will define S-subsemirings and S-ideals using interval matrix semirings.

**DEFINITION 3.1.4:** *Let S = {All n × n interval matrices with intervals of the form [0, a] / a ∈ $Z^+ \cup \{0\}$} or $Q^+ \cup \{0\}$ or $R^+ \cup \{0\}$} be an interval matrix semiring. Let P ⊆ S, be such that P is* *a interval matrix subsemiring of S. If P contains a subset T, $\phi \not\subseteq T \not\subset P$ such that T is a interval matrix semifield, then we define P to a Smarandache interval matrix subsemiring (S-interval matrix subsemiring).*

We will illustrate this situation by some examples.

***Example 3.1.38:*** Let S = {all 5 × 5 interval matrices with intervals of the form [0, a] where a ∈ $Q^+ \cup \{0\}$} be an interval matrix semiring. Consider P = {all 5 × 5 diagonal interval matrices with intervals of the form [0, a] where a ∈ $Q^+ \cup \{0\}$} $\subseteq$ S, P is a S-interval matrix subsemiring as P contains

$$T = \left\{ \begin{bmatrix} [0,a] & 0 & 0 & 0 & 0 \\ 0 & 0 & 0 & 0 & 0 \\ 0 & 0 & 0 & 0 & 0 \\ 0 & 0 & 0 & 0 & 0 \\ 0 & 0 & 0 & 0 & 0 \end{bmatrix} \,\middle|\, a \in Q^+ \cup \{0\} \right\} \subseteq P,$$

so that T is a semifield.

***Example 3.1.39:*** Let V = {All 6 × 6 upper triangular interval matrices with intervals of the form [0, a] where a ∈ $Z^+ \cup \{0\}$} be the interval matrix semiring. Take W = {all 6 × 6 diagonal interval matrices with intervals of the form [0, a] where a ∈ $Z^+ \cup \{0\}$} $\subseteq$ V; W is an interval matrix subsemiring of V.
Now

$$T = \left\{ \begin{bmatrix} 0 & 0 & 0 & 0 & 0 & 0 \\ 0 & [0,a] & 0 & 0 & 0 & 0 \\ 0 & 0 & 0 & 0 & 0 & 0 \\ 0 & 0 & 0 & 0 & 0 & 0 \\ 0 & 0 & 0 & 0 & 0 & 0 \\ 0 & 0 & 0 & 0 & 0 & 0 \end{bmatrix} \,\middle|\, a \in Z^+ \cup \{0\} \right\} \subseteq W;$$

is an interval matrix semifield. Thus W is a S-matrix interval subsemiring of V.

It is left as an exercise for the reader to prove the following theorem.

**THEOREM 3.1.2:** *Let V be an interval matrix semiring. If V has a S-interval matrix subsemiring then V itself is a S-interval matrix semiring.*
*However if V is a S-interval matrix semiring then in general all interval matrix subsemirings in V need not be a S-matrix interval subsemirings.*

We prove this by the following examples.

*Example 3.1.40:* Let V = {All 3 × 3 interval matrices with intervals of the form [0, a] with $a \in Z^+ \cup \{0\}$} be a interval matrix semiring.
Choose

$$W = \left\{ \begin{bmatrix} 0 & 0 & 0 \\ 0 & 0 & 0 \\ 0 & 0 & [0,a] \end{bmatrix} \middle| a \in Z^+ \cup \{0\} \right\} \subseteq V$$

be an interval matrix semifield. Thus V is a S-interval matrix semiring. Now consider

$$T = \left\{ \begin{bmatrix} [0,a] & [0,b] & [0,c] \\ 0 & [0,d] & [0,e] \\ 0 & 0 & [0,f] \end{bmatrix} \middle| a,b,c,d,e,f \in 7Z^+ \cup \{0\} \right\} \subseteq V;$$

T is a interval matrix subsemiring of V. Clearly T has no proper interval matrix semifield, so T is not a S-interval matrix subsemiring of V.

However V is a S-interval matrix semiring. Hence the claim.

Now we proceed onto define S-ideals of a interval matrix semiring.

**DEFINITION 3.1.5:** *Let V be a S-matrix interval semiring. A non empty subset P of V is said to be a Smarandache interval matrix ideal (S-interval matrix ideal) of V if the following conditions are satisfied.*
1. *P is an S-interval matrix subsemiring.*
2. *For every $p \in P$ and $A \subseteq P$ where A is the semifield of P we have for all $a \in A$ and $p \in P$, pa and ap is in A.*

We will illustrate this situation by some examples.

*Example 3.1.41:* Let V = {all 4 × 4 interval matrices with intervals of the form [0, a] where $a \in Q^+ \cup \{0\}$} be an interval matrix semiring.
Choose

$$W = \left\{ \begin{bmatrix} [0,a] & 0 & 0 & 0 \\ 0 & [0,b] & 0 & 0 \\ 0 & 0 & [0,c] & 0 \\ 0 & 0 & 0 & 0 \end{bmatrix} \middle| a,b,c \in Q^+ \cup \{0\} \right\} \subseteq V.$$

W is a S-matrix interval ideal of V.

*Example 3.1.42:* Let V = {all 3 × 3 interval matrices with intervals of the form [0, a] where $a \in Z^+ \cup \{0\}$} be a matrix interval semiring. Choose

$$W = \left\{ \begin{bmatrix} [0,a] & 0 & 0 \\ 0 & [0,b] & 0 \\ 0 & 0 & [0,c] \end{bmatrix} \middle| a,b,c \in Z^+ \cup \{0\} \right\} \subseteq V$$

be a S-ideal of V.

The following result is left as an exercise for the reader to prove.

**THEOREM 3.1.3:** *Let V be a S-matrix interval semiring. Every S-matrix interval ideal of V is a S-matrix interval subsemiring of V but every S-matrix interval subsemiring of V in general need not be a S-ideal of V.*

We now proceed onto define the notion of S-pseudo interval matrix subsemiring.

**DEFINITION 3.1.6:** *Let S be a interval matrix semiring. A non empty proper subset A of S is said to be Smarandache pseudo interval matrix subsemiring (S-pseudo interval matrix subsemiring) if the following condition is true.*

*If there exists a subset P of S such that $A \subseteq P$, where P is a S-interval matrix subsemiring, that is P has a subset B such that B is a semifield under the operations of S or P itself is a semifield under the operations of S.*

We will illustrate this situation by some examples.

*Example 3.1.43:* Let S = {All 4 × 4 interval matrices with intervals of the form [0, a] where $a \in Q^+ \cup \{0\}$} be a matrix interval semiring.
Let
$$A = \left\{ \begin{bmatrix} [0,a] & [0,b] & 0 & 0 \\ 0 & [0,c] & 0 & [0,b] \\ 0 & 0 & 0 & [0,d] \\ 0 & 0 & 0 & 0 \end{bmatrix} \middle| a,b,c,d \in Q^+ \cup \{0\} \right\} \subseteq S$$
be a proper subset of S. Now

$$P = \left\{ \begin{bmatrix} [0,a] & [0,b] & [0,c] & [0,d] \\ 0 & [0,e] & [0,f] & [0,h] \\ 0 & 0 & [0,m] & [0,n] \\ 0 & 0 & 0 & [0,t] \end{bmatrix} \middle| \begin{array}{l} a,b,c,d,e,f,h,m, \\ n,t \in Q^+ \cup \{0\} \end{array} \right\} \subseteq S$$

is a S-matrix interval subsemiring of S. Thus A is a S-pseudo interval matrix subsemiring of S.

*Example 3.1.44:* Let V = {set of all 5 × 5 interval matrices with intervals of the form [0, a] where a ∈ $Z^+ \cup \{0\}$} be an interval matrix semiring.
  Let

$$A = \left\{ \begin{bmatrix} [0,a] & 0 & [0,b] & 0 & [0,c] \\ 0 & [0,a] & 0 & [0,b] & 0 \\ 0 & 0 & [0,b] & 0 & [0,a] \\ 0 & 0 & 0 & [0,a] & 0 \\ 0 & 0 & 0 & 0 & [0,b] \end{bmatrix} \middle| a,b,c \in Z^+ \cup \{0\} \right\}$$

$\subseteq$ V be a subset of interval 5 × 5 matrices of V. Now A $\subseteq$ {all collection of 5 × 5 upper triangular interval matrices with entries of the form [0, a], a ∈ $Z^+ \cup \{0\}$} = P $\subseteq$ V, P is a S-matrix interval subsemiring of V. Now

$$T = \left\{ \begin{bmatrix} 0 & 0 & 0 & 0 & 0 \\ 0 & [0,a] & 0 & 0 & 0 \\ 0 & 0 & 0 & 0 & 0 \\ 0 & 0 & 0 & 0 & 0 \\ 0 & 0 & 0 & 0 & 0 \end{bmatrix} \middle| a \in Z^+ \cup \{0\} \right\} \subseteq P$$

is a matrix interval semifield in V. Thus A is a S-pseudo interval matrix subsemiring of V.

Now we proceed onto define the notion of Smarandache pseudo matrix interval ideal of an interval matrix semiring.

**DEFINITION 3.1.7:** *Let V = {collection of all n × n interval matrices with intervals of the form [0, a] where a ∈ $Z^+$ ∪ {0} or $Q^+$ ∪ {0} or $R^+$ ∪ {0}} be an interval matrix semiring. Let P ⊆ V be a S-pseudo subsemiring of V and A ⊆ P, A a matrix interval semifield or a S-interval matrix subsemiring in V. If for every p ∈ P and for every a ∈ A ap and pa ∈ P then we define P to be a Smarandache pseudo matrix interval ideal of V. (S-pseudo matrix interval ideal of V).*

We will illustrate this situation by an example.

*Example 3.1.45:* Let V = {4 × 4 interval matrices with intervals of the form [0, a] with a ∈ $R^+$ ∪ {0}} be an interval matrix semiring. Choose

$$P = \left\{ \begin{bmatrix} [0,a] & 0 & [0,a] & 0 \\ 0 & [0,b] & 0 & [0,b] \\ 0 & 0 & [0,c] & 0 \\ 0 & 0 & 0 & [0,d] \end{bmatrix} \mid a, b, c, d \in R^+ \cup \{0\} \right\} \subseteq V.$$

P is a S-pseudo matrix interval subsemiring of V as P ⊆

$$A = \left\{ \begin{bmatrix} [0,a_1] & [0,a_2] & [0,a_3] & [0,a_4] \\ 0 & [0,a_5] & [0,a_6] & [0,a_7] \\ 0 & 0 & [0,a_8] & [0,a_9] \\ 0 & 0 & 0 & [0,a_{10}] \end{bmatrix} \bigg| \begin{array}{l} a_i \in R^+ \cup \{0\}; \\ 1 \leq i \leq 10 \end{array} \right\} \subseteq V$$

and A is a S-interval matrix subsemiring as A contains

$$T = \left\{ \begin{bmatrix} [0,a] & 0 & 0 & 0 \\ 0 & 0 & 0 & 0 \\ 0 & 0 & 0 & 0 \\ 0 & 0 & 0 & 0 \end{bmatrix} \middle| a \in R^+ \cup \{0\} \right\} \subseteq A \subseteq V$$

is a interval matrix semifield in V. Thus P is a S-pseudo interval matrix ideal of V.

Now we can as in case of other interval semirings define Smarandache zero divisors and Smarandache anto zero divisors [ ]. The task of defining these notions are left as an exercise to the reader.

However we will illustrate this situation by an example.

*Example 3.1.46:* Let

$$V = \left\{ \begin{bmatrix} [0,a] & [0,c] \\ [0,b] & [0,d] \end{bmatrix} \middle| a,b,c,d \in Z^+ \cup \{0\} \right\}$$

be an interval matrix semiring. It is easily verified V has zero divisors and right or left zero divisor.
For instance

$$x = \begin{bmatrix} 0 & 0 \\ [0,a] & 0 \end{bmatrix} \text{ and } y = \begin{bmatrix} 0 & 0 \\ 0 & [0,b] \end{bmatrix}$$

in V is such that

$$x.y = \begin{bmatrix} 0 & 0 \\ [0,a] & 0 \end{bmatrix} \begin{bmatrix} 0 & 0 \\ 0 & [0,b] \end{bmatrix} = \begin{bmatrix} 0 & 0 \\ 0 & 0 \end{bmatrix}.$$

$$y.x = \begin{bmatrix} 0 & 0 \\ 0 & [0,b] \end{bmatrix} \begin{bmatrix} 0 & 0 \\ [0,a] & 0 \end{bmatrix} = \begin{bmatrix} 0 & 0 \\ [0,ab] & 0 \end{bmatrix} \neq \begin{bmatrix} 0 & 0 \\ 0 & 0 \end{bmatrix}.$$

Consider $p = \begin{bmatrix} [0,a] & 0 \\ 0 & 0 \end{bmatrix}$ and $q = \begin{bmatrix} 0 & 0 \\ 0 & [0,b] \end{bmatrix}$ in V

$$p \cdot q = \begin{bmatrix} [0,a] & 0 \\ 0 & 0 \end{bmatrix} \begin{bmatrix} 0 & 0 \\ 0 & [0,b] \end{bmatrix} = \begin{bmatrix} 0 & 0 \\ 0 & 0 \end{bmatrix}.$$

Now

$$q \cdot p = \begin{bmatrix} 0 & 0 \\ 0 & [0,b] \end{bmatrix} \begin{bmatrix} [0,a] & 0 \\ 0 & 0 \end{bmatrix} = \begin{bmatrix} 0 & 0 \\ 0 & 0 \end{bmatrix}.$$

Thus $p.q = (0)$ is a zero divisor in V.
Now choose

$$x = \begin{bmatrix} 0 & 0 \\ [0,t] & 0 \end{bmatrix} \in V$$

is such that

$$x.p = \begin{bmatrix} 0 & 0 \\ [0,t] & 0 \end{bmatrix} \begin{bmatrix} [0,a] & 0 \\ 0 & 0 \end{bmatrix}$$

$$= \begin{bmatrix} 0 & 0 \\ [0,ta] & 0 \end{bmatrix} \neq (0)$$

and

$$px = \begin{bmatrix} [0,a] & 0 \\ 0 & 0 \end{bmatrix} \begin{bmatrix} 0 & 0 \\ [0,t] & 0 \end{bmatrix}$$

$$= \begin{bmatrix} 0 & 0 \\ 0 & 0 \end{bmatrix}.$$

Now

$$y = \begin{bmatrix} 0 & [0,m] \\ 0 & 0 \end{bmatrix} \in V$$

is such that

$$qy = \begin{bmatrix} 0 & 0 \\ 0 & [0,b] \end{bmatrix} \begin{bmatrix} 0 & [0,m] \\ 0 & 0 \end{bmatrix} = \begin{bmatrix} 0 & 0 \\ 0 & 0 \end{bmatrix}$$

and
$$yq = \begin{bmatrix} 0 & [0,m] \\ 0 & 0 \end{bmatrix} \begin{bmatrix} 0 & 0 \\ 0 & [0,b] \end{bmatrix} = \begin{bmatrix} 0 & [0,mb] \\ 0 & 0 \end{bmatrix} \neq (0).$$

Now
$$xy = \begin{bmatrix} 0 & 0 \\ [0,t] & 0 \end{bmatrix} \begin{bmatrix} 0 & [0,m] \\ 0 & 0 \end{bmatrix} = \begin{bmatrix} 0 & 0 \\ 0 & [0,tm] \end{bmatrix} \neq (0).$$

$$yx = \begin{bmatrix} 0 & [0,m] \\ 0 & 0 \end{bmatrix} \begin{bmatrix} 0 & 0 \\ [0,t] & 0 \end{bmatrix} = \begin{bmatrix} [0,mt] & 0 \\ 0 & 0 \end{bmatrix} \neq (0).$$

Thus p, q ∈ V is a Smarandache zero divisor of V.

The following observation is important.
    If an interval matrix semiring V has Smarandache interval matrix zero divisors then V has zero divisors, but however the question is, will every zero divisor in an interval matrix semiring will always be a S-zero divisor; is open.

Now we proceed onto give an example of a Smarandache anti zero divisors of an interval matrix semiring V.

*Example 3.1.47:* Let
$$V = \left\{ \begin{bmatrix} [0,a] & [0,c] \\ [0,b] & [0,d] \end{bmatrix} \middle| a,b,c,d \in Q^+ \cup \{0\} \right\}$$
be an interval matrix semiring.
    Let
$$x = \begin{bmatrix} [0,a] & [0,b] \\ 0 & 0 \end{bmatrix} \text{ and } y = \begin{bmatrix} 0 & 0 \\ [0,t] & [0,u] \end{bmatrix}$$
be in V. We see
$$x.y = \begin{bmatrix} [0,a] & [0,b] \\ 0 & 0 \end{bmatrix} \begin{bmatrix} 0 & 0 \\ [0,t] & [0,u] \end{bmatrix}$$

$$= \begin{bmatrix} [0,bt] & [0,bu] \\ 0 & 0 \end{bmatrix} \neq (0).$$

$$y.x = \begin{bmatrix} 0 & 0 \\ [0,t] & [0,u] \end{bmatrix} \begin{bmatrix} [0,a] & [0,b] \\ 0 & 0 \end{bmatrix} = \begin{bmatrix} 0 & 0 \\ [0,at] & [0,tb] \end{bmatrix}.$$

Consider

$$p = \begin{bmatrix} 0 & 0 \\ 0 & [0,t] \end{bmatrix} \text{ and } q = \begin{bmatrix} [0,m] & 0 \\ 0 & 0 \end{bmatrix} \text{ in V.}$$

Now

$$p.x = \begin{bmatrix} 0 & 0 \\ 0 & [0,t] \end{bmatrix} \begin{bmatrix} [0,a] & [0,b] \\ 0 & 0 \end{bmatrix}$$

$$= \begin{bmatrix} 0 & 0 \\ 0 & 0 \end{bmatrix}.$$

$$x.p = \begin{bmatrix} [0,a] & [0,b] \\ 0 & 0 \end{bmatrix} \begin{bmatrix} 0 & 0 \\ 0 & [0,t] \end{bmatrix}$$

$$= \begin{bmatrix} 0 & [0,bt] \\ 0 & 0 \end{bmatrix} \neq (0).$$

$$q.y = \begin{bmatrix} [0,m] & 0 \\ 0 & 0 \end{bmatrix} \begin{bmatrix} 0 & 0 \\ [0,t] & [0,u] \end{bmatrix}$$

$$= \begin{bmatrix} 0 & 0 \\ 0 & 0 \end{bmatrix} \text{ and}$$

$$y.q = \begin{bmatrix} 0 & 0 \\ [0,t] & [0,u] \end{bmatrix} \begin{bmatrix} [0,m] & 0 \\ 0 & 0 \end{bmatrix}$$

$$= \begin{bmatrix} 0 & 0 \\ [0,tm] & 0 \end{bmatrix} \neq (0).$$

But pq = (0). Thus x, y in V is a Smarandache zero divisor of V.

It is still interesting to note that in general if x in V, V a matrix interval semiring is a Smarandache matrix interval anti zero divisor in V then in general x need not be a zero divisor in G.

Further it is important to note these interval matrix semirings can have right zero divisors or left zero divisors.

The author is left with the task of defining and giving examples of Smarandache idempotents in interval matrix semirings.

The study of Smarandache units in interval matrix semirings can be carried out as a matter of routine [ ]. Now in the following section we proceed onto define and describe interval polynomial semirings.

## 3.2 Interval Polynomial Semirings

In this section also we only concentrate on intervals of the form [0, a] where a is from $Z^+ \cup \{0\}$, $Q^+ \cup \{0\}$ or $R^+ \cup \{0\}$.

If x is any variable or an indeterminate we define interval polynomial semiring in an analogous way.

**DEFINITION 3.2.1:** *Let $S = \left\{ \sum_{i=0}^{\infty} [0, a_i] x^i \,\middle|\, a_i \in Z^+ \cup \{0\} \text{ or } R^+ \cup \{0\} \text{ or } Q^+ \cup \{0\} \right\}$ be the collection of all interval polynomials. S under addition is an abelian semigroup with $0 = [0, 0]$ as the zero polynomial which has the representation of the form $0 = [0, 0] + [0, 0]x + \ldots + [0, 0]x^n$ acts as the additive identity. S under multiplication is a semigroup. It is easily verified $p(x) \cdot [q(x) + r(x)] = p(x) \cdot q(x) = p(x) \cdot r(x)$ and $[a(x) + b(x)] \cdot c(x) = a(x) \cdot c(x) + b(x) \cdot c(x)$ for all $p(x), q(x), r(x), a(x), b(x)$ and $c(x)$ in S.*

*Thus (S, +, .) is a semiring known as the polynomial interval semiring in the variable x.*

Before we proceed to give examples of this structure we show how the addition and multiplication of interval polynomials are carried out.

Let p(x) = [0, 5] + [0, 3] $x^2$ + [0, 4] $x^3$ + [0, 9] $x^7$ and q (x) = [0, 3] + [0, 7] x + [0, 12] $x^2$ + [0, 40] $x^3$ + [0, 18] $x^5$ be any two interval polynomials in the variable x and the interval coefficients from $Z^+ \cup \{0\}$.

p (x) + q (x)
= ([0, 5] + [0, 3] $x^2$ + [0, 4] $x^3$ + [0, 9] $x^7$) + ([0, 3] + [0, 7] x + [0, 12] $x^2$ + [0, 40] $x^3$ + [0, 18] $x^5$)
= ([0, 5] + [0, 3]) + ([0, 3] + [0, 0]) $x^2$ + ([0, 4] + [0, 40]) $x^3$ + ([0, 0] + [0, 7]) x + ([0, 0] + [0, 18]) $x^5$ + ([0, 9] + [0, 0]) $x^7$
= [0, 8] + [0, 3] $x^2$ + [0, 44] $x^3$ + [0, 7] x + [0, 18] $x^5$ + [0, 9] $x^7$.

Now we will show how the product is defined p (x). q (x)
= ([0, 5] + [0, 3] $x^2$ + [0, 4]$x^3$ + [0, 9] $x^7$) ([0, 3] + [0, 7] x + [0, 12] $x^2$ + [0, 40] $x^3$ + [0, 18] $x^5$)

= [0, 5] [0, 3] + [0, 3] [0, 3] $x^2$ + [0, 4] [0, 3] $x^3$ + [0, 9] [0, 3] $x^7$ + [0, 5] [0, 7] x + [0, 3] [0, 7] $x^2$.x + [0, 4] [0, 7] $x^3$.x + [0, 9] [0, 7] $x^7$.x + [0, 5] [0, 12]$x^2$ + [0, 3] [0, 12] $x^2$.$x^2$ + [0, 4] [0, 12] $x^3$. $x^2$ + [0, 9] [0, 12] $x^7$. $x^2$ + [0, 5] [0, 40] $x^3$ + [0, 3] [0, 40] $x^2$. $x^3$ + [0, 4] [0, 40] $x^3$. $x^3$ + [0, 9] [0, 40] $x^7$.$x^3$ + [0, 5] [0, 18] $x^5$ + [0, 3] [0, 18] $x^2$. $x^5$ + [0, 4] [0, 18] $x^3$.$x^5$ + [0, 9] [0, 18] $x^7$.$x^5$

= [0, 15] + [0, 9] $x^2$ + [0, 12] $x^3$ + [0, 27] $x^7$ + [0, 35] x + [0, 21] $x^3$ + [0, 28] $x^4$ + [0, 63] $x^8$ + [0, 60] $x^2$ + [0, 36] $x^4$ + [0, 48] $x^5$ + [0, 108] $x^9$ + [0, 200] $x^3$ + [0, 120] $x^5$ + [0, 160] $x^6$ + [0, 360] $x^{10}$ + [0, 90] $x^5$ + [0, 54] $x^7$ + [0, 72] $x^8$ + [0, 162] $x^{12}$ = [0, 15] + [0, 35]x + ([0, 9] + [0, 60]) $x^2$ + ([0, 12] + [0, 21] + [0, 200]) $x^3$ + ([0, 28] + [0,

36]) $x^4$ + ([0, 48] + [0, 120] + [0, 90]) $x^5$ + [0, 160] $x^6$ + ([0, 27] + [0, 54]) $x^7$ + ([0, 63] + [0, 72]) $x^8$ + [0, 108] $x^9$ + [0, 360]$x^{10}$ + [0, 162]$x^{12}$

= [0, 15] + [0, 35] x + [0, 69]$x^2$ + [0, 233]$x^3$ + [0, 64]$x^4$ + [0, 258]$x^5$ + [0, 160]$x^6$ + [0, 81]$x^7$ + [0, 135] $x^8$ + [0, 108] $x^9$ + [0, 360] $x^{10}$ + [0, 162] $x^{12}$.

We see the interval polynomial multiplication as in case of polynomial multiplication is commutative and [0, 1] = [0,1] + [0,0]x + [0,0] $x^2$ + …+ [0,0] $x^n$, acts as the unit element.

We will illustrate the distributive law.

Let
$$a(x) = [0, 5] + [0, 8] x + [0, 1] x^8$$
$$b(x) = [0, 1] x + [0, 3] x^2 \text{ and}$$
$$c(x) = [0, 7] + [0, 3] x^4 + [0, 5] x^7$$

be three interval polynomials of an interval polynomial semiring.

To show
$$a(x) [b(x) + c(x)] = a(x) b(x) + a(x) c(x).$$

Consider
a (x) [b (x) + c (x)]

= ([0, 5] + [0, 8] x + [0, 1] $x^8$) [([0, 1]x + [0, 3] $x^2$ ) + ([0, 7] + [0, 3] $x^4$ + [0, 5] $x^7$)]

= ([0, 5] + [0, 8] x + [0, 1] $x^8$) + ([0, 1] x + [0, 3] $x^2$ + [0, 7] + [0, 3] $x^4$ + [0, 5] $x^7$)

= [0, 5] [0, 1] x + [0, 8] [0, 1]x.x + [0, 1] [0, 1] $x^8$.x + [0, 5] [0, 3] $x^2$ + [0, 8] [0, 3] x.$x^2$ + [0, 1] [0, 3] $x^8$.$x^2$ + [0, 5] [0, 7] + [0, 8] [0, 7] x + [0, 1] [0, 7] $x^8$ + [0, 5] [0, 3] $x^4$ + [0, 8] [0, 3] x.$x^4$ + [0, 1] [0, 3] $x^8$.$x^4$ + [0, 5] [0, 5] $x^7$ + [0, 8] [0, 5] x.$x^7$ + [0, 1] [0, 5] $x^8$.$x^7$

= [0, 35] + ([0, 5] + [0, 56]) x + ([0, 8] + [0, 15])$x^2$ + [0, 24] $x^3$ + [0, 15]$x^4$ + [0, 24] $x^5$ + [0, 25] $x^7$ + ([0, 7] + [0, 40]) $x^8$ + [0, 1] $x^9$ + [0, 3] $x^{10}$ + [0, 3] $x^{12}$ + [0, 5] $x^{15}$

=    [0, 35] + [0, 61] x + [0, 23] $x^2$ + [0, 24] $x^3$ + [0, 15] $x^4$ + [0, 24] $x^5$ + [0, 25] $x^7$ + [0, 47] $x^8$ + [0, 1] $x^9$ + [0, 3] $x^{10}$ + [0, 3] $x^{12}$ + [0, 5] $x^{15}$    ---- I

Consider
a(x) b(x) + a(x) c(x) =

=    ([0, 5] + [0, 8]x + [0, 1] $x^8$) ([0, 1]x + [0, 3] $x^2$) + ([0, 5] + [0, 8]x + [0, 1]$x^8$) ([0, 7] + [0, 3] $x^4$ + [0, 5] $x^7$)

=    [0, 5] [0, 1] x + [0, 8] [0, 1] x.x + [0, 1] [0, 1] $x^8$.x + [0, 5] [0, 3] $x^2$ + [0, 8] [0, 3] x.$x^2$ + [0, 1] [0, 3] $x^8$. $x^2$ + [0, 5] [0, 7] + [0, 8] [0, 7] x + [0, 1] [0, 7] $x^8$ + [0, 5] [0, 3] $x^4$ + [0, 8] [0, 3] x.$x^4$ + [0, 1] [0, 3] $x^8$.$x^4$ + [0, 5] [0, 5] $x^7$ + [0, 8] [0, 5] x.$x^7$ + [0, 1] [0, 5] $x^8$.$x^7$

=    [0, 35] + ([0, 5] [0, 1] + [0, 8] [0, 7]) x + ([0, 8] + [0, 15]) $x^2$ + ([0, 24] $x^3$ + [0, 15]) $x^4$ + [0, 24] $x^5$ + [0, 25] $x^7$ + ([0, 7] + [0, 40]) $x^8$ + [0, 1] $x^9$ + [0, 3] $x^{10}$ + [0, 3] $x^{12}$ + [0, 15] $x^{15}$

=    [0, 35] + [0, 61] x + [0, 23] $x^2$ + [0, 24] $x^3$ + [0, 15] $x^4$ + [0, 24] $x^5$ + [0, 25] $x^7$ + [0, 47] $x^8$ + [0, 1] $x^9$ + [0, 3] $x^{10}$ + [0, 3] $x^{12}$ + [0, 5] $x^{15}$    ---- II

It is easily verified I and II are equal, thus S is an interval polynomial semiring.
Now we will give examples of interval polynomial semirings.

*Example 3.2.1:* Let

$$S = \left\{ \sum_{i=0}^{\infty} [0, a_i] x^i \mid a_i \in R^+ \cup \{0\} \right\}$$

be an interval polynomial semiring.

***Example 3.2.3:*** Let $P = \left\{ \sum_{i=0}^{\infty} [0, a_i] x^i \mid a_i \in Z^+ \cup \{0\} \right\}$ be a polynomial interval semiring. Now we can construct polynomial interval semirings using $Z_n$; $n < \infty$ also.

***Example 3.2.4:*** Let $P = \left\{ \sum_{i=0}^{\infty} [0, a_i] x^i \mid a_i \in Z_9 \right\}$ be a polynomial interval semiring.

***Example 3.2.5:*** Let $R = \left\{ \sum_{i=0}^{\infty} [0, a_i] x^i \mid a_i \in Z_{19} \right\}$ be an interval polynomial semiring.

***Example 3.2.6:*** Let $T = \left\{ \sum_{i=0}^{20} [0, a_i] x^i \mid x^{21} = 1 \right.$ be an interval polynomial with coefficients from $Z_3 \}$ be an interval polynomial semiring.

***Example 3.2.7:*** Let $P = \left\{ \sum_{i=0}^{3} [0, a_i] x^i \mid x^4 = 1 \text{ and } a_i \in Z_2 \right\}$ be an interval polynomial semiring, P is finite. For P = {0, [0, 1], [0, 1] x, [0, 1] $x^2$ [0, 1] $x^3$, [0, 1] + [0, 1]x, [0, 1] + [0, 1]$x^2$, [0, 1] + [0, 1]$x^3$, [0, 1] + [0, 1] x + [0, 1] $x^2$ + [0, 1] $x^3$, [0, 1] + [0, 1] x + [0, 1] $x^2$, [0, 1] + [0, 1]x + [0, 1]$x^3$, [0, 1] + [0, 1]$x^2$ + [0, 1]$x^3$, [0, 1]x + [0, 1]$x^2$ + [0, 1] $x^3$ + [0, 1]x + [0, 1]$x^2$, [0, 1]x + [0, 1]$x^3$ [0, 1]$x^2$ + [0, 1]$x^3$} and |P| = 16.

Now we can define substructures in interval polynomial semirings.

**DEFINITION 3.2.2:** *Let $S = \left\{ \sum_{i=0}^{\infty} [0, a_i] x^i \mid a_i \in Z^+ \cup \{0\}, \text{ or } Z_n \text{ or } R^+ \cup \{0\} \text{ or } Q^+ \cup \{0\} \right\}$ be an interval polynomial semiring.*

Let $T = \left\{ \sum_{i=0}^{\infty} [0, a_i] x^{2i} \,\middle|\, a_i \in Z^+ \cup \{0\} \text{ or } Z_n \text{ or } R^+ \cup \{0\} \text{ or } Q^+ \cup \{0\} \right\} \subseteq S$; $T$ is an interval polynomial subsemiring of $S$.

**Example 3.2.8:** Let
$$S = \left\{ \sum_{i=0}^{\infty} [0, a_i] x^i \,\middle|\, a_i \in Q^+ \cup \{0\} \right\}$$
be an interval polynomial semiring.
Take
$$T = \left\{ \sum_{i=0}^{\infty} [0, a_i] x^i \,\middle|\, a_i \in Q^+ \cup \{0\} \right\} \subseteq S;$$
$T$ is an interval polynomial subsemiring of $S$.

**Example 3.2.9:** Let
$$V = \left\{ \sum_{i=0}^{\infty} [0, a_i] x^i \,\middle|\, a_i \in R^+ \cup \{0\} \right\}$$
be an interval polynomial semiring.
Take
$$W = \left\{ \sum_{i=0}^{\infty} [0, a_i] x^{3i} \,\middle|\, a_i \in R^+ \cup \{0\} \right\} \subseteq V,$$
$W$ is an interval polynomial subsemiring of $V$.

**Example 3.2.10:** Let
$$P = \left\{ \sum_{i=0}^{\infty} [0, a_i] x^i \,\middle|\, a_i \in Q^+ \cup \{0\} \right\}$$
be an interval polynomial semiring.
Take
$$T = \left\{ \sum_{i=0}^{\infty} [0, a_i] x^i \,\middle|\, a_i \in 3Z^+ \cup \{0\} \right\} \subseteq P$$
is an interval polynomial subsemiring of $P$.

Now having seen examples of interval polynomial semirings we now proceed onto define Smarandache interval polynomial semirings.

**DEFINITION 3.2.3**: *Let V be any interval polynomial semiring. A proper subset T of V, is such that T is an interval semifield then we define V to be a Smarandache polynomial interval semiring (S – polynomial interval semiring).*

We will illustrate this situation by some examples.

*Example 3.2.11:* Let

$$P = \left\{ \sum_{i=0}^{\infty} [0, a_i] x^i \mid a_i \in Z^+ \cup \{0\} \right\}$$

be a polylnomial interval semiring. Take $T = \{[0, a_i] / a_i \in Z^+ \cup \{0\}\} \subseteq P$, T is an interval semifield, so P is a S-polynomial interval semiring.

*Example 3.2.12:* Let

$$W = \left\{ \sum_{i=0}^{\infty} [0, a_i] x^i \mid a_i \in Q^+ \cup \{0\} \right\}$$

be an interval polynomial semiring. Choose $G = \{[0, a_i] / a_i \in Q^+ \cup \{0\}\} \subseteq W$, G is an S-interval polynomial semiring.

**THEOREM 3.2.1:** *Let*

$$V = \left\{ \sum_{i=0}^{\infty} [0, a_i] x^i \mid a_i \in Z^+ \cup \{0\} \text{ or } Q^+ \cup \{0\} \text{ or } R^+ \cup \{0\} \right\}$$

*be an interval polynomial semiring, then V is a S-interval polynomial semiring.*

The proof is direct and hence is left as an exercise for the reader to prove.

We can define S-polynomial interval subsemiring in an analogous way. We leave this task to the reader, however give some examples of S-interval polynomial subsemirings.

***Example 3.2.13:*** Let

$$T = \left\{ \sum_{i=0}^{\infty} [0, a_i] x^i \,\middle|\, a_i \in Q^+ \cup \{0\} \right\}$$

be an interval polynomial semiring.
Choose

$$W = \left\{ \sum_{i=0}^{\infty} [0, a_i] x^i \,\middle|\, a_i \in Z^+ \cup \{0\} \right\} \subseteq T,$$

W is a S-interval polynomial subsemiring as $G = \{[0, a_i] / a_i \in Z^+ \cup \{0\}\} \subseteq W$ is an interval semifield of W.

***Example 3.2.14:*** Let $V = \{\Sigma [0, a_i] x^i / a_i \in Z^+ \cup \{0\}\}$ be an interval polynomial semiring. $G = \{\Sigma [0, a_i] x^{2i} / a_i \in Z^+ \cup \{0\}\} \subseteq V$, G is a S-interval polynomial subsemiring of V for $M = \{\Sigma [0, a_i] / a_i \in Z^+ \cup \{0\}\} \subseteq G$ is an interval semifield.

Now we have the following theorem.

**THEOREM 3.2.2:** *Let V be an interval polynomial semiring. If V has a S-interval polynomial subsemiring then V is a S-interval polynomial semiring. But every interval polynomial subsemiring in general is not an S-interval polynomial subsemiring.*

The proof is direct and hence is left as an exercise for the reader.

We will an example of a S-interval polynomial semiring which has interval polynomial subsemiring which is not a S-interval polynomial subsemiring.

*Example 3.2.15:* Let

$$S = \left\{\sum_{i=0}^{\infty}[0,a_i]x^i \,\bigg|\, a_i \in Z^+ \cup \{0\}\right\}$$

be an interval polynomial semiring. Clearly S is a S-interval polynomial semiring.
Take

$$W = \left\{\sum_{i=0}^{\infty}[0,a_i]x^i \,\bigg|\, a_i \in 5Z^+ \cup \{0\}\right\} \subseteq S,$$

W is only an interval polynomial subsemiring which is not a S-interval polynomial subsemiring.

Now we give only examples of ideals in interval polynomial semirings, however the task of defining ideals in these semirings are left as an exercise for the reader.

*Example 3.2.16:* Let

$$G = \left\{\sum[0,a_i]x^i \,\bigg|\, a_i \in Z^+ \cup \{0\}\right\}$$

be an interval polynomial semiring.

Consider

$$W = \left\{\sum[0,a_i]x^i \,\bigg|\, a_i \in 3Z^+ \cup \{0\}\right\} \subseteq G$$

be an interval polynomial subsemiring of G. Clearly W is an ideal of the interval polynomial semiring G.

Infact this interval polynomial semiring G has infinitely many ideals.

We will give examples of interval polynomial semirings which has interval polynomial subsemirings which are not ideals. We know every ideal of an interval polynomial semiring is an interval polynomial subsemiring.

*Example 3.2.17:* Let

$$P = \left\{ \sum [0, a_i] x^i \mid a_i \in Q^+ \cup \{0\} \right\}$$

be an interval polynomial semiring.
Consider

$$G = \left\{ \sum_{i=0}^{\infty} [0, a_i] x^i \mid a_i \in Z^+ \cup \{0\} \right\} \subseteq P,$$

G is only an interval polynomial subsemiring which is not an ideal of P.

Now we will give examples of S-subsemirings which are not S-ideals of an interval polynomial semiring.

*Example 3.2.18:* Let

$$V = \left\{ \sum_{i=0}^{\infty} [0, a_i] x^i \mid a_i \in Q^+ \cup \{0\} \right\}$$

be an interval polynomial semiring.

$$P = \left\{ \sum_{i=0}^{\infty} [0, a_i] x^{2i} \,\middle|\, a_i \in Q^+ \cup \{0\} \right\} \subseteq V$$

be a S-interval polynomial subsemiring of V. P is not an S-ideal of V.

Concept of S-pseudo ideal and S-pseudo subsemiring in case of interval polynomial semirings can be studied, defined and analysed by the reader with appropriate examples.

However the concept of S-zero divisors, S-anti zero divisors or S-idempotents have no relevance in case of interval polynomial semirings built using $Z^+ \cup \{0\}$ or $Z_p$ (p a prime) or $Q^+ \cup \{0\}$ or $R^+ \cup \{0\}$. However when $Z_n$ is used n a compositive number these concepts have relevance.

**Chapter Four**

# GROUP INTERVAL SEMIRINGS

In this chapter we for the first time introduce the new notions of group interval semirings, semigroup interval semiring, loop interval semiring and groupoid interval semiring and discuss a few of the properties associated with them. This chapter has two sections. In section one we introduce the concept of group interval semirings and semigroup interval semirings and discuss a few of the properties associated with them. In section two we define loop interval semirings and study a few properties associated with them.

## 4.1 Group interval semirings and semigroup interval semirings

In this section we introduce the two new concepts viz. group interval semirings and semigroup interval semirings and study some properties related with them.

**DEFINITION 4.1.1**: *Let $S = \{[0, a_i] / a_i \in Z_n$ or $Z^+ \cup \{0\}$ or $Q^+ \cup \{0\}$ or $R^+ \cup \{0\}\}$ be the commutative semiring with unit $[0, 1]$. G be any group. The group semiring of the group G over the semiring S consists of all finite formal sums of the form $\sum_i [0, \alpha_i] g_i$ (i runs over finite number) where $[0, \alpha_i] \in S$ and $g_i \in G$ satisfying the following conditions.*

(i) $\sum_{i=1}^{n} [0, \alpha_i] g_i = \sum_{i=1}^{n} [0, \beta_i] g_i$ *if and only if* $\alpha_i = \beta_i$ *for* $i = 1, 2, \ldots, n$.

(ii) $\sum_{i=1}^{n} [0, \alpha_i] g_i + \sum_{i=1}^{n} [0, \beta_i] g_i = \sum_{i=1}^{n} [0, \alpha_i + \beta_i] g_i$

(iii) $\left( \sum_{i=1}^{n} [0, \alpha_i] g_i \right) \left( \sum_{i=1}^{n} [0, \beta_i] g_i \right) = \sum_{j=1}^{t} [0, \gamma_i] m_j$ *where* $g_i h_k = m_j$ *and* $\gamma_j = \Sigma [0, \alpha_i \beta_k]$

(iv) *For* $[0, \alpha_i] \in S$ *and* $g_i \in G$; $g_i [0, \alpha_i] = [0, \alpha_i] g_i$

(v) $[0, a] \sum_{i=1}^{n} [0, b_i] g_i = \sum_{i=1}^{n} [0, a][0, b_i] g_i = \sum_{i=1}^{n} [0, ab_i] g_i$ *for* $[0, a] \in S$ *and* $g_i \in G$.

*SG is an associative interval semiring with $0 \in S$ as its additive identity. Since $[0,1] \in S$ we have $G = [0, 1] G \subseteq SG$ and $S.e = S \subseteq SG$ (for $[0, 1] e = [0, 1]$ acts as the identity).*

*If we replace the group G by a semigroup with identity P, we get SP the semigroup interval semiring in the place of group interval semiring SG.*

We will illustrate this by some examples.

***Example 4.1.1:*** Let $S = \{[0, a] \mid a \in Z^+ \cup \{0\}\}$ be an interval semiring. $G = \langle g \mid g^2 = 1 \rangle$ be a cyclic group. $SG = \{[0, a] + [0, b]g \mid [0, a], [0, b] \in S; g \in G\}$ is the group interval semiring.

Clearly SG is of infinite order and is a commutative interval semiring with unit [0, 1].

***Example 4.1.2:*** Let $S = \{[0, a] \mid a \in Q^+ \cup \{0\}\}$ be an interval semiring with unit [0, 1]. $G = D_{2.7} = \{1, a, b \:/\: a^2 = 1 = b^7, bab = a\}$ be the group.

$$SG = \left\{ \sum_i [0, a_i]g_i \:\middle|\: a_i \in Q^+ \cup \{0\}, g_i \in D_{2.7} \right\},$$

the group interval semiring is of infinite order and is non commutative.

***Example 4.1.3:*** Let $S = \{[0, a] \:/\: a \in Z_7\}$ be an interval semiring and $G = \langle g \mid g^5 = 1 \rangle$ be the group, SG group interval semiring is of finite order and is commutative.

$SG = \{[0, a_0] + [0, a_1] g + [0, a_2] g^2 + [0, a_3] g^3 + [0, a_4] g^4 \:/\: g^5 = 1, [0, a_i] \in S, 0 \le i \le 4\}$. [0, 1] acts as the identity of SG.

***Example 4.1.4:*** Let $S = \{[0, a] \:/\: a \in Z_6 = \{0, 1, 2, 3, 4, 5\}\}$ be a finite commutative semiring with identity [0, 1]. $G = S_3 = \{1, p_1, p_2, p_3, p_4, p_5\}$ be the permutation group where

$$1 = \begin{pmatrix} 1 & 2 & 3 \\ 1 & 2 & 3 \end{pmatrix}, p_1 = \begin{pmatrix} 1 & 2 & 3 \\ 1 & 3 & 2 \end{pmatrix}, p_2 = \begin{pmatrix} 1 & 2 & 3 \\ 3 & 2 & 1 \end{pmatrix},$$

$$p_3 = \begin{pmatrix} 1 & 2 & 3 \\ 2 & 1 & 3 \end{pmatrix}, p_4 = \begin{pmatrix} 1 & 2 & 3 \\ 2 & 3 & 1 \end{pmatrix} \text{ and } p_5 = \begin{pmatrix} 1 & 2 & 3 \\ 3 & 1 & 2 \end{pmatrix}.$$

SG is the group interval semiring of finite order and it is non commutative.

$SG = \{[0, a_0] + [0, a_1] p_1 + [0, a_2] p_2 + [0, a_3] p_3 + [0, a_4] p_4 + [0, a_5] p_5 \mid a_i \in Z_6; 0 \le i \le 5\}$.

***Example 4.1.5:*** Let $S = \{[0, a] \,/\, a \in Q^+ \cup \{0\}\}$ be a commutative interval semiring with unit. Choose $G = \{g \,/\, g^{125} = 1\}$ be a cyclic group of finite order.

$$SG = \left\{\sum_i [0,a]g^i \,\middle|\, a_i \in Q^+ \cup \{0\}\right\}; g^i \in G, 1 \le i \le 125\}$$

is the group interval semiring of the group G over the interval semiring S. Clearly SG is a commutative interval semiring of infinite order with the unit [0, 1].

It is important to mention here that the method of constructing group interval semirings leads to get an infinite class of both finite and infinite and commutative and non commutative interval semirings.

We can define substructures and special elements in group interval semirings.

We see $S = \{[0, a] \,/\, a \in Z^+ \cup \{0\}$ or $Q^+ \cup \{0\}$ or $Z_p$ or $R^+ \cup \{0\}\}$ (p a prime) are interval semirings which are interval semifields.

We will describe some zero divisors, units and idempotents in these group interval semirings.

***Example 4.1.6:*** Let $S = \{[0, a] \,/\, a \in Z_4\}$ be an interval semiring. $G = \langle g \,/\, g^4 = 1 \rangle$ be any group. The group interval semiring SG has zero divisors, units and idempotents.

Consider $\alpha = [0, 1] + [0, 1]g + [0, 1]g^2 + [0, 1]g^3 \in SG$. Consider

$$\begin{aligned}\alpha^2 &= ([0, 1] + [0, 1]g + [0, 1]g^2 + [0, 1]g^3)^2 \\ &= ([0, 1]) + ([0, 1])g^2 + ([0, 1]g^2)^2 + ([0, 1]g^3)^2 + 2 [0, 1] \\ &\quad [0, 1]g + 2 [0, 1] [0, 1]g^2 + 2 [0, 1] [0, 1]g^3 + 2 [0, 1] [0, \\ &\quad 1] g.g^2 + 2 [0, 1] [0, 1] g.g^3 + 2 [0, 1] [0, 1] g^2 . g^3 \end{aligned}$$

$$\begin{aligned}
= \quad & [0, 1] + [0, 1]g^2 + [0, 1]g^4 + [0, 1]g^6 + 2\,[0, 1]g + 2\,[0, 1]g^2 + 2[0, 1]g^3 + 2\,[0, 1]g^3 + 2\,[0, 1]g^4 + 2\,[0, 1]g^5 \\
= \quad & [0, 1] + [0, 1]g^2 + [0, 1] + [0, 1]g^2 + 2\,[0, 1]g + 2\,[0, 1]g^2 + 2\,[0, 1]g^3 + 2\,[0, 1]g^3 + 2\,[0, 1] + 2\,[0, 1]g \\
= \quad & 0.
\end{aligned}$$

Thus $\alpha^2 = 0$ in SG is nilpotent. Also if $b = [0, 2] + [0, 2]g$ then
$$b^2 = ([0, 2])^2 + 2\,([0, 2][0, 2])g + ([0, 2])g)^2$$
$= 0$ is again an nilpotent element of order two. Take $x.y = ([0, 2]g + [0, 2])([0, 2]g^2 + [0, 2]g^3) = 0$ where $x = [0, 2]g + [0, 2]$ and $y = [0, 2]g^2 + [0, 2]g^3$ are in SG. Thus $x, y \in$ SG is a zero divisor in SG.

***Example 4.1.7:*** Let $S = \{[0, a] \,/\, a \in Z_2\}$ and $G = \{g \,/\, g^4 = 1\}$ be the interval semiring and group respectively. SG the group interval semiring of G over the interval semiring S.

SG = $\{0, [0, 1], [0, 1]g, [0, 1]g^2, [0, 1]g^3, [0, 1] + [0, 1]g + [0, 1]g^2 + [0, 1]g^3, [0, 1] + [0, 1]g, [0, 1] + [0, 1]g^2, [0, 1] + [0, 1]g^3, [0, 1]g + [0, 1]g^2, [0, 1]g + [0, 1]g^3, [0, 1]g^2 + [0, 1]g^3, [0, 1] + [0, 1]g + [0, 1]g^2, [0, 1] + [0, 1]g + [0, 1]g^3, [0, 1] + [0, 1]g^2 + [0, 1]g^3 + [0, 1]g + [0, 1]g^2 + [0, 1]g^3\}$.

Now
$$\alpha = [0, 1] + [0, 1]g + [0, 1]g^2 + [0, 1]g^3$$
in SG is such that
$$\begin{aligned}
\alpha^2 = \quad & [0, 1] + [0, 1]g^2 + [0, 1]g^4 + [0, 1]g^6 + 2\,[0, 1]g + 2\,[0, 1]g^2 + 2\,[0, 1]g^3 + 2\,[0, 1]g^4 \\
= \quad & [0, 1] + [0, 1]g^2 + [0, 1]g^4 + [0, 1]g^6 + 2[0, 1]g + 2[0, 1]g^2 + 2[0, 1]g^3 + 2\,[0, 1]g4 \\
= \quad & [0, 1] + [0, 1]g^2 + [0, 1]g^0 + [0, 1]g^2 \\
= \quad & 0
\end{aligned}$$

$\alpha$ is an nilpotent element of G. This SG has several zero divisors.

Take $\alpha = [0, 1] + [0, 1]g + [0, 1]g^2 + [0, 1]g^3$ and $\beta = [0, 1] + [0, 1]g$ in SG. Clearly $\alpha.\beta = 0$.

*Example 4.1.8:* Let $S = \{[0, a] / a \in Z_3\}$ be an interval semiring and $G = \langle g \mid g^2 = 1 \rangle$ be the cyclic group of order two. The interval group semiring or group interval semiring SG = $\{[0, 1]$, 0, $[0, 2]$, $[0, 1]g$, $[0, 2]g$, $[0, 1] + [0, 1]g$, $[0, 1] + [0, 2]g$, $[0, 2] + [0, 1]g$, $[0, 2] + [0, 2]g\}$.

Now $([0, 2]g)^2 = [0, 1]$ is a unit in SG. Consider $\alpha = [0, 2] + [0, 2]g$ in SG.

$$\begin{aligned}\alpha^2 &= [0, 2]^2 + [0, 2]2\,g^2 + 2\,[0, 2]\,[0, 2]\,g \\ &= [0, 1] + [0, 1] + 2\,[0, 1]g \\ &= [0, 2] + [0, 2]\,g \\ &= \alpha.\end{aligned}$$

Thus $\alpha$ is an idempotent in SG. Consider $\beta = [0, 2] + [0, 1]g$ in SG.

$$\begin{aligned}\beta^2 &= ([0, 2] + [0, 1]g)^2 \\ &= [0, 2]^2 + [0, 1]g^2 + 2\,[0, 2]\,[0, 1]g \\ &= [0, 1] + [0, 1] + [0, 1]g \\ &= [0, 2] + [0, 1]g = \beta.\end{aligned}$$

Thus $\beta$ is an idempotent in SG.

*Example 4.1.9:* Let $S = \{[0, a] / a \in Z_3\}$ be an interval semiring and $G = \{g / g^3 = 1\}$ be a cyclic group of order three. SG = $\{0$, $[0, 1]$, $[0, 2]$, $[0, 1]g$, $[0, 2]g$, $[0, 1]g^2$, $[0, 2]g^2$, $[0, 1] + [0, 1]g$, $[0, 1] + [0, 2]g$, $[0, 2] + [0, 1]g$, $[0, 2] + [0, 2]g + [0, 1] + [0, 1]g^2$, $[0, 1] + [0, 2]g^2$, $[0, 2] + [0, 2]g^2$, $[0, 2] + [0, 1]g^2$, $[0, 1]g + [0, 1]g^2$, $[0, 1]g + [0, 2]g^2$, $[0, 2]g + [0, 1]g^2$, $[0, 2]g + [0, 2]g^2$ $[0, 1] + [0, 1]g + [0, 1]g^2$, $[0, 2] + [0, 1]g + [0, 1]g^2$, $[0, 1] + [0, 2]g + [0, 1]g^2$, $[0, 1] + [0, 1]g + [0, 2]g^2$, $[0, 2] + [0, 1]g + [0, 1]g^2$, $[0, 2] + [0, 2]g + [0, 1]g^2 + [0, 2] + [0, 1]g + [0, 2]\,g^2$, $[0, 2] + [0, 2]g + [0, 2]g^2\}$.

Let $\alpha = [0, 2] + [0, 2]g + [0, 2]\,g^2 \in$ SG is such that $\alpha^2 = 0$.

Now we work with infinite group interval semiring SG.

**Example 4.1.10:** Let $S = \{[0, a] / a \in Z^+ \cup \{0\}\}$ be an interval semiring. $G = \langle g / g^2 = 1 \rangle$ be a cyclic group of order two. SG the group interval semiring. $SG = \{[0, a] + [0, b] g / [0, a], [0, b] \in S\}$.

Take $\alpha = [0, 1]g \in SG$; $\alpha^2 = [0, 1]$ is a self inversed element of SG. We see if $a = [0, a] + [0, b] g$ then $a^n = [0, x] + [0, y]g$ for all n and $[0, a] < [0, x]$ and $[0, b] < [0, y]$ thus as $n \to \infty$ $a^n \to [0, t] + [0, u] g$ where t and u are very large.

Clearly SG has no zero divisors.

**Example 4.1.11:** Let $S = \{[0, a] / a \in Z^+ \cup \{0\}\}$ be an interval semiring. $G = \langle g / g^3 = 1 \rangle$ be a group of order three. The group interval semiring SG has no zero divisors.

Now we proceed onto define substructures in group interval semirings SG.

Let us illustrate them by examples and the definition of substructures can be defined as a matter of routine.

**Example 4.1.12:** Let $S = \{[0, a] / a \in Z^+ \cup \{0\}\}$ be an interval semiring. $G = S_3$ the symmetric group of degree three, $SG = \{\Sigma [0, a]g / g \in G, a \in Z^+ \cup \{0\}\}$ is the group interval semiring of the group G over the interval semiring S.
Consider

$$H = \left\{1 = \begin{pmatrix} 1 & 2 & 3 \\ 1 & 2 & 3 \end{pmatrix}, p_1 = \begin{pmatrix} 1 & 2 & 3 \\ 1 & 3 & 2 \end{pmatrix}\right\} \subseteq G.$$

H is a subgroup of G. $SH = \{\Sigma [0, a] h / h \in H, a \in Z^+ \cup \{0\}\} \subseteq SG$, SH is an interval subsemiring of SG.
Take

$$p = \left\{1 = \begin{pmatrix} 1 & 2 & 3 \\ 1 & 2 & 3 \end{pmatrix}, p_4 = \begin{pmatrix} 1 & 2 & 3 \\ 2 & 3 & 1 \end{pmatrix}, p_5 = \begin{pmatrix} 1 & 2 & 3 \\ 3 & 1 & 2 \end{pmatrix}\right\} \subseteq G,$$

SP is an interval subsemiring of SG.

Now we can define interval subsemiring also as follows.
TG = {Σ [0, a] g / g ∈ G, a ∈ $2Z^+ \cup \{0\}$} ⊆ SG; TG is an interval subsemiring. (This subsemiring has no unit)

***Example 4.1.13:*** Let S = {[0, a] / a ∈ $Z^+ \cup \{0\}$} be an interval semiring. G = ⟨g | $g^p$ = 1⟩ be a cyclic group of order p, p a prime.

$$SG = \left\{ \sum_i [0,a]g^i \,\middle|\, g^i \in G; a \in Z^+ \cup \{0\} \right\}$$

is the group interval semiring of the group G over the interval semiring S. G is simple, as G has no proper subgroups. However

$$PG = \left\{ \sum_i [0,a]g^i \,\middle|\, a \in 5Z^+ \cup \{0\} \right\} \subseteq SG.$$

Now P = {[0, a] / a ∈ $5Z^+ \cup \{0\}$} ⊆ S = {[0, a] / a ∈ $Z^+ \cup \{0\}$}, so PG ⊆ SG is an group interval subsemiring of SG which has no unit.

***Example 4.1.14:*** Let S = {[0, a] / a ∈ $Z_4$} be an interval semiring. Consider G = ⟨g / $g^5$ = 1⟩ be a cyclic group of degree 5. SG be the group interval semiring of the group G over the semiring S. SG = {[0, 1], [0, 2], [0, 3], 0, $\sum_{i=0}^{4}[0,a]g^i$ | a ∈ $Z_4$}.
Consider P = {0, [0, 2], [0, 2]g, [0, 2]$g^2$, [0, 2]$g^3$, [0, 2]$g^4$, [0, 2] + [0, 2]g, [0, 2] + [0, 2]$g^2$, [0, 2] + [0, 2]$g^3$, [0, 2] + [0, 2]$g^4$, [0, 2]g + [0, 2]$g^2$, [0, 2]g + [0, 2]$g^3$, [0, 2]g + [0, 2]$g^4$, [0, 2]$g^2$, [0, 2]$g^3$, [0, 2]$g^2$ + [0, 2]$g^4$, [0, 2]$g^3$ + [0, 2]$g^4$, [0, 2] + [0, 2]g + [0, 2]$g^2$, …, [0, 2] + [0, 2]g + [0, 2]$g^2$ + [0, 2]$g^4$ + [0, 2]$g^3$} ⊆ SG. P is an ideal of SG. P is also an interval subsemiring of S.

One can define left ideals and right ideals of group interval semirings SG, when the group G is non commutative.

It is interesting to note that none of the group interval semirings SG when S = {[0, a] | a ∈ $Z^+$ ∪ {0} or $R^+$ ∪ {0} or $Q^+$ ∪ {0}} have zero divisors for any group G. Infact they are also strict interval semirings and S-interval semirings. When G is a commutative group SG turns out to be an interval semifield.

But when SG, where S = {[0, a] | a ∈ $Z_n$, n < ∞} is taken SG can have zero divisors, nilpotents and idempotents. Infact S is not an interval semifield as [0, a] + [0, b] = 0 (mod n) can occur in which a ≠ 0 and b ≠ 0.

Thus question of SG becoming an interval semifield is ruled out for S built using $Z_n$, n < ∞.

Infact all interval semirings S = {[0, a] | a ∈ $Z^+$ ∪ {0} or $Q^+$ ∪ {0} or $R^+$ ∪ {0}} happen to be a complete partially ordered set with least element [0, 0] = 0 such that the sum and product operations are continuous. Interested author can study these interval semirings under * operation and characterize those interval * semirings.

A study of group interval semirings SG where G is a group of order n (n a non prime) and S = {[0, a] | a ∈ $Z_p$} p a prime such that (1) where p / n. (2) when p ∤ n.

Also the question when G is a cyclic group of order p, p a prime and S = {[0, a] / a ∈ $Z_n$}, n a composite number and p / n and p ∤ n can be studied by interested reader.

Now we proceed onto study semigroup interval semirings SP where S is an interval semiring and P is a multiplicative semigroup.

In the definition of group interval semiring if we replace the group by the semigroup we get the notion of semigroup interval semiring.

We give some examples before we proceed onto study the properties associated with them.

***Example 4.1.15:*** Let S = {[0, a] / a ∈ $Z^+$ ∪ {0}} be an interval semiring. Let G = ($Z_4$, ×) be the semigroup given by the following table.

| × | 0 | $\bar{1}$ | $\bar{2}$ | $\bar{3}$ |
|---|---|---|---|---|
| 0 | 0 | 0 | 0 | 0 |
| $\bar{1}$ | 0 | $\bar{1}$ | $\bar{2}$ | $\bar{3}$ |
| $\bar{2}$ | 0 | $\bar{2}$ | 0 | $\bar{2}$ |
| $\bar{3}$ | 0 | $\bar{3}$ | $\bar{2}$ | $\bar{1}$ |

Now SG = {[0, a] + [0, b] $\bar{2}$ + [0, c] $\bar{3}$ / a, b, c ∈ $Z^+$ ∪ {0}} is an interval semigroup semiring. SG has zero divisors.

For take [0, b] $\bar{2}$ . [0, c] $\bar{2}$ = 0

For all a, b ∈ $Z^+$. [0, 1] $\bar{3}$ . [0, 1] $\bar{3}$ = [0, 1] is unit. It is to be noted that SG is an infinite semigroup interval semiring. SG has nilpotent element of order two. For all ([0, a] $\bar{2}$ )² = 0 for all a ∈ $Z^+$.

Now we proceed onto give another example of a finite semigroup interval semiring.

***Example 4.1.16:*** Let S = {[0, a]| a ∈ $Z_6$} be an interval semiring. Let G = ($Z_3$, ×) be the semigroup, under multiplication modulo three given by the table.

| × | 0 | $\bar{1}$ | $\bar{2}$ |
|---|---|---|---|
| 0 | 0 | 0 | 0 |
| $\bar{1}$ | 0 | $\bar{1}$ | $\bar{2}$ |
| $\bar{2}$ | 0 | $\bar{2}$ | $\bar{1}$ |

SG = {[0, a] + [0, b] $\bar{2}$ / a, b ∈ $Z_6$} is the semigroup interval semiring of finite order. Further SG is commutative.

SG = {0, [0, 1], [0, 2], [0, 3], [0, 4], [0, 5], [0, 1] $\bar{2}$ , [0, 2] $\bar{2}$ , [0, 4] $\bar{2}$ , [0, 3] $\bar{2}$ , [0, 5] $\bar{2}$ , [0, 1] + [0, 2] $\bar{2}$ , [0, 1] + [0, 1] $\bar{2}$ , [0, 1] + [0, 3] $\bar{2}$ , [0, 1] + [0, 4] $\bar{2}$ , …, [0, 5] + [0, 5] $\bar{2}$ }.

SG has units and zero divisors.

[0, 2] $\bar{2}$ . [0, 3] $\bar{2}$ = [0, 0] $\bar{1}$ = 0

[0, 3] . [0, 2] $\bar{2}$ = 0 and so on.

Also [0, 5] [0, 5] = [0, 1], [0, 5] $\bar{2}$ / [0, 5] $\bar{2}$ = [0, 1] and so on. [0, 3]. [0, 3] = [0, 3] is an interval idempotent of SG.

*Example 4.1.17:* Let S = {[0, a]| a ∈ $Z_2$} be an interval semiring. G = ($Z_6$, ×) be the semigroup under multiplication.

SG = {[0, 1] + [0, 1] $\bar{2}$ + [0, 1] $\bar{3}$ + [0, 1] $\bar{4}$ + [0, 1] $\bar{5}$, 0, [0, 1], [0, 1] + [0, 1] $\bar{2}$, …, [0, 1] $\bar{3}$ + [0, 1] $\bar{4}$ + [0, 1] $\bar{2}$ … } be the semigroup interval semiring.

Clearly SG is of finite order and |SG| = 32.

Consider x = [0, 1] +[0, 1] $\bar{2}$, + [0, 1] $\bar{4}$ in SG. We see $x^2$ = [0, 1], so x is a unit in SG.

Consider y = [0, 1] + [0, 1] $\bar{3}$ ∈ SG; $y^2$ = [0, 1] + [0, 1] $\bar{3}$ = y, so y is an idempotent of SG.

Take p = [0, 1] + [0, 1] $\bar{5}$ ∈ SG.

$p^2$ = (0). Several other properties related with SG can be obtained by the interested reader.

*Example 4.1.18:* Let S = {[0, a] / a ∈ $Z_{12}$} be an interval semiring. G = {$Z_6$, ×} be the semigroup under multiplication. SG be the semigroup interval semiring of G over S.

SG has zero divisors, units and idempotents in it.

Now we give examples of semigroup interval subsemirings.

*Example 4.1.19:* Let S = {[0, a] | a ∈ $Z_2$} be the interval semiring. Choose G = {$Z_{12}$, ×} the semigroup under multiplication modulo 12. SG is the semigroup interval semiring of the semigroup G over the interval semiring S. Consider H = {0, $\bar{2},\bar{1},\bar{4},\bar{6},\bar{8},\bar{10}$, ×} ⊆ G; H is a subsemigroup of G. The semigroup interval semiring SH of the semigroup H over the interval semiring S is a semigroup interval subsemiring of SG.

*Example 4.1.20:* Let $S = \{[0, a] \mid a \in Z_{24}\}$ be an interval semiring. $G = \{Z_5, \times\}$ be the semigroup under multiplication modulo 5. Let

$$SG = \left\{\sum_g [0,a]g \,\middle|\, a \in Z_{24}\right\}$$

be the semigroup interval semiring of the semigroup G over the interval semiring S.

Now consider $W = \{\Sigma [0, a] g \mid g \in G$ and $a \in \{0, 2, 4, 6, 8, \ldots, 22\} \subseteq Z_{24}\} \subseteq SG$; W is a semigroup interval subsemiring of SG.

*Example 4.1.21:* Let $S = \{[0, a] \mid a \in Q^+ \cup \{0\}\}$ be an interval semiring. Let $G = \{Z_{19}, \times\}$ be a semigroup under multiplication modulo 19. SG be the semigroup interval semiring of the semigroup G over the interval semiring S.

$$SG = \left\{\sum_g [0,a]g \,\middle|\, a \in Q^+ \cup \{0\} \text{ and } g \in G\right\}.$$

Now let $M = \{\Sigma [0, a] g \mid a \in Z^+ \cup \{0\}$ and $g \in G\} \subseteq SG$, clearly M is an interval semiring and is the semigroup interval subsemiring of SG. Infact SG has infinite number of semigroup interval subsemirings. The order of SG is infinite and the order of every semigroup interval subsemiring is also infinite.

*Example 4.1.22:* Let $S = \{[0, a] \mid a \in Z^+ \cup \{0\}\}$ be an interval semiring. $G = \{Z_{71}, \times\}$ be the semigroup under multiplication. SG be the semigroup interval semiring of the semigroup G over the interval semiring S. Now let $I = \{\Sigma [0, a] g \mid g \in G$ and $a \in 5Z^+ \cup \{0\}\} \subseteq SG$; I is an ideal of SG.

*Example 4.2.23:* Let $S = \{[0, a] \mid a \in Z_8\}$ be an interval semiring. $G = \{Z_5, \times\}$ be the semigroup under multiplication modulo 5. $SG = \{\Sigma[0, a]g \mid a \in Z_8$ and $g \in G\}$ be the semigroup interval semiring of the semigroup G over the interval semiring

S. Now let $J = \{\Sigma\, [0, a]\, g\, /\, g \in G$ and $a \in \{0,\, \overline{2}, \overline{4}, \overline{6}\,\} \subseteq Z_8\} \subseteq$ SG, J is an ideal of SG. Take $\alpha = [0, a] + [0, b]\, \overline{2} + [0, b']\, \overline{4} + [0, c]\, \overline{6}$ in J. It is easily verified J is an ideal of SG.

Interested reader can study more about ideals and subsemirings in a semigroup interval semiring. We now consider other types of semigroups and study them.

*Example 4.1.24*: Let S(3) be the set of all mappings of (123) to itself S(3) is the semigroup under the operation of composition of maps. Let $S = \{[0, a]\, /\, a \in Z_3\}$ be the interval semiring. SS (3) is the semigroup interval semiring of the semigroup S(3) over the interval semiring S. Clearly order of SS (3) is finite and SS(3) is a non commutative interval semiring.

Further it is interesting to note that SS (3) contains a group interval semiring $SS_3$ which is a semigroup interval subsemiring. We can also call SS (3) as a S-semigroup interval semiring as S (3) is a S-semigroup for $S_3 \subseteq S$ (3).

*Example 4.1.25:* Let S (3) be the symmetric semigroup [16]. $S = \{[0, a]\, |\, a \in Z^+ \cup \{0\}\}$ be an interval semiring. $SS(3) = \{\Sigma[0, a]g\, |\, g \in S(3)$ and $a \in Z^+ \cup \{0\}\}$ is a semigroup interval semiring of the semigroup S(3) over the interval semiring S. Clearly SS(3) is of infinite cardinality and is a noncommutative interval semiring. Infact SS(3) is a S-interval semiring.

*Example 4.1.26:* Let S(7) be a symmetric semigroup. $S = \{[0, a]\, |\, a \in Z_2\}$ be the interval semiring. SS(7) is the symmetric semigroup interval semiring of the symmetric semigroup S(7) over the interval semiring S. Clearly over of SS(7) is finite.

Further SS(7) is S-interval semiring. $SS_7 \subseteq SS(7)$ is a symmetric group interval semiring and is a interval subsemiring of S. Clearly SS (7) has zero divisors. Further SS (7) has non trivial idempotents. Several other interesting features about this semiring can be analysed by an interested reader.

However SS (7) is not a semigroup interval S-semiring.

*Example 4.1.27:* Let $S = \{[0, a] \mid a \in Z_7\}$ be an interval semiring and S(m) be a symmetric semigroup (m, n $< \infty$). SS(m) the semigroup interval semiring of the semigroup S(m) over the interval semiring S.S is not a S-semiring as SS (m) is non commutative and not a strict semiring.

Now we give examples of interval semidivision ring.

*Example 4.1.28:* Let $S = \{[0, a] \mid a \in Z^+ \cup \{0\}$ or $Q^+ \cup \{0\}$ or $R^+ \cup \{0\}\}$ be an interval semiring. Take S (n) to be the symmetric semigroup. The semigroup interval semiring SS (n) is non commutative but is a strict semiring and SS (n) has no zero divisors. Hence SS (n) is a division interval semiring. Also SS (n) is a S-semiring.

Now having defined and seen the concept of semigroup interval semigroup and group interval semiring we now proceed onto define the notion of interval semigroup semiring and interval group semiring which we choose to call as special interval semirings.

**DEFINITION 4.1.2:** *Let S be any commutative semiring with unit which is not an interval semiring. P an interval semigroup. SP $= \{\Sigma\, a_i\, [0, p_i] \mid a_i \in S$ and $[0, p_i] \in P\}$ under usual addition and multiplication is a semiring defined as the interval semigroup semiring or special interval semiring.*

We will give examples of these special interval semirings.

*Example 4.1.29:* Let $S = \{Z^+ \cup \{0\}, +, \times\}$ be a semiring. Choose $P = \{[0, a] \mid a \in Z_{12}\}$ to be an interval semigroup under multiplication. SP be the interval semigroup semiring of the interval semigroup P over the semiring S. SP = $\{\Sigma a_i [0, p_i] \mid p_i \in Z_{12}, a_i \in Z^+ \cup \{0\}\}$. The following observations are direct

(i) SP is not a semifield for x = 5 $[0, \overline{6}]$ and y = 9 $[0, \overline{2}]$ in SP are such that xy = 0. So SP has zero divisors.
(ii) SP has idempotents, for take a = $[0, \overline{4}] \in$ SP is such that $a^2 = a$.
(iii) SP is a S-interval semiring as S ⊆ SP for 1. [0,1] = 1 is the identity element of SP; so S = S [0,1] = S ⊆ SP and S is a semifield.
(iv) S is strict interval semiring as s + p = 0 if and only if s = 0 and p = 0.
(v) SP has nilpotents for consider t = n [0, 6] ∈ SP $t^2 = n^2$ [0, 0] = 0.

**Example 4.1.30:** Let S = $\{Q^+ \cup \{0\}\}$ be a semifield. Consider P = $\{[0, a] \mid a \in Z^+ \cup \{0\}\}$ an interval semiring. Clearly SP the interval semigroup semiring is of infinite order. SP = {finite formal sums of the form $\Sigma a_i [0, x_i]$ where $a_i \in$ S and $[0, x_i] \in$ P. $a_i [0, 0] = 0$ for all $a_i \in Q^+ \cup \{0\}$ and $a_i [0, 1] \neq a_i$ but 1 [0, 1] is assumed to be the multiplicative identity of SP.

It is easily verified SP is an infinite special interval semiring. Further SP is commutative. Also SP is a strict interval special semiring.

Thus SP itself is a semifield.

Thus we by this method construct semifields of infinite cardinality.

**Example 4.1.31:** Let S = $\{R^+ \cup \{0\}\}$ be a semiring. W = $\{[0, a] / a \in R^+ \cup \{0\}\}$ be an interval semigroup. Clearly SW = {finite formal sums of the form $\Sigma a_i [0, w_i] / a_i \in$ S and $[0, w_i] \in$ W} under addition and multiplication is a semiring, called the special interval semiring or the interval semigroup semiring of the interval semigroup W over the semiring S. SW is a semifield of infinite dimension bigger than S.

Now by this method we can use distributive lattices to construct interval semigroup semirings of finite order. We know all distributive lattices are semirings. We will using these lattices to construct interval semigroup rings.

*Example 4.1.32:* Let

$$S = \{\begin{matrix} \circ\, 1 \\ | \\ \circ\, 0 \end{matrix} \text{ be the distributive lattice}\}$$

which is a semiring. $P = \{[0, a] \;/\; a \in Z_7\}$ be the interval semigroup under multiplication. SP is an interval semigroup semiring where $SP = \{0, 1\,[0, a], \Sigma\, x\, [0, a]$ where $x = 0$ or 1 and $a \in Z_7\}$. Clearly SP is of finite order. SP has no zero divisors or idempotents but SP has units for $1\,[0, 6]$. $1\,[0, 6] = 1\,[0, 1]$, $1\,[0, 4]$. $1\,[0, 2] = 1\,[0, 8] = 1\,[0, 1]$ and so on.

*Example 4.1.33:* Let $S = \{$ the chain lattice given by $0, a_1, a_2, a_3, a_4, 1, a_i^2 = a_i, a_i\, 1 = a_i, a_i\, 0 = 0, 0 < a_1 < a_2 < a_3 < a_4 < 1; 1 \le i \le 4\}$ is a semiring. $P = \{[0, a] \;/\; a \in Z_6\}$ is the interval semigroup under multiplication. $SP = \{\Sigma\, x\, [0, a]$ finite formal sums where $x \in S$ and $a \in Z_6\}$ is the interval semigroup semiring of finite order. This special interval semiring has zero divisors units and idempotents. However SP is a S-semiring.

*Example 4.1.34:* Let $S = \{0, a, b, 1$, a lattice given by the diagram

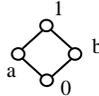

$\}$ which is a semiring. Let $P = \{[0, a] \;/\; a \in Z^+ \cup \{0\}\}$ be an interval semigroup. $SP = \{\Sigma\, x\, [0, a] \;/\; a \in Z^+ \cup \{0\}$ and $x \in \{0, 1, a, b\} = S\}$ is a interval semigroup P over the semiring S. SP has zero divisors, for $a\,[0, x]\,.\, b\,[0, y] = a\,.\,b\,[0, xy]$ where $a, b \in S$ and $x, y \in P$; as $a\,.\,b = 0$ we get $a\,[0, x]\,.\, b\,[0, y] = 0$. This interval semigroup semiring is of infinite order.

***Example 4.1.35:*** Let S = {0, 1, a, b, c, d, e, f |

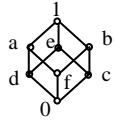

be the semiring (which is a distributive lattice) and P = {[0, a] / a ∈ R$^+$ ∪ {0}} be an interval semigroup. Consider SP = {Σ x [0, a], finite formal sums with x ∈ S and a ∈ R$^+$ ∪ {0} or [0, a] ∈ P} is the interval semigroup semiring of infinite order. Clearly SP has zero divisors and units.

Thus using the special type of semirings which are finite lattices we can construct interval semigroup semirings of finite order as well as infinite order.

Now we proceed onto define interval group under addition and multiplication. However we will be using only interval groups under multiplication.

**DEFINITION 4.1.3:** *Let G = {[0, a] | a ∈ $Z_n$; n < ∞}, G is a group under addition of finite order defined as interval group under addition.*

***Example 4.1.36:*** Let G = {[0, a] | a ∈ $Z_8$}, G is an interval group under addition.

***Example 4.1.37:*** Let G = {[0, a] | a ∈ $Z_{17}$}, G is an interval group under addition of order 17. Clearly G is a commutative group of finite order.

However it is important and interesting to note when we consider intervals from Z$^+$ ∪ {0} or Q$^+$ ∪ {0} or R$^+$ ∪ {0} that interval set is not an interval group under addition. We see S = {[0, a] / a ∈ Z$^+$ ∪ {0}}. Likewise S = {[0, a] / a ∈ Q$^+$ ∪ {0} or $R^+ \cup \{0\}\}$ is not an interval group under addition. Thus we now proceed onto define interval group under multiplication.

**DEFINITION 4.1.3:** *Let $S = \{[0, a] \mid a \in Q^+\}$ be an interval set. S is an interval group under multiplication of infinite order.*

*If we replace $Q^+$ by $R^+$ then also S is an interval group under multiplication of infinite order.*

*However if we replace $Q^+$ by $Z^+$, S is not an interval group under multiplication only an interval semigroup under multiplication.*

*Now if we replace $Q^+$ by $Z_p \setminus \{0\}$ p any prime $p \geq 3$ then S is an interval group under multiplication of finite order. Thus using $Z_p \setminus \{0\}$ we get a collection of interval groups under multiplication of finite order which is commutative.*

We will illustrate this situation by some examples.

***Example 4.1.38:*** Let $S = \{[0, a] \mid a \in Z_5 \setminus \{0\}\} = \{[0, 1], [0, 2], [0, 3], [0, 4]\}$ is an interval group under multiplication modulo five. We see $[0, 4] [0, 4] = [0, 1] [0, 1]$ is the identity element of S.

We will illustrate this by the following table.

| ×      | [0, 1] | [0, 2] | [0, 3] | [0, 4] |
|--------|--------|--------|--------|--------|
| [0, 1] | [0, 1] | [0, 2] | [0, 3] | [0, 4] |
| [0, 2] | [0, 2] | [0, 4] | [0, 1] | [0, 3] |
| [0, 3] | [0, 3] | [0, 1] | [0, 4] | [0, 2] |
| [0, 4] | [0, 4] | [0, 3] | [0, 2] | [0, 1] |

S is an interval abelian group of finite order under multiplication.

***Example 4.1.39:*** Let $S = \{[0, a] \mid a \in Z_7 \setminus \{0\}\}$ be an interval group under multiplication of order 6. Clearly S is of order six and is commutative. [0, 2] [0, 4] = [0, 1], i.e., [0, 2] is the inverse of [0, 4].

[0, 3] [0, 5] = [0, 1] and [0, 6] [0, 6] = [0, 1]. This group is given by the following table.

| × | [0, 1] | [0, 2] | [0, 3] | [0, 4] | [0, 5] | [0, 6] |
|---|--------|--------|--------|--------|--------|--------|
| [0, 1] | [0, 1] | [0, 2] | [0, 3] | [0, 4] | [0, 5] | [0, 6] |
| [0, 2] | [0, 2] | [0, 4] | [0, 6] | [0, 1] | [0, 3] | [0, 5] |
| [0, 3] | [0, 3] | [0, 6] | [0, 2] | [0, 5] | [0, 1] | [0, 4] |
| [0, 4] | [0, 4] | [0, 1] | [0, 5] | [0, 2] | [0, 6] | [0, 3] |
| [0, 5] | [0, 5] | [0, 3] | [0, 1] | [0, 6] | [0, 4] | [0, 2] |
| [0, 6] | [0, 6] | [0, 5] | [0, 4] | [0, 3] | [0, 2] | [0, 1] |

Clearly S is cyclic and generated by <[0, 3]> S is only commutative.
Now we give yet another example of interval group of order 10.

*Example 4.1.40:* Let G = {[0, 1], [0, 2], [0, 3], [0, 4], [0, 5], [0, 6], [0, 7], [0, 8], [0, 9], [0, 10], ×} where × is modulo 11. G is an interval group of order 10 and <[0, 6]> generates G and G is a cyclic interval group of order 10.

Now having seen examples of interval groups of finite order and infinite order now we proceed onto define the interval group semirings.

*Example 4.1.41:* Let S = {$Z^+ \cup \{0\}$} be a semiring with unit, G = {[0, a] | a ∈ $Z_{19} \setminus \{0\}$} be an interval group of order 18. Let

$$SG = \left\{ \sum_i \alpha_i [0, a_i] \,\middle|\, a_i \in Z_{19} \setminus \{0\} \text{ and } \alpha_i \in Z^+ \cup \{0\} \right\}$$

be the interval group semiring of the interval group G over the semiring G. Clearly SG is a semiring. We see SG is not finite but SG is commutative. Further SG is a strict semiring.

**Example 4.1.42:** Let $S = \{R^+ \cup \{0\}\}$ be a semifield. Suppose $G = \{[0, a] / a \in Q^+\}$ be an interval group. Let

$$SG = \left\{\sum_i \alpha_i [0, a_i] \right\}$$

we consider only finite formal sums of the interval elements from the interval group and coefficients from the semifield $R^+ \cup \{0\}\}$, SG is a semiring which is the interval group semiring of the interval group G over the semiring S. SG is commutative of infinite order infact SG is a semifield. SG is also known as the special interval semiring. Suppose $\alpha = 8 [0, 4] + 19 [0, 3/7] + 2 [0, 11/2]$ and $\beta = [0, 3] + 3 [0, 5/9] + 8 [0, 3/7]$ be in SG. Then

$$\begin{aligned}\alpha + \beta &= (8 [0, 4] + 19 [0, 3/7] + 2 [0, 11/2]) + [0, 3] + 3 [0, 5/9] + 8 [0, 3/7] \\ &= (8 [0, 4] + 27 [0, 3/7] + 2 [0, 11/2]) + 1[0, 3] + [0, 5/9].\end{aligned}$$

Now
$$\begin{aligned}\alpha\beta &= (8 [0, 4] + 19 [0, 3/7] + 2 [0, 11/2]) \times (1[0, 3] + 1 [0, 5/9] + 8 [0, 3/7]) \\ &= (8.1 [0, 4] [0, 3] + 19.1 [0, 3/7] [0, 5/9] + 8.1 [0, 4] [0, 5/9] + 19.1 [0, 3/7] [0, 3] + \\ & \quad + 8.8 [0, 4] [0, 5/9] + 19.1 [0, 3/7] [0, 3] + 2.1 [0. 11/2] [0, 3] + 2.1 [0, 11/2] [0, 5/9] + 2.8 [0, 11/2] [0, 3/7] \\ &= 8 [0, 12] + 19 [0, 15/63] + 8 [0, 20/9] + 19 [0, 9/7] + 64 [0, 12/7] + 152 [0, 9/49] + 2 [0, 33/2] + 2 [0, 55/18] + 16 [0, 33/14].\end{aligned}$$

**Example 4.1.43:** Let $S = \{Q^+ \cup \{0\}\}$ be an infinite semifield. $G = \{[0, a] / a \in Z_{11} \setminus \{0\}\}$ be an interval group. Then $SG = \{\Sigma a_i [0, g_i] / a_i \in S \text{ and } g_i \in Z_{11} \setminus \{0\}\}$ is the interval group semiring of the interval group G over the semiring S. Suppose $\alpha = 9 [0,$

3] + 2 [0, 9] + 4 [0, 6] + 3/7 [0, 7] and β = 3 [0, 8] + [0, 10] + 3/2 [0, 6] be in SG. Now

$$\begin{aligned}\alpha + \beta &= (9\,[0, 3] + 2\,[0, 9] + 4\,[0, 6] + 3/7\,[0, 7]) + (3\,[0, 8] \\ &\quad + [0, 10] + 3/2\,[0, 6]) \\ &= 9\,[0, 3] + 2\,[0, 9] + 3/7\,[0, 7] + (4 + 3/2)\,[0, 6] + 3\,[0, 8] + [0, 10].\end{aligned}$$

$$\begin{aligned}\alpha\beta &= (9\,[0, 3] + 2\,[0, 9] + 4\,[0, 6] + 3/7\,[0, 7]) \times (3\,[0, 8] + [0, 10] + 3/2\,[0, 6]) \\ &= 9.3\,[0, 3]\,[0, 8] + 9.1\,[0, 3]\,[0, 10] + 9.3/2\,[0, 3]\,[0, 6] + \\ &\quad 2.3\,[0, 9]\,[0, 8] + 2.1\,[0, 9]\,[0, 10] + 2.3/2\,[0, 9]\,[0, 6] + \\ &\quad 4.3\,[0, 6]\,[0, 8] + 4.1\,[0, 6]\,[0, 10] + 4.3/2\,[0, 6]\,[0, 6] + \\ &\quad 3/7.3\,[0, 7]\,[0, 8] + 3/7.1\,[0, 7]\,[0, 10] + 3/7.\,3/2\,[0, 7]\,[0, 6] \\ &= 27\,[0, 2] + 9\,[0, 8] + 27/2\,[0, 7] + 6\,[0, 6] + 2\,[0, 2] + 3\,[0, 10] + 12\,[0, 4] + 4\,[0, 5] + 6\,[0, 3] \\ &= (27\,[0, 2] + 2\,[0, 2]) + 9\,[0, 8] + 27/2\,[0, 7] + 6\,[0, 6] + 3\,[0, 10] + 12\,[0, 4] + 4\,[0, 5] + 6\,[0, 3] \\ &= 29\,[0, 2] + 9\,[0, 8] + 27/2\,[0, 7] + 6\,[0, 6] + 3\,[0, 10] + 12\,[0, 4] + 4\,[0, 5] + 6\,[0, 3].\end{aligned}$$

SG is a special interval semiring of infinite order.

*Example 4.1.44:* Let $S = \{R^+ \cup \{0\}\}$ be a semiring. $G = \{[0, a] / a \in R^+\}$ be an interval group under multiplication. SG be the interval group semiring of the interval group G over the semiring S. Suppose $\alpha = \sqrt{2}\,[0, 5] + 3\,[0, \sqrt{19}\,] + 2\,[0, 4/3]$ and $\beta = \sqrt{41}\,[0, \sqrt{5}\,] + \sqrt{3}\,[0, 1] + \sqrt{17}\,[0, \sqrt{8}\,]$ be in SG.

Then

$$\begin{aligned}\alpha + \beta &= \sqrt{2}\,[0, 5] + 3\,[0, \sqrt{19}\,] + 2\,[0, 4/3] + \sqrt{41}\,[0, \sqrt{5}\,] \\ &\quad + \sqrt{3}\,[0, 1] + \sqrt{17}\,[0, \sqrt{8}\,].\end{aligned}$$

Now

$\alpha \cdot \beta = \sqrt{2}\ [0, 5] + 3\ [0, \sqrt{19}\ ] + 2\ [0, 4/3]\ (\sqrt{41}\ [0, \sqrt{5}\ ] + \sqrt{3}\ [0, 1] + \sqrt{17}\ [0, \sqrt{8}\ ]$

$= \sqrt{2}.\sqrt{41}\ [0, 5], [0, \sqrt{5}\ ] + \sqrt{2}.\sqrt{3}\ [0, 5]\ [0, 1] + \sqrt{2}.\sqrt{17}\ [0, 5], [0, \sqrt{8}\ ] + 3\ \sqrt{41}\ [0, \sqrt{19}\ ]\ [0, \sqrt{5}\ ] + 3.\ \sqrt{3}\ [0, \sqrt{19}\ ]\ [0, 1] + 3.\ \sqrt{17}\ [0, \sqrt{19}\ ]\ [0, \sqrt{8}\ ] + 2.\ \sqrt{41}\ [0, 4/3], [0, \sqrt{5}\ ] + 2.\ \sqrt{3}\ [0, 4/3]\ [0, 1] + 2.\ \sqrt{17}\ [0, 4/3]\ [0, \sqrt{8}\ ]$

$= \sqrt{82}\ [0, 5\sqrt{5}\ ] + \sqrt{6}\ [0, 5] + \sqrt{34}\ [0, 5\sqrt{8}\ ] + 3\ \sqrt{41}\ [0, \sqrt{95}\ ] + 3\ \sqrt{3}\ [0, \sqrt{19}\ ] + 3\ \sqrt{17}\ [0, \sqrt{152}\ ] + 2\ \sqrt{41}\ \left[0, \dfrac{4\sqrt{5}}{3}\right] + 2\ \sqrt{3}\ [0, 4/3] + 2\ \sqrt{7}\ [0, 4\sqrt{8}/3].$

Clearly SG is an infinite commutative interval group semiring. Clearly SG has no zero divisors. Further SG is a strict interval semiring hence is an interval semifield.

Now we proceed onto give examples of interval group semirings of finite order using distributive finite lattices.

***Example 4.1.45:*** Let S = {0, 1, a, b, c, d, e, f, g} be a finite chain lattice given by the following diagram. That is S is a finite semiring.

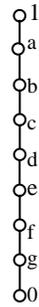

Choose G = {[0, a] | a ∈ $Z_3$ \ {0}} be an interval group. SG = {$a_1$ [0, 1] + $a_2$ [0, 2] / $a_i$ ∈ S 1≤ i ≤ 2} be the interval group semiring of the interval group G over the semiring S. Clearly SG is of finite order. SG has no zero divisor.

*Example 4.1.46:* Let S = {0, a, b, 1} be a semiring, given in example 4.1.34. G = {[0, a] / a ∈ $Z_5$ \ {0}} be an interval group of order four. SG be the interval group semiring of the interval group G over the semiring S. Clearly SG is of finite order SG has zero divisors and units.
Take

a [0, 4] . b [0, 2]     =     a.b [0,4] [0, 2]
                              =     a.b [0, 3]
                              =     0 as a.b = 0.
1 [0, 4] 1 [0, 4]     =     1. [0, 1]
and
1 [0, 2] 1 [0, 3]     =     1 [0, 1].

But SG is a strict commutative interval group semiring which has zero divisor.

Now having seen special types of interval semirings we now proceed onto define interval groupoid semiring, groupoid interval semirings, loop interval semirings and interval loop semirings.

## 4.2 Interval loop semirings

In this section we for the first time introduce non associative interval semirings constructed using groupoids and interval loops. We define these concepts and illustrate them by some simple examples.

**DEFINITION 4.2.1**: *Let G be any groupoid and S be an interval semiring with unit [0,1].*

$$SG = \left\{ \sum_i [0,a]g_i \,\middle|\, [0,a] \in S \text{ and } g_i \in G \right\}$$

*where the sum* $\sum_i [0,a]g_i$ *is a finite formal sum with coefficients from S. The product and sum are defined as in case of semigroup interval semiring. SG is a non associative interval semiring known as the groupoid interval semiring. Thus SG is the groupoid interval semiring of the groupoid G over the interval semiring S.*

We will illustrate this situation by some simple examples.

***Example 4.2.1:*** Let $S = \{[0, a] \,/\, a \in Z^+ \cup \{0\}\}$ be an interval semiring with unit. $G = \{Z_5, *, a * b = ta + ub \pmod 5\ t = 3$ and $u = 2\}$ be a groupoid. $SG = \{\Sigma\ [0, a]\ g \,/\, [0, a] \in S$ and $g \in G\}$, be the groupoid interval semiring of the groupoid G over the interval semiring S. Consider $\alpha = [0, 3]\overline{4} + [0, 8]\overline{3} + [0, 4]\overline{1}$, $\beta = [0, 5]\overline{2} + [0, 7]\overline{3}$ in SG.

$\alpha + \beta\ =\ [0, 3]\overline{4} + ([0, 8] + [0, 7])\overline{3} + [0, 4]\overline{1} + [0, 5]\overline{2}$
$=\ [0, 3]\overline{4} + ([0, 15])\overline{3} + [0, 4]\overline{1} + [0, 5]\overline{2}.$

Now

$\alpha \times \beta\ =\ ([0, 3]\overline{4} + [0, 8]\overline{3} + [0, 4]\overline{1}) \times [0, 5]\overline{2} + [0, 7]\overline{3}$
$=\ [0, 3][0, 5]\overline{4} * \overline{2} + [0, 8][0, 5]\overline{3} * \overline{2} + [0, 4][0, 5]$
$\overline{1} * \overline{2} + [0, 3][0, 7]\overline{4} * \overline{3} + [0, 8][0, 7]\overline{3} * \overline{3} + [0, 4][0, 7]\overline{1} * \overline{3}$
$=\ [0, 15][12+4 \pmod 5) + [0, 40](9 + 4 \pmod 5) + [0, 20](3 + 4 \pmod 5) + [0, 21](12 + 6 \pmod 5) + [0, 56](9 + 6 \pmod 5) + [0, 28](3 + 6 \pmod 5)$
$=\ [0, 15]\overline{1} + [0, 40]\overline{3} + [0, 20]\overline{2} + [0, 21]\overline{3} + [0, 56]\overline{0} + [0, 28]\overline{4}$
$=\ [0, 15]\overline{1} + ([0, 40] + [0, 21])\overline{3} + [0, 20]\overline{2} + + [0, 56]\overline{0} + [0, 28]\overline{4}$
$=\ [0, 15]\overline{1} + [0, 61]\overline{3} + [0, 20]\overline{2} + [0, 56]\overline{0} + [0, 28]\overline{4}.$

We show SG is not associative under the product operation in the groupoid interval semiring.

Consider $\alpha = [0, 2]\,\overline{3}$, $\beta = [0, 5]\,\overline{2}$ and $\gamma = [0, 7]\,\overline{4}$ in SG.
Now

$$
\begin{aligned}
(\alpha\beta)\gamma &= ([0, 2]\,\overline{3} \cdot [0, 5]\,\overline{2}) \times ([0, 7]\,\overline{4}) \\
&= ([0, 2] \cdot [0, 5]\,(9 + 4 \pmod 5))\,[0, 7]\,\overline{4} \\
&= [0, 10]\,\overline{3} \times [0, 7]\,\overline{4} \\
&= [0, 70]\,(9 + 8 \pmod 5) \\
&= [0, 70]\,\overline{2}. \qquad\qquad (I)
\end{aligned}
$$

Consider

$$
\begin{aligned}
\alpha(\beta\gamma) &= [0, 2]\,\overline{3}\,([0, 5]\,\overline{2} \cdot ([0, 7]\,\overline{4}) \\
&= [0, 2]\,\overline{3}\,([0, 5]\,[0, 7]\,(6 + 8 \pmod 5)) \\
&= [0, 2]\,\overline{3}\,([0, 35]\,\overline{4}) \\
&= ([0, 2] \cdot [0, 35]\,(9 + 8 \pmod 5) \\
&= [0, 70]\,\overline{2}. \qquad\qquad (II)
\end{aligned}
$$

I and II are equal hence this set of $\alpha\beta\gamma$ is associative. But in general $\alpha(\beta\gamma) \neq (\alpha\beta)\gamma$. For take $\alpha = [0, 7]\,\overline{4}$, $\beta = [0, 12]\,\overline{2}$ and $\gamma = [0, 10]\,\overline{3}$.
Consider

$$
\begin{aligned}
\alpha(\beta)\gamma &= ([0, 7]\,\overline{4} \cdot [0, 12]\,\overline{2})\,([0, 10]\,\overline{3}) \\
&= [0, 84]\,(12 + 4 \pmod 5)\,([0, 10]\,\overline{3}) \\
&= [0, 84]\,\overline{1} \cdot [0, 10]\,\overline{3} \\
&= [0, 840]\,(3 + 6 \pmod 5) \\
&= [0, 840]\,\overline{4}. \qquad\qquad (I)
\end{aligned}
$$

$$
\begin{aligned}
\alpha(\beta\gamma) &= [0, 7]\,\overline{4},\,([0, 12]\,\overline{2} \cdot [0, 10]\,\overline{3}) \\
&= ([0, 7]\,\overline{4}\,([0, 120]\,(6 + 6 \pmod 5)) \\
&= [0, 7]\,\overline{4} \cdot [0, 120]\,\overline{2} \\
&= [0, 840]\,(12 + 4 \pmod 5) \\
&= [0, 840]\,\overline{1}. \qquad\qquad (II)
\end{aligned}
$$

Clearly I and II are not equal. Thus the operation of product of two elements in the groupoid interval semiring in general is non associative. We have shown for certain triple α, β, γ in SG (αβ) γ = α(βγ) but there exists triples a, b, c in SG with (ab)c ≠ a(bc). Thus SG is a non associative interval semiring of infinite order.

***Example 4.2.2:*** Let $S = \{[0, a] / a \in R^+ \cup \{0\}\}$ be an interval semiring. $G = \{Z_8, *, (3, 4)\}$ be a groupoid. SG be the groupoid interval semiring of the groupoid G over the interval semiring S. $SG = \{\Sigma [0, a_i] g_i \mid g_i \in Z_8 \text{ and } a_i \in R^+ \cup \{0\}\}$. For $\alpha = [0, 5]\overline{3} + [0, 7]\overline{5}$ and $\beta = [0, \sqrt{2}]\overline{6} + [0, \sqrt{5}]\overline{2} + [0, 3]\overline{4}$ in SG

$$\alpha + \beta = [0, 5]\overline{3} + [0, 7]\overline{5} + [0, \sqrt{2}]\overline{6} + [0, \sqrt{5}]\overline{2} + [0, 3]\overline{4}.$$

$$\alpha\beta = ([0, 5]\overline{3} + [0, 7]\overline{5})([0, \sqrt{2}]\overline{6} + [0, \sqrt{5}]\overline{2} + [0, 3]\overline{4})$$

$$= ([0, 5]([0, \sqrt{2}]\overline{3} * \overline{6} + [0, 5][0, \sqrt{5}]\overline{3} * \overline{2} + [0, 5][0, 3]\overline{3} * \overline{4} + [0, 7][0, \sqrt{5}]\overline{5} * \overline{2}$$

$$= [0, 5\sqrt{2}][9 + 24 \pmod 8) + [0, 5\sqrt{5}](9 + 8 \pmod 8) + [0, 15](9 + 16 \pmod 8) + ([0, 7\sqrt{2}](15 + 24 \pmod 8) + ([0, 21](15 + 16 \pmod 8) + ([0, 7\sqrt{5}](15 + 8 \pmod 8)$$

$$= [0, 5\sqrt{2}]\overline{1} + [0, 5\sqrt{5}]\overline{1} + [0, 15]\overline{1} + [0, 7\sqrt{2}]\overline{7} + [0, 21]\overline{7} + [0, 7\sqrt{5}]\overline{7}$$

$$= ([0, 5\sqrt{2}] + [0, 5\sqrt{5}] + [0, 15]\overline{1} + ([0, 7\sqrt{2}] + [0, 21] + [0, 7\sqrt{5}]\overline{7}$$

$$= [0, 5\sqrt{2} + 5\sqrt{5} + 15]\overline{1} + [0, 7\sqrt{2}, 21 + 7\sqrt{5}]\overline{7}.$$

The reader is left with the task of verifying that SG is a non associative groupoid interval semiring.

***Example 4.2.3:*** Let $S = \{[0, a] \mid a \in Q^+ \cup \{0\}\}$ be an interval semiring. Let $G = \{Z^+ \cup \{0\}, *, (3, 5)\}$ be a groupoid of infinite order (for any $a, b \in Z^+ \cup \{0\}$, $a * b = 3a + 5b$; clearly $a * b \in Z^+ \cup \{0\}$).

Consider

$$SG = \left\{ \sum_i [0, a_i] g_i \;\middle|\; a_i \in Q^+ \cup \{0\} \text{ and } g_i \in Z^+ \cup \{0\} \right\}$$

where i run over only finite index; that is all the sums are taken as finite formal sums. SG is the groupoid interval semiring of the groupoid G over the interval semiring S.

It is easily verified SG is a non associative interval semiring of infinite order.

We study a few properties enjoyed by these groupoid interval semirings like ideals and subsemiring. We do not define these concepts the interested reader is requested to refer [12-17].

We only give a few examples related with these concepts.

***Example 4.2.4:*** Let $S = \{[0, a] \mid a \in Z^+ \cup \{0\}\}$ be an interval semiring. $P = \{Z_7, *, (3, 4)\}$ be a groupoid. The groupoid interval semiring SG is a non commutative groupoid interval semiring.

***Example 4.2.5:*** Let $S = \{[0, a] \mid a \in Q^+ \cup \{0\}\}$ be an interval semiring. $G = \{Z_4, *, (2, 3)\}$ be a groupoid. SG the groupoid interval semiring. Choose $W = \{SP \text{ where } P = \{0, 2\} \subseteq Z_4\} \subseteq SG$, W is a left ideal in the groupoid interval semiring. It is easily verified W is not a right ideal of SG. Also $V = \{ST \text{ where } T = \{1, 3\} \subseteq Z_4, S \text{ the semiring}\} \subseteq SG$ is again a left ideal of SG and is not a right ideal of SG.

***Example 4.2.6:*** Let $S = \{[0, a] / a \in Q^+ \cup \{0\}\}$ be an interval semiring. $G = \{Z_4, *, (3, 2)\}$ be a groupoid. SG be the groupoid interval semiring of the groupoid G over the interval semiring S.

Now $P = \{0, 2\} \subseteq Z_4$, SP is a groupoid interval semiring and $SP \subseteq SG$ so SP is a groupoid interval subsemiring of SG. SP is a right ideal and is not a left ideal of SG.

$T = \{1, 3\} \subseteq Z_4$, ST is a groupoid interval semiring and $ST \subseteq SG$ so ST is a groupoid interval subsemiring of SG. ST is a right ideal, both ST and SP are only right ideals of SG and are not left ideal of SG.

***Example 4.2.7:*** Let $G = \{Z_{12}, *, (2, 10)\}$ be a groupoid. $S = \{[0, a] / a \in Q^+ \cup \{0\}\}$ be an interval semiring. SG is the groupoid interval semiring of the groupoid G over the interval semiring S.

$SP = \{\Sigma [0, a] g_i / g_i \in P \{0, 2, 4, 6, 8, 10\} \subseteq Z_{12}, a \in Q^+ \cup \{0\}\} \subseteq SG$; SP is a groupoid interval subsemiring.

Now having seen examples of groupoid interval subsemiring and ideal of groupoid interval semiring we now define S-groupoid interval semiring.

**DEFINITION 4.2.2:** *Let G be any groupoid. S an interval semiring. If G is a S-groupoid then G has a subset P of G where P is a semigroup. SG is a groupoid interval semiring of the groupoid G over the interval semiring S. SP is a semigroup interval semiring and SP is a proper subset of SG. We call SG a S-groupoid interval semiring and G is a S-groupoid.*

We will give examples of such structures.

***Example 4.2.8:*** Let $G = \{Z_{10}, *, (1, 5)\}$ be a groupoid and $S = \{[0, a] / a \in Q^+ \cup \{0\}\}$ be an interval semiring. SG the groupiod interval semiring of the groupoid G over the interval semiring S. Take $P_1 = \{0, 5\} \subseteq Z_{10}$, $P = \{P_1, *, (1, 5)\}$ is a semigroup interval semiring and $SP \subseteq SG$ so SG is a S-groupoid interval semiring.

It is interesting and important to note that all groupoid interval semirings are not S-groupoid interval semirings.

Interested reader is expected to give examples of such structures. Now we get different classes of groupoid interval semirings like P-groupoid interval semiring, right or left alternative groupoid interval semiring, Moufang groupoid interval semiring, idempotent groupoid interval semiring and so on.

We call a groupoid interval semiring SG of a groupoid G over an interval semiring S to be a P-groupoid interval semiring if G is a P-groupoid likewise all other structures. We will give only examples of them.

***Example 4.2.9:*** Let $S = \{[0, a] / a \in Z^+ \cup \{0\}\}$ be an interval semiring. $G = \{Z_{12}, *, (7, 7)\}$ be a P-groupoid then SG the groupoid interval semiring is a P-groupoid interval semiring.

***Example 4.2.10:*** Let $S = \{[0, a] / a \in Q^+ \cup \{0\}\}$ be an interval semiring $G = \{Z_6, *, (3, 3)\}$ be a groupoid. Clearly G is an alternative groupoid. SG be the groupoid interval semiring of the groupoid G over the semiring S. SG is an alternative groupoid interval semiring.

All groupoid interval semirings are not alternative groupoid interval semirings. For take $G = \{Z_p, *, (3, 3); p = 19\}$ be a groupoid. $S = \{[0, a] / a \in Z^+ \cup \{0\}\}$ be an interval semiring SG is a groupoid interval semiring but is not an alternative groupoid interval semiring as G is not an alternative groupoid.

Now we can also have infinite groupoids G and get groupoid interval semirings which we illustrate by some examples.

***Example 4.2.11:*** Let $G = \{Z^+ \cup \{0\}, *, (3, 7)\}$ be a groupoid and $S = \{[0, a] / a \in Q^+ \cup \{0\}\}$ be an interval semiring. SG is an infinite groupoid interval semiring.

***Example 4.2.12:*** Let $G = \{Q^+ \cup \{0\}, *, (4, 7)\}$ be a groupoid. $S = \{[0, a] / a \in R^+ \cup \{0\}\}$ be an interval semiring. SG the groupoid interval semiring is an infinite groupoid interval semiring.

***Example 4.2.13:*** Let $G = \{R^+ \cup \{0\}, *, (3, 11)\}$ be a groupoid. $S = \{[0, a] / a \in R^+ \cup \{0\}\}$ be an interval semiring. SG is the groupoid interval semiring of infinite order.

Now having seen infinite groupoid interval semirings we now proceed onto give examples of interval groupoid semirings of an interval groupoid over a semiring.

We will illustrate this situation by some examples.

However the reader is requested to refer [17] for information about interval groupoid.

***Example 4.2.14:*** Let $G = \{[0, a], *, (2, 3), a \in Z_{11}\}$ be an interval groupoid. For any $[0, a]$ and $[0, b]$ in G we have $[0, a] * [0, b] = [0, 2a+3b \pmod{11}]$ $S = \{Z^+ \cup \{0\}\}$ be a semiring. SG $= \{\Sigma\, x\, [0, a]$ where $x \in S$ and $[0, a] \in G$ denotes the collection of all finite formal sums$\}$; SG under addition and multiplication is a non associative semiring known as the interval groupoid semiring.

For $\alpha = 9\, [0, 8] + 3\, [0, 7] + 2\, [0, 3]$ and $\beta = 5\, [0, 7] + 12\, [0, 3] + [0, 5]$ in SG define $\alpha + \beta = 9\, [0, 8] + 8\, [0, 7] + 14\, [0, 3] + [0, 5]$ in SG.

$\alpha\beta$ = $(9\, [0, 8] + 3\, [0, 7] + 2\, [0, 3]) \times (5\, [0, 7] + 12\, [0, 3] + [0, 5])$

= $9.5\, [0, 8]\, [0, 7] + 3.5\, [0, 7]\, [0, 7] + 2.5\, [0, 3]\, [0, 7] + 9.12\, [0, 8]\, [0, 3] + 12.3\, [0, 7]\, [0, 3] + 2.12\, [0, 3]\, [0, 3] + 9.1\, [0, 8]\, [0, 5] + 3.1\, [0, 7]\, [0, 5] + 2.1\, [0, 3]\, [0, 5]$

= $45\, [0, 16 + 2 \pmod{11}] + 15\, (0, 14 + 21 \pmod{11}] + 10\, [0, 6 + 21 \pmod{11}] + 108\, [0, 16 + 9 \pmod{11}] + 36\, [0, 14+9 \pmod{11}] + 24\, [0, 6 + 9 \pmod{11}] + 9\, [0, 16 + 15$

(mod 11)] + 3 [0, 14 + 15 (mod 11)] + 2 [0, 6 + 15 (mod 11)]

= 45 [0, 4] + 15[0, 2] + 10 [0, 5] + 108 [0, 3] + 9 [0, 9] + 3 [0, 7] + 36 [0, 1] + 24 [0, 4]

= 69 [0, 4] + 36 [0, 1] + 15 [0, 2] + 10 [0, 5] + 9 [0, 9] + 3 [0, 7] + 108 [0, 3] $\in$ SG.

***Example 4.2.15:*** Let G = {[0, a] / a $\in$ $Z_{10}$, *, [0, a] * [0, b] = [0, (1a + 5b) mod 10], (1, 5)} be an interval groupoid. G is a S-interval groupoid. Take S = $Q^+ \cup \{0\}$ to be the semiring. SG = {finite formal sums of the form $\Sigma$ a [0, x] / a $\in$ S and [0, x] $\in$ G} is an interval groupoid semiring of infinite order. $\alpha$ = 8 [0, 3] + 2 [0, 7] and $\beta$ = 3 [0, 7] + 4 [0, 9] + 2 [0, 5] in SG we have

$\alpha + \beta$ = 8 [0, 3] + 2 [0, 7] + 3 [0, 7] + 4 [0, 9] + 2 [0, 5]
  = 8 + 5 [0, 7] + 4 [0, 9] + 2 [0, 5] and
$\alpha\beta$ = (8 [0, 3] + 2 [0, 7]) (3 [0, 7] + 4 [0, 9] + 2 [0, 5])
= 8.3 [0, 3] [0, 7] + 2.3 [0, 7] [0, 7] + 8.4 [0, 3] [0, 9] + 2.4 [0, 7] [0, 9] + 2.2 [0, 7] [0, 5] + 8.2 [0, 3] [0, 5]
= 24 [0, 3 + 35 (mod 10)] + 6 [0, 7 + 35 (mod 10)] + 32 [0, 3 + 45 (mod 10)] + 8 [0, 7 + 45 (mod 10)] + 4 [0, 7 + 25 (mod 10) + 16 [0, 3+25 (mod 10)]
= 24 [0, 8] + 6 [0, 2] + 32 [0, 8] + 8 [0, 2] + 4 [0, 2] + 16 [0, 8]
= 72 [0, 8] + 18 [0, 2].

We call SG the S-interval groupoid interval semiring. This semiring is non associative but has a subsemiring which is associative given by SP = {$\Sigma$ a [0, g] where g $\in$ {0, 5} $\subseteq$ $Z_{10}$} $\subseteq$ SG.

We will call such non associative interval semirings which contains interval subsemiring which are associative as special Smarandache interval groupoid semiring.

***Example 4.2.16:*** Let $S = \{Z^+ \cup \{0\}\}$ be a semiring. $G = \{[0, a] / a \in Z_{12}, *, (3, 9)\}$ be an interval groupoid. G is a Smarandache strong Moufang interval groupoid. Consider SG the interval groupoid semiring of the interval groupoid G over the semiring S. Clearly SG is a special Smarandache interval semiring also SG in this case is a Smarandache strong interval Moufang groupoid semiring. Clearly SG is an infinite non commutative non associative interval semiring.

It is pertinent at this stage to mention that by these methods we can construct infinite non associative interval semirings which are also commutative as well as they have interval subsemirings which are associative. Now we proceed onto give ways to construct groupoid interval semirings of finite order.

To this end we will be using mainly distributive lattices as semirings. So we will be making use of these semirings and construct finite interval groupoid semirings.

***Example 4.2.17:*** Let S = {be the chain lattice of order 20 that is $S = \{0, a_1, a_2, \ldots, a_{18}, 1\}$ $\{0 < a_1 < a_2 < \ldots < a_{18} < 1\}$ be the semiring of finite order. Take $G = \{[0, a], a \in Z_8, *, (2, 6)\}$ be an interval groupoid of finite order.

Now $SG = \{a_1 [0, \overline{1}] + a_2 [0, \overline{2}] + \ldots + a_7 [0, \overline{7}] + a_8 [0, \overline{0}] \mid a_i \in S$ and $[0, a] \in G\}$ be a interval groupoid semiring of the interval groupoid G over the semiring S.

Clearly SG is non commutative non associative and is of finite order. SG is a special Smarandache groupoid semiring. Infact SG has two sided ideal given by $I = \{\Sigma a_i [0, x] / a_i \in S$ and $x \in \{\overline{0}, \overline{2}, \overline{4}, \overline{6}\} \subseteq Z_8\}$.

***Example 4.2.18:*** Let $S = \{C_7 = \{0 < a_1 < a_2 < a_3 < a_4 < a_5 < 1\}\}$ be a semiring of finite order. $G = \{[0, a], a \in Z_6, *, (3, 5)\}$ be an interval groupoid of finite order. $SG = \{\Sigma a_i[0, x] \mid [0, x] \in G; a_i \in S\}$ is the interval groupoid semiring of finite order. SG is a special Smarandache interval groupoid of finite order.

***Example 4.2.19:*** Let $S = \{C_9 = \{0 < a_1 < a_2 < \ldots < a_7 < 1\}\}$ be the semiring $G = \{[0, a], a \in Z_{14}, *, (7, 8)\}$ be an interval groupoid with product defined by for any $[0, a]$, $[0, b]$ in G. $[0, a] * [0, b] = [0, 7a + 8b \pmod{14}]$. SG the interval groupoid semiring is of finite order and is a Smarandache special interval semiring.

Having seen examples of finite non associative interval semirings we now proceed onto define the notion of loop interval semirings and interval loop semirings using loops and interval loops using $Z_n$.

**DEFINITION 4.2.3:** *Let $S = \{[0, a] \,/\, a \in Q^+ \cup \{0\}$ or $Z^+ \cup \{0\}$ or $R^+ \cup \{0\}\}$ be an interval semiring. L be any loop. The loop interval semiring SL consists of all finite formal sums of the form $\sum_i [0,a]m_i$ where $[0, a] \in S$ and $m_i \in L$ with '+' and '.' defined by the following way.*

1.  $\sum_i [0,\alpha_i]m_i = \sum_i [0,\beta_i]m_i$ *if and only if* $[0, \alpha_i] = [0, \beta_i]$ *for all i.*

2.  $\left(\sum_i [0,\alpha_i]m_i\right) + \left(\sum_i [0,\beta_i]m_i\right)$
    
    $= \left(\sum_i ([0,\alpha_i] + [0,\beta_i])m_i\right)$
    
    $= \left(\sum_i ([0,\alpha_i + \beta_i])m_i\right)$

3.  $\left(\sum_i [0,\alpha_i]m_i\right)\left(\sum_i [0,\beta_i]m_i\right) = \left(\sum_j [0,\gamma_j]t_j\right)$
    
    *where* $t_j = m_i n_k$ *and* $\gamma_j = \left(\sum \alpha_i, \beta_k\right)$

4.  $[0, r_i] m_i = m_i [0, r_i]$ *for all* $[0, r_i] \in S$ *and* $m_i \in L$
5.  $[0, r] \Sigma [0, r_i] m_i = \Sigma [0, r] [0, r_i] m_i$
    $= \Sigma [0, rr_i] m_i$ *for all* $[0, r] [0, r_i]$ *in S and* $m_i \in L$.

*SG is a non associative interval semiring with $0 \in S$ as the additive identity and $[0, 1] \in S$ and $e \in L$ (e identity of L).*

*$[0, 1] L \subseteq SL$ and $S.e = S \subseteq SL$.*

We will first illustrate this situation by some examples before we proceed onto discuss other related properties with them.

***Example 4.2.20:*** Let $S = \{[0, a] \mid a \in Z^+ \cup \{0\}\}$ be an interval semiring. Let L be a loop given by the following table.

| 'o' | e | $a_1$ | $a_2$ | $a_3$ | $a_4$ | $a_5$ |
|---|---|---|---|---|---|---|
| e | e | $a_1$ | $a_2$ | $a_3$ | $a_4$ | $a_5$ |
| $a_1$ | $a_1$ | e | $a_3$ | $a_5$ | $a_2$ | $a_4$ |
| $a_2$ | $a_2$ | $a_5$ | e | $a_4$ | $a_1$ | $a_3$ |
| $a_3$ | $a_3$ | $a_4$ | $a_1$ | e | $a_5$ | $a_2$ |
| $a_4$ | $a_4$ | $a_3$ | $a_5$ | $a_2$ | e | $a_1$ |
| $a_5$ | $a_5$ | $a_2$ | $a_4$ | $a_1$ | $a_3$ | e |

$SL = \{\sum [0, x] a_i / x \in Z^+ \cup \{0\}; a_i \in L, 1 \leq i \leq 5\}$ is the loop interval semiring of the loop L over the interval semiring S. Clearly SL is of infinite order and SL is a non associative non commutative interval semiring.

***Example 4.2.21:*** Let $S = \{[0, a] / a \in Q^+ \cup \{0\}\}$ be the interval semiring. $L = \{e, a, b, c, d, g\}$ be the loop given by the following table.

| 'o' | e | a | b | c | d | g |
|---|---|---|---|---|---|---|
| e | e | a | b | c | d | g |
| a | a | e | d | b | g | c |
| b | b | d | e | g | c | a |
| c | c | b | g | e | a | d |
| d | d | g | c | a | e | b |
| g | g | c | a | d | b | e |

$SL = \{\sum [0, x] y \;/\; y \in L \text{ and } x \in Q^+ \cup \{0\}\}$ be the loop interval semiring of the loop L over the interval semiring S. Clearly SL is an infinite non associative interval semiring. SL is commutative interval semiring.

*Example 4.2.22:* Let $S = \{[0, a] \;/\; a \in R^+ \cup \{0\}\}$ be an interval semiring. Let $L = \{e, g_1, g_2, \ldots, g_7\}$ be a loop given by the following table.

| 0 | e | $g_1$ | $g_2$ | $g_3$ | $g_4$ | $g_5$ | $g_6$ | $g_7$ |
|---|---|---|---|---|---|---|---|---|
| e | e | $g_1$ | $g_2$ | $g_3$ | $g_4$ | $g_5$ | $g_6$ | $g_7$ |
| $g_1$ | $g_1$ | e | $g_5$ | $g_2$ | $g_6$ | $g_3$ | $g_7$ | $g_4$ |
| $g_2$ | $g_2$ | $g_5$ | e | $g_6$ | $g_3$ | $g_7$ | $g_4$ | $g_1$ |
| $g_3$ | $g_3$ | $g_2$ | $g_6$ | e | $g_7$ | $g_4$ | $g_1$ | $g_5$ |
| $g_4$ | $g_4$ | $g_6$ | $g_3$ | $g_7$ | e | $g_1$ | $g_5$ | $g_2$ |
| $g_5$ | $g_5$ | $g_3$ | $g_7$ | $g_4$ | $g_1$ | e | $g_2$ | $g_6$ |
| $g_6$ | $g_6$ | $g_7$ | $g_4$ | $g_1$ | $g_5$ | $g_2$ | e | $g_3$ |
| $g_7$ | $g_7$ | $g_4$ | $g_1$ | $g_5$ | $g_2$ | $g_6$ | $g_3$ | e |

Consider

$$SL = \left\{ \sum_{i=0}^{7} [0,a] g_i \;/\; a \in R^+ \cup \{0\} \text{ and } g_i \in L, 1 \leq i \leq 7 \text{ or } g_i = e \right\}$$

be the loop interval semiring of loop L over the interval semiring S. Clealry SL is a non associative, commutative infinite interval semiring.

We can define substructures like subsemirings and ideals of the loop interval semiring. The definition is routine and is left as an exercise for the reader; however we give a few examples of them.

*Example 4.2.23:* Let $S = \{[0, a] \;|\; a \in Q^+ \cup \{0\}\}$ be an interval semiring and L be a loop given by the following table.

| 'o' | e | a | b | c | d |
|---|---|---|---|---|---|
| e | e | a | b | c | d |
| a | a | e | c | d | b |
| b | b | d | a | e | c |
| c | c | b | d | a | e |
| d | d | c | e | b | a |

L is clearly non commutative. Consider $SL = \{\Sigma[0, x] g_i \mid x \in Q^+ \cup \{0\}$ and $g_i \in L\}$ be the loop interval semiring of the loop L over the interval semiring S. Take $P = \{\Sigma[0, x] g_i \mid x \in Z^+ \cup \{0\}$ and $g_i \in L\} \subseteq SL$ be a subset of SL. Clearly P is an interval subsemiring of SL of infinite order.

*Example 4.2.25:* Let $S = \{[0, a] \mid a \in Z^+ \cup \{0\}\}$ be an interval semiring. $L_7(3) = \{e, 1, 2, 3, 4, 5, 6, 7\}$ be a loop given by the following table.

| 'o' | e | 1 | 2 | 3 | 4 | 5 | 6 | 7 |
|---|---|---|---|---|---|---|---|---|
| e | e | 1 | 2 | 3 | 4 | 5 | 6 | 7 |
| 1 | 1 | e | 4 | 7 | 3 | 6 | 2 | 5 |
| 2 | 2 | 6 | e | 5 | 1 | 4 | 7 | 3 |
| 3 | 3 | 4 | 7 | e | 6 | 2 | 5 | 1 |
| 4 | 4 | 2 | 5 | 1 | e | 7 | 3 | 6 |
| 5 | 5 | 7 | 3 | 6 | 2 | e | 1 | 4 |
| 6 | 6 | 5 | 1 | 4 | 7 | 3 | e | 2 |
| 7 | 7 | 3 | 6 | 2 | 5 | 1 | 4 | e |

$SL_7(3)$ be the loop interval semiring of the loop $L_7(3)$ over the interval semiring S. Suppose $x = [0, 5] 2 + [0, 7] 3 + [0, 20] 5 + [0, 10]e$ and $y = [0, 3]e + [0, 4] 5$ be in $SL_7(3)$ we find
   $x + y = ([0, 5) 2 + [0, 7] 3 + [0, 24] 5 + [0, 10] e + [0, 3]6$;
clearly $x + y \in SL_7(3)$.
   $x.y = ([0, 5]2 + [0, 7]3 + [0, 20]5 + [0, 10]e) \times$
          $([0, 3]6 + [0, 4]5).$

The product is performed using the table given for $L_7(3)$.

$$\begin{aligned} x.y &= [0, 5] [0, 3] 2.6 + [0, 7] [0, 3] 3.6 + [0, 20] [0, 3] 5.6 + \\ &\quad [0, 10] [0, 3] e.6 + [0, 5] [0, 4] 2.5 + [0, 7] [0, 4] 3.5 + \\ &\quad [0, 20] [0, 4] 5.5 + [0, 10] [0, 4] e.5 \\ &= [0, 15] 7 + [0, 21] 5 + [0, 60] 1 + [0, 30] 6 + [0, 20] 4 + \\ &\quad [0, 28] 2 + [0, 80] e + [0, 40] 5 \\ &= [0, 15] 7 + [0, 61] 5 + [0, 60] 1 + [0, 30] 6 + [0, 20] 4 + \\ &\quad [0, 28] 2 + [0, 80] e. \end{aligned}$$

Thus having seen how multiplication in $SL_7(3)$ is performed now we proceed onto define substructures.

Consider $P = \{\sum [0, x] t_i / x \in Z^+ \cup \{0\}$ and $t \in \{e, 5\} \subseteq L_7(3)\} \subseteq SL_7(3)$; it is easily verified P is a loop interval subsemiring of $SL_7(3)$.

Also P is a strict interval subsemiring. Thus P is an interval semifield contained in $SL_7(3)$, so $SL_7(3)$ is a Smarandache loop interval semiring. $SL_7(3)$ has atleast 7 distinct interval semifields.

**THEOREM 4.2.1:** *Let $S = \{[0, a] / a \in Q^+ \cup \{0\}$ or $R^+ \cup \{0\}$ or $Z^+ \cup \{0\}\}$ be an interval semiring. $L_n(m) = \{e, 1, 2, …, n\}$ be a set where $n > 3$, n odd and m is a positive integer such that $(m, n) = 1$ and $(m-1, n) = 1$ with $m < n$.*

*Define on $L_n(m)$ a binary operation * as follows.*
1. *$e * i = i * e = i$ for all $i \in L_n(m)$*
2. *$i^2 = i * i = e$ for all $i \in L_n(m)$*
3. *$i * j = t$ where $t = (mj - (m-1)i) \pmod{n}$*

*for all $i, j \in L_n(m)$ $i \neq j$, $i \neq e$ and $j \neq e$, then $L_n(m)$ is a loop under the binary operation *, $SL_n(m)$ is a S-semiring.*

The proof is direct and hence is left as an exercise for the reader to prove.

Further we have n number of interval semifields contained in $SL_n(m)$. This method gives us interval semifields different from $Q^+ \cup \{0\}$ or $Z^+ \cup \{0\}$ or $R^+ \cup \{0\}$ but of infinite order containing any one of these subsemifields.

**THEOREM 4.2.2:** *Let $L_n$ (n+1/2) be a loop. $S = \{[0, a] / a \in Z^+ \cup \{0\}$ or $Q^+ \cup \{0\}$ or $R^+ \cup \{0\}\}$ be an interval semiring. The loop interval semiring $SL_n$ (n+1/2) is a commutative interval semiring which is an interval semifield of infinite order.*

The proof is direct and follows from the fact $L_n$ (m) is commutative when m = n+1/2 for more refer [14].

For the first time we proceed onto define interval loops using $Z_n$ (n > 3; n odd).

**DEFINITION 4.2.4:** *Let $G = \{[0, a] / a \in L_n(m); L_n(m)$ the loop built using $Z_n$, n > 3, n odd, with m < n (m, n) = 1 (m – 1, n) = 1\} G under the operations of $L_n(m)$ is an interval loop.*
*For if $[0, a], [0, b] \in G$ then*
*$[0, a] * [0, b] = [0, a * b] = [0, (mb – (m – 1) a) (mod n)]$*
*(if $a \neq e$, $b \neq a$ and $b \neq e$); if $a = b$ then $[0, a] [0, a] = [0, e]$ for $[0, a]$. $[0, e] = [0, e] [0, a] = [0, a]$ for all $[0, a] \in G$.*

It is left as an exercise for the reader to prove G is a loop and is called an interval loop.
We will proceed onto give examples of them.

*Example 4.1.25 :* Let $G = \{[0, a] \mid a \in L_5(2)\}$ G is an interval loop of order six given by the following table.

| * | [0, e] | [0, 1] | [0, 2] | [0, 3] | [0, 4] | [0, 5] |
|---|---|---|---|---|---|---|
| [0, e] | [0, e] | [0, 1] | [0, 2] | [0, 3] | [0, 4] | [0, 5] |
| [0, 1] | [0, 1] | [0, e] | [0, 3] | [0, 5] | [0, 2] | [0, 4] |
| [0, 2] | [0, 2] | [0, 5] | [0, e] | [0, 4] | [0, 1] | [0, 3] |
| [0, 3] | [0, 3] | [0, 4] | [0, 1] | [0, e] | [0, 5] | [0, 2] |
| [0, 4] | [0, 4] | [0, 3] | [0, 5] | [0, 2] | [0, e] | [0, 1] |
| [0, 5] | [0, 5] | [0, 2] | [0, 4] | [0, 1] | [0, 3] | [0, e] |

We see any [0, 3] * [0, 4] = [0, 5] got from the table.

We can have infinite number of interval loops constructed in this way where n > 3 and n odd is the only criteria. All these interval loops are of finite order and their order is even.

Further the reader is given the task of proving the following theorem.

**THEOREM 4.2.3**: *Let*

$$G = \{[0, a] \mid a \in L_n\left(\frac{n+1}{2}\right)\}$$

*be an interval loop, G is a commutative interval loop.*

For more information about loops please refer [14].

***Example 4.1.26:*** Let $G = \{[0, a] \mid a \in L_7(4)$, * operation on $L_7(4)\}$ be an interval loop. It can be easily verified $L_7(4)$ is a commutative loop of order 8.

We see several nice properties enjoyed by these interval loops. These interval loops constructed by this way are not Moufang interval loops. They are also not interval Bol loop or interval Bruck loop. However we have interval loop which satisfy the weak inverse property. That is we have WIP interval loops. The following theorem the proof of which is direct gurantees the existence of WIP interval loops.

**THEOREM 4.2.4:** *Let $(G, *) = \{[0, a] \mid a \in L_n(m), *\}$ be an interval loop. $L_n(m)$ is a weak inverse property loop if and only if $(m^2 - m + 1) \equiv 0 \pmod{n}$.*

We will give an example of a interval WIP loop or WIP interval loop.

***Example 4.1.27:*** Let $G = \{[0, a] \mid a \in L_7(3) = \{e, 1, 2, 3, 4, 5, 6,$ *\}$ be the interval loop. It is easily verified G is an WIP interval loop.

However these class of interval loops contains an left alternative interval loop and right alternative interval loop. This is evident from the following theorem, the proof of which is left as an exercise to the reader. For more information refer [14].

**THEOREM 4.2.5:** *Let $G = \{[0, a] / a \in (L_n(2), *)\}$ be an interval loop. G is a left alternative interval loop.*

**THEOREM 4.2.6:** *Let $G = \{[0, a] / a \in (L_n(n - 1), *)\}$ be an interval loop. G is a right alternative interval loop.*

**THEOREM 4.2.7:** *No interval loop G constructed using $L_n(m)$ is alternative.*

*Hint:* Both the left alternative interval loop G or the right alternative interval loop G are built using $L_n(2)$ and $L_n(n - 1)$ respectively and both of them are non commutative.

We can define associator of an interval loop exactly as in case of the loops. The following theorem is also direct and hence is left as an exercise for the reader.

**THEOREM 4.2.8:** *Let $G = \{[0, a] / a \in (L_n(m), *)\}$ be an interval loop built using $L_n(m)$. The associator $A(G) = G$.*

The proof directly follows from the fact $A(L_n(m)) = L_n(m)$. Also these interval loops built using $L_n(m)$ are G-interval loops [14].

**THEOREM 4.2.9:** *All interval loops $G = \{[0, a] \mid a \in L_n(m), *\}$ are Smarandache interval loops (S-interval loops).*

Follows from the fact that all loops $L_n(m)$ are S-loops.

**THEOREM 4.2.10:** *The interval loops $G = \{[0, a] \mid a \in (L_n(m), *)\}$ are Smarandache simple interval loops (S-simple interval loops).*

We can define Smarandache Lagrange interval loop (S-Lagrange interval loop) as in case of general loops [14]. Also the notion of Smarandache weakly Lagrange interval loop can be defined [14]. This task is left as an exercise to the reader.

We will however give examples of these and results connected with them.

***Example 4.2.28:*** Let G = {[0, a] | a ∈ $L_{15}$ (2)} be an interval loop. G is a S-Lagrange interval loop as G has only interval subgroups of order 2 and 4.

We give the following theorems and the proofs of which are direct and the reader is expected to prove.

**THEOREM 4.2.11:** *Every S Lagrange interval loop is a S-weakly Lagrange interval loop.*

**THEOREM 4.2.12:** *Let G = {[0, a] | a ∈ $L_n$ (m) where n is a prime} be an interval loop. Then G is a S-Lagrange interval loop for all m such that (m-1, n) = 1 and (n, m) = 1, m < n.*

**THEOREM 4.2.13:** *Let G = {[0, a] | a ∈ $L_n$ (m)} be an interval loop. Every interval loop G is a S-weakly Lagrange loop.*

We say an interval loop L of finite order is a Smarandache pseudo Lagrange interval loop if order of every S-interval subloop divides order of L. If the order of atleast one S-interval subloop of L divides the order of the interval loop L then we call L to be a S-weakly pseudo Lagrange interval loop.

We can as in the case of finite loops define the notion of Smarandache p-Sylow interval subloops and Smarandache p-Sylow interval subgroups [14]. Also the notion of Smarandache strong p-Sylow loop.

The following theorem is direct and hence the proof is left as an exercise for the reader.

**THEOREM 4.2.14:** *Let G = {[0, a] | a ∈ $L_n(m)$} be an interval loop where n is a prime. Then G is a S-strong 2-Sylow interval loop for all m < n.*

$G = \{[0, a] \mid a \in L_n (m)\}$ *gives an interval loop for every* $1 < m < n$, *as n is given to be a prime. Now every interval loop constructed in this manner are Smarandache strong 2-Sylow loop as every element in G is of order two and $2/p+1$ where $p = n$ is a prime.*

We will illustrate this by examples.

***Example 4.2.29:*** Let $G = \{[0, a] \mid a \in L_{11} (m); 1 < m < 11\}$ be an interval loop. We can construct 9 distinct interval loops for varying m, $m = 2, 3, \ldots, 10$. It is easily verified every $x = [0, a]$ in G is such that $x^2 = [0,e]$ hence G is a S-strong 2-Sylow interval loop.

Infact we have infinite class of S-strong 2-Sylow interval loops using $L_n (n)$ where n is a prime; as the number of primes is infinite. We can define as in case of S loops for S-interval loops the notion of Smarandache interval loop homomorphism [14].

Now with these concepts of interval loops now we can define the new notion of interval loop semirings and interval loop interval semirings.

**DEFINITION 4.2.5:** *Let $G = \{[0, a] \mid a \in L_n (m)\}$ be an interval loop. $S = Z^+ \cup \{0\}$ be a semiring.*

$$SG = \left\{ \sum_i a_i [0, x_i] \middle| i \right.,$$

*runs over finite formal sums and $a_i \in S$, $x_i \in L_n (m)\}$ is an interval semiring under usual addition and multiplication, defined as the interval loop interval semiring.*

***Example 4.2.30:*** Let $S = Q^+ \cup \{0\}$ be a semiring. $L = \{[0, a] \mid a \in L_7 (3)\}$ be an interval loop. SL the interval loop semiring is an interval semiring.

***Example 4.2.31:*** Let $S = Z^+ \cup \{0\}$ be a semiring. $L = \{[0, a] \mid a \in L_{13} (5)\}$ be an interval loop. SL the interval loop semiring where

$$SL = \left\{ \sum_i a_i[0, x_i] \,\middle|\, x_i \in L_{13}(5) \text{ and } a_i \in S \right\}$$

is of infinite order. Clearly this interval semiring is non associative and non commutative.

**Example 4.2.32:** Let $G = \{[0, a] \mid a \in L_{15}(8)\}$ be an interval loop and $S = R^+ \cup \{0\}$ be a semiring.

$$SG = \left\{ \sum_i a_i[0, x_i] \,\middle|\, a_i \in S \text{ and } x_i \in L_{15}(8) \right\}$$

be the interval loop semiring of the interval loop over the semiring S. SG is of infinite order. SG is a commutative interval loop semiring, which is a non associative interval loop semiring.

**Example 4.2.33:** Let $G = \{[0, a] \mid a \in L_{17}(8)\}$ be an interval loop and $S = Z^+ \cup \{0\}$ be a semiring.

$$SG = \left\{ \sum_i a_i[0, x_i] \,\middle|\, a_i \in S \text{ and } x_i \in G \right\}$$

is the interval loop semiring of the interval loop G over the semiring S. SG is non commutative and is of infinite order. $1.[0, e]$ is the identity element of SG.

We can construct 15 such interval loop semirings using the interval loops $L_{17}(2)$, $L_{17}(3)$, $L_{17}(4)$, …, $L_{17}(16)$ and the semiring $G = Z^+ \cup \{0\}$. Depending on m the interval loop semirings enjoy different properties.

We will find the substructures like ideals and subsemirings in them.

**Example 4.2.34:** Let $S = \{Q^+ \cup \{0\}\}$ be a semiring. $G = \{[0, a] \mid a \in L_9(7)\}$ be an interval loop. SG the loop interval semiring of the interval loop G over the semiring S.
Consider

$$P = \left\{ \sum_i a_i[0, x_i] \,\middle|\, a \in Z^+ \cup \{0\} \text{ and } x_i \in L_9(7) \right\} \subseteq SG;$$

P is a loop interval subsemiring of the interval loop semiring SG. Clearly P is not an ideal of SG.

Consider

$$T = \left\{ \sum_i a_i[0, x_i] \,\middle|\, x_i = 8, a \in Q^+ \cup \{0\} \right\} \subseteq SG,$$

T is an interval group subsemiring which is also a semifield.

***Example 4.2.35:*** Let $G = \{[0, a] \mid a \in L_{19}(10)\}$ be a interval loop and $S = \{Z^+ \cup \{0\}\}$ be a semiring;

$$SG = \left\{ \sum_i a_i[0, x_i] \,\middle|\, a \in Z^+ \cup \{0\} \text{ and } x_i \in L_{19}(10) \right\}$$

be the interval loop semiring of the interval loop G over the semiring S.

Clearly SG is a commutative interval semiring. Consider

$$P = \left\{ \sum_i x_i[0, a] \,\middle|\, x_i \in 5Z^+ \cup \{0\} \text{ and } a \in L_{19}(10) \right\} \subseteq SG;$$

P is an interval loop subsemiring of SG as well as an ideal of SG.

***Example 4.2.36:*** Let $G = \{[0, a] \mid a \in L_{41}(21)\}$ be a loop interval and $S = \{Z^+ \cup \{0\}\}$ be a semiring. SG be the interval loop semiring of G over the semiring S.

$$SG = \left\{ \sum_i x_i[0, a_i] \,\middle|\, x_i \in Z^+ \cup \{0\} \text{ and } a \in L_{41}(21) \right\}.$$

Consider

$$P = \left\{ \sum_i a_i[0, x_i] \,\middle|\, a_i \in 5Z^+ \cup \{0\} \text{ and } x_i \in L_{19}(10) \right\} \subseteq SG;$$

P is an interval loop subsemiring of SG as well as an ideal of SG.

*Example 4.2.37:* Let $G = \{[0, a] \mid a \in L_{23}(12)\}$ be an interval loop. $S = Z^+ \cup \{0\}$ be the semiring.

$$SG = \left\{ \sum x_i[0, a_i] \,\middle|\, x_i \in S \text{ and } a_i \in L_{23}(12) \right\}$$

be the interval loop semiring of the interval loop G over the semiring S.

Consider

$$P = \left\{ \sum x_i[0, a_i] \,\middle|\, x_i \in 9Z^+ \cup \{0\}, a_i \in L_{23}(12) \right\} \subseteq SG;$$

P is an ideal of SG for SG is a commutative interval semiring.

We can as in case of loops build in case of interval loops also the notion of first normalizer and second normalizer [14].

First normalizer of $G = \{[0, a] \mid a \in L_n(m)\}$ is given by;

$N_1(H_1(t)) = \{[0, a] \in L_n(m)\} \mid [0, a] H_i(t) = (H_i(t) [0, a]\}$ and second normalizer of G is $N_2(H_i(t)) = \{[0, x] \in L_n(m) \mid [0, x] (H_i(t)) [0, x] = H_i(t)\}$.

We illustrate this situation by an example.

*Example 4.2.38:* Let $G = \{[0, a] \mid a \in L_{25}(2)\}$. Consider $H_i(3) = \{[0, e], [0, 1], [0, 4], [0, 7], [0, 10], [0, 13]\} \subseteq G$, $H_i(3)$ is an interval subloop of G.

Further $N_1(H_1(3)) = G$.

*Example 4.2.39:* Let $G = \{[0, a] \mid a \in L_{45}(8)\}$ be an interval loop $H_1(15) = \{[0,a] \mid a \in L_{45}(8)\} \subseteq G$, but $N_2(H_1(15)) =$ $H_1(15)$. So the two normalizers of the interval subloops are not always equal. Now using interval subloops of G also we can construct in interval loop subsemirings and ideals. This take is left as an exercise to the reader.

***Example 4.2.40:*** Let $G = \{[0,a] \,/\, a \in L_{21}(20)\}$ be $S = \{C_9; 0 < a_1 < a_2 < \ldots < a_7 < 1\}$ be a chain lattice which is a semiring.

$$SG = \left\{\sum x_i[0,a_i] \,\Big|\, a_i \in S,\, x_i \in L_{21}(20)\right\}$$

be the interval loop semiring. Clearly SG is of finite order; SG has non trivial idempotents but no zero divisors.

***Example 4.2.41:*** Let $G = \{[0, a] \mid a \in L_{25}(8)\}$ be an interval loop.

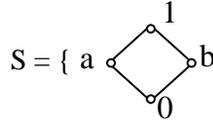

$S = \{ a,\ b,\ 0,\ 1$

$a.b = 0,\ a+b = 1,\ a.0 = 0 = b0,\ 1+a = 1+b = 1,\ 1.a = a,\ b.1 = b\}$ be a semiring.

$$SG = \left\{\sum x_i[0,a_i] \,\Big|\, a_i \in S,\, x_i \in L_{25}(8)\right\}$$

be an interval loop semiring. SG has zero divisors idempotents.
For take $x = a\,[0, 9]$ and $y = b\,[0, 12]$ ;
$x.y =\ \ a.b\,([0, 9] \cdot [0.12])$
$\ \ \ \ =\ \ 0.$
Thus $x, y \in SG$ is a zero divisor.
$x^2\ =\ \ ([0, e] + [0, 9])^2$
$\ \ \ =\ \ [0, e] + [0, 9] + [0, e] + [0, 9]$
$\ \ \ =\ \ [0, e] + [0, 9] = x.$

Thus SG has non trivial idempotents.
Further $x = 1.\,[0, t];\ t \in L_{25}(8)$ are units in SG.

*Example 4.2.42:* Let $G = \{[0, a] \mid a \in L_7(4)\}$ be an interval loop.

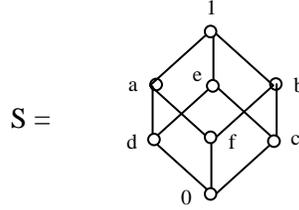

$S =$

/ $de = ef = 0$, $a+c = b+d = 1$ and so on$\}$ be a semiring. SG be the interval loop semiring of the interval loop G over the semiring S. SG is of finite order.

$x = [0, e] + [0, 1] + [0, 2] + [0, 3] + [0, 4] + [0, 5] + [0, 6] \in SG$
$x^2 = ([0, e] + [0, 1] + [0, 2] + [0, 3] + [0, 4] + [0, 5] + [0, 6])^2$
$\phantom{x^2} = [0, e] + [0, 1] + [0, 2] + [0, 3] + [0, 4] + [0, 5] + [0, 6]$
$\phantom{x^2} = x.$

As $1 + 1 = 1$ and $1.1 = 1$ in S and $[0, a][0, a] = [0, e]$ in G.

Thus SG has idempotents. It is easily verified SG has zero divisors also.

*Example 4.2.43:* Let $G = \{[0, a] \mid a \in L_{33}(7)\}$ be an interval loop. $S = \{C_2 = 0 < 1\}$ be a semiring. SG is a division semiring and is non commutative have no zero divisors but has idempotents.

*Example 4.2.44:* Let $G = \{[0, a] \mid a \in L_{23}(12)\}$ be an interval loop.

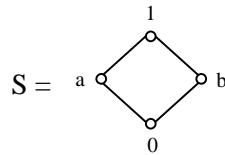

$S =$

be the semiring. SG the interval loop semiring. SG has zero divisor, idempotents but commutative and is not semifield.

*Example 4.2.45:* Let $G = \{[0, a] \mid a \in L_{13}(7)\}$ be an interval loop $S = \{C_2 = 0, 1; 0 < 1\}$ be the semiring SG is a semifield.

In view of this we have the following interesting theorem.

**THEOREM 4.2.15:** *Let $G = \{[0, a] \mid a \in L_p\left(\frac{p+1}{2}\right), p$ a prime$\}$ be the interval loop. $S = \{C_n, 0 < a_1 < \ldots < a_{n-2} < 1\}$ be the semiring. SG the interval loop semiring is a semifield.*

The proof is left as an exercise for the reader.

*Example 4.2.46:* Let $G = \{[0, a] \mid a \in L_n(m)\}$ be an interval loop $m \neq n+1/2$. $S = \{C_n; 0 < a_1 < \ldots < a_{n-1} < 1\}$ be the semiring. SG be the interval loop semiring. SG is a division semiring; that is SG has no zero divisors and is strict.

In view of this we have the following theorem.

**THEOREM 4.2.16:** *Let $G = \{[0, a] \mid a \in L_n(m); m \neq n+1/2\}$ be an interval loop. $S = \{C_d; 0 < a_1 < \ldots < a_{d-1} < 1\}$ be a semiring. The interval loop semiring SG is an interval division semiring.*

The proof is direct and is left as an exercise for the reader to prove.

It is important to mention here that all interval loop, $G = \{[0, a] \mid a \in L_n(m)\}$ are S-simple. Hence to find ideals in them are not very easy that too when the semiring is $Q^+ \cup \{0\}$ or $R^+ \cup \{0\}$. Now we proceed on to define interval subsemiring and ideals of finite interval semirings.

*Example 4.2.47:* Let $G = \{[0, a] \mid a \in L_{17}(15)\}$ be an interval loop. $S = \{C_5, 0 < a_1 < a_2 < a_3 < 1\}$ be the semiring. SG the interval loop semiring of the interval loop G over the semiring S.

$$SG = \left\{\sum a_i[0, x_i] \mid x_i \in L_{17}(15) \text{ and } a_i \in S\right\}.$$

Consider

$$P = \left\{\sum a[0, x_i] \mid a \in \{0,1\} \text{ and } x_i \in L_{17}(15)\right\} \subseteq SG;$$

P is an interval loop subsemiring of SG but clearly P is not an ideal of SG.

*Example 4.2.48:* Let $G = \{[0, a] \mid a \in L_{21}(11)\}$ be interval loop. $S = \{0 < a_1 < a_2 < a_3 < a_4 < a_5 < a_6 < 1\}$ be the semiring. SG the interval loop semiring. SG has several interval loop subsemirings.

Now we will proceed onto define interval loop interval semiring.

**DEFINITION 4.2.17:** *Let $G = \{[0, a] \mid a \in L_n(m)\}$ be an interval loop. $S = \{[0, x] \mid x \in Z^+ \cup \{0\} \text{ or } Q^+ \cup \{0\} \text{ or } R^+ \cup \{0\}\}$ be an interval semiring. SG the interval loop interval semiring as*

$$SG = \left\{ \sum [0,x][0,a] \;\middle|\; [0,x] \in S \right\}$$

*and $[0, a] \in G$ where the sums are finite formal sums with addition and product respecting S and G} is a semiring.*

*For instance if $x = [0, a_1][0, g_1] + [0, a_2][0, g_2] + [0, a_3][0, g_3]$ and $y = [0, b_1][0, h_1] + [0, b_2][0, h_2])$ are in SG then*

$x.y \;\; = \;\; ([0, a_1][0, g_1] + [0, a_2][0, g_2] + [0, a_3][0, g_3])$
$\qquad\quad ([0, b_1][0, h_1] + [0, b_2][0, h_2])$
$\quad = \;\; [0, a_1][0, b_1][0, g_1][0, h_1] + [0, a_1][0, b_2][0, g_1][0, h_2] + [0, a_2][0, b_1][0, g_1][0, h_1] + \ldots + [0, a_3][0, b_2][0, g_3][0, h_2]$
$\quad = \;\; [0, a_1 b_1][0, b_1 h_1] + [0, a_1 b_2][0, g_1 h_2] + \ldots + [0, a_3 b_2][0, g_3 h_2]$.

We will illustrate this situation by some examples.

*Example 4.2.49:* Let $G = \{[0, a] \mid a \in L_7(4)\}$ be an interval loop. $S = \{[0, x] \mid x \in Z^+ \cup \{0\}\}$ be an interval semiring. $SG = \{\Sigma [0, x][0, a]$ with $a \in L_7(4)$ and $x \in Z^+ \cup \{0\}\}$ be the interval loop interval semiring.

***Example 4.2.50:*** Let $G = \{[0, a] \mid a \in L_{11}(3)\}$ be an interval loop and $S = \{[0, x] \mid x \in Q^+ \cup \{0\}\}$ be an interval semiring. SG be the interval loop interval semiring.

***Example 4.2.51:*** Let $G = \{[0, a] \mid a \in L_{13}(8)\}$ and $S = \{[0, x] \mid x \in R^+ \cup \{0\}\}$ be the interval loop and interval semiring respectively. SG is the interval loop interval semiring.

We will now give examples of substructures of these interval loop interval semirings which are of infinite order.

***Example 4.2.52:*** Let $G = \{[0, a] \mid a \in L_{19}(7)\}$ and $S = \{[0, x] \mid x \in Z^+ \cup \{0\}\}$ be interval loop and interval semiring respectively. $SG = \{\Sigma\, [0, x]\, [0, a] \mid x \in Z^+ \cup \{0\}$ and $a \in L_{19}(7)\}$ be the interval loop interval semiring of the interval loop G over the interval semiring S. Take $P = \{\Sigma\, [0, x]\, [0, a] \mid x \in 13Z^+ \cup \{0\}$ and $a \in L_{19}(7)\} \subseteq SG$; P is an interval loop interval subsemiring; infact P is also an ideal of SG.

***Example 4.2.53:*** Let $G = \{[0, a] \mid a \in L_{27}(13)\}$ be an interval loop and $S = \{[0, x] \mid x \in Q^+ \cup \{0\}\}$ be the interval semiring. $SG = \{$finite formal sums of the form $\Sigma\, [0, x]\, [0, a] \mid a \in L_{27}(13)$ and $x \in Q^+ \cup \{0\}\}$ be the interval loop interval semiring. Let $H = \{\Sigma\, [0, x]\, [0, a] \mid a \in L_{27}(13)$ and $x \in 4Z^+ \cup \{0\}\} \subseteq SG$ be the interval loop interval subsemiring of SG. Clearly H is not an interval loop interval ideal of SG.

***Example 4.2.54:*** Let $G = \{[0, a] \mid a \in L_5(4) = \{e, \overline{1}, \overline{2}, \overline{3}, \overline{4}, \overline{5}\} = \{e, x_0, x_1, x_2, \ldots, x_5\}$ that is $x_n = \overline{n}$ ; $0 \leq n \leq \overline{5}\}\}$ be an interval loop. $S = \{[0, a] \mid a \in Q^+ \cup \{0\}\}$ be the interval semiring.

$$SG = \left\{\sum_{i=0}^{5}[0, a_i][0, x_i] \,\middle|\, x_i \in L_5(4) \text{ and } a_i \in Q^+ \cup \{0\}\right\}$$

be the interval loop interval semiring.

Consider $P = \langle x_0 + x_1 + x_2 + x_3 + x_4 + x_5 \rangle \subseteq SG$, P is an interval loop interval subsemiring as well as ideal.

Now we proceed onto illustrate interval loop interval semiring of finite order; the definition is left as an exercise for the reader.

*Example 4.2.55:* Let G = {[0, a] | a ∈ $L_7(4)$} be an interval loop. S = {[0, x] | x ∈ $Z_{12}$} be the interval semiring. SG = {Σ [0, x] [0, a] | a ∈ G and x ∈ $Z_{12}$} be the interval loop interval semiring. Clearly SG is of finite order. SG has zero divisors, subsemirings and ideals.

Consider P = {Σ [0, x] [0, a] | x ∈ {0, 2, 4, 6, 8, 10} ⊆ $Z_{12}$, a ∈ $L_7(4)$} ⊆ SG; P is an interval loop interval semiring which is also an ideal.

Take x = [0, 2] [0, 3] + [0, 6] [0, 5] + [0, 4] [0, 6] and y = [0, 6] [0, 2] in SG.

x . y =  ([0, 2] [0, 3] + [0, 6] [0, 5] + [0, 4] [0, 6]) ([0, 6] [0, 2])
     =  [0, 2] [0, 6] [0, 3] [0, 2] + [0, 6] [0, 6] [0, 6] [0, 2] + [0, 4] [0, 6] [0, 6] [0, 2]
     =  0. (0, 3.2]) + 0 ([0, 5.2]) + 0 ([0, 6.2])
     =  0.

Thus SG has non trivial zero divisors. Further SG is not a strict interval loop interval semiring.

For [0, 3] [0, 5] + [0, 9] [0, 6]= 0 but none of the two terms is zero. [0, 1] [0, e] + [0, 1] [0, 1] + [0, 1] [0, 2] + … + [0, 1] [0, 7] = x ∈ SG is such that $x^2 = 0$.

*Example 4.2.56:* Let S = {[0, a] | a ∈ $Z_{30}$} be an interval semiring. G = {[0, x] | x ∈ $L_5(2)$} be an interval loop. SG the interval loop interval semiring is of finite order. SG is non commutative, SG has zero divisors and SG has ideals and subsemirings.

Now having seen examples of interval loop interval semiring now we proceed onto give examples of interval groupoid interval semirings. However the task of defining these terms are left as an exercise as it is direct.

**Example 4.2.57:** Let $G = \{[0, a], *, (3, 4), a \in Z_{12}\}$ be an interval groupoid. $S = \{[0, x] \mid x \in Z^+ \cup \{0\}\}$ be an interval semiring. SG is the interval groupoid interval semiring of infinite order. This has subsemirings given by $P = \{\Sigma [0, x] [0, g]$ where $x \in 3Z^+ \cup \{0\}$ and $[0, g] \in G\} \subseteq SG$.

**Example 4.2.58:** Let $G = \{[0, a], *, (3, 11), a \in Z_{12}\}$ be an interval groupoid. Let $S = \{[0, a] \mid a \in Q^+ \cup \{0\}\}$ be an interval semiring. SG the interval groupoid interval semiring is of infinite order. $P = \{\Sigma [0, x] [0, g] \mid x \in Z^+ \cup \{0\}$ and $g \in Z_{12}\} \subseteq$ SG is an interval groupoid interval subsemiring. However P is not an ideal of SG.

**Example 4.2.59:** Let $G = \{[0, a], *, (8, 12), a \in Z_{35}\}$ be an interval groupoid. $S = \{[0, x] \mid x \in Z^+ \cup \{0\}\}$ be an interval semiring. SG the interval groupoid interval semiring. $P = \{\Sigma [0, x] [0, a] \mid x \in 3Z^+ \cup \{0\}$ and $a \in Z_{35}\} \subseteq$ SG is an ideal of SG. Infact SG has infinite number of ideals.

Now having seen examples of interval groupoid interval semirings of infinite order we now proceed onto give examples of interval groupoid interval semirings of finite order.

**Example 4.2.60:** Let $G = \{[0, a], *, (13, 7), Z_{20}\}$ be an interval groupoid. $S = \{[0, x] \mid x \in Z_{10}\}$ be an interval semiring. SG is the interval groupoid interval semiring of finite order. Consider $P = \{\Sigma[0, x] [0, g] \mid x \in \{0, 2, 4, 5, 6, 8\} \subseteq Z_{10}$ and $g \in Z_{20}\} \subseteq$ SG, P is not only a subsemiring but also an ideal of SG.

**Example 4.2.61:** Let $G = \{[0, g], *, (23, 17), a \in Z_{29}\}$ be an interval groupoid. $S = \{[0, x] \mid x \in Z_{24}\}$ be an interval semiring. SG is the interval groupoid interval semiring of finite order.
  Consider $P = \{\Sigma [0, x] [0, g] \mid x \in \{0, 2, 4, 5, 6, 8, 10, \ldots, 22\} \subseteq Z_{24}$ and $g \in Z_{29}\} \subseteq$ SG is an interval groupoid interval subsemiring which is also an ideal of SG.

Now we proceed onto give examples of interval semigroup interval semirings.

***Example 4.2.62:*** Let $G = \{[0, g] \mid g \in Z_{32}\}$ be an interval semiring under multiplication modulo 32.

Let $S = \{[0, a] \mid a \in Z^+ \cup \{0\}\}$ be an interval semiring. $SG = \{\Sigma\, [0, a]\, [0, g] \mid a \in Z^+ \cup \{0\}$ and $g \in Z_{32}\}$ be the interval semigroup interval semiring.

SG is of infinite order. SG has zero divisors.

$P = \{\Sigma\, [0, a]\, [0, g] \mid a \in 3Z^+ \cup \{0\}$ and $g \in \{0, 4, 8, 12, 16, 20, 24, 28\} \subseteq Z_{32}\} \subseteq SG$ is a interval semigroup interval subsemiring of SG. It is easily verified P is an ideal of SG.

***Example 4.2.63:*** Let $G = \{[0, g] \mid g \in Z_{45}\}$ be an interval semigroup. $S = \{[0, a] \mid a \in Q^+ \cup \{0\}\}$ be an interval semiring. SG be the interval semigroup interval semiring. $P = \{\Sigma\, [0, x]\, [0, g] \mid g \in Z_{45}, a \in Z^+ \cup \{0\}\} \subseteq SG$; is an interval semigroup interval subsemiring of SG, but P is not an ideal of SG. However SG has zero divisors.

***Example 4.2.64:*** Let $G = \{[0, g] \mid g \in Q^+ \cup \{0\}\}$ be an interval semigroup. $S = \{[0, x] \mid x \in Z^+ \cup \{0\}\}$ be an interval semiring. SG be the interval semigroup interval semiring. SG is of infinite order. Take $P = \{\Sigma\, [0, x]\, [0, g] \mid x \in 3Z^+ \cup \{0\}$ and $g \in Q^+ \cup \{0\}\} \subseteq SG$, P is an interval semigroup interval subsemiring infact an ideal of SG.

If we take $T = \{\Sigma\, [0, x]\, [0, h] \mid x \in 5Z^+ \cup \{0\}$ and $h \in 3Z^+ \cup \{0\}\} \subseteq SG$, P is only a interval semigroup interval subsemiring and not an ideal of SG.

***Example 4.2.65:*** Let $G = \{[0, g] \mid g \in Z_{40}\}$ be an interval semigroup and $S = \{[0, a] \mid a \in Z_{20}\}$ be an interval semiring. $SG = \{\Sigma\, [0, x]\, [0, g] \mid g \in Z_{40}, x \in Z_{20}\}$ is the interval semigroup interval semiring of finite order. Consider $P = \{\Sigma\, [0, y]\, [0, h] \mid y \in \{0, 10\} \subseteq Z_{20}, h \in Z_{40}\} \subseteq SG$; P is a interval semigroup interval subsemiring of SG of finite order.

***Example 4.2.66:*** Let $S = \{[0, a] \mid a \in Z_{120}\}$ be the interval semiring. $G = \{[0, x] \mid x \in Z_{12}\}$ be the interval semigroup. SG be the interval semigroup interval semiring of finite order. $P =$

$\{\Sigma [0, a] [0, g] \mid a \in \{0, 10, 20, 30, \ldots, 110\} \subseteq Z_{120}, g \in \{0, 6\}\}$ $\subseteq$ SG, P is an interval semigroup interval subsemigroup and also an ideal of SG. This interval semigroup interval semiring has non trivial zero divisors and it is not a strict semiring.

Now we proceed onto give examples of interval group interval semirings.

***Example 4.2.67:*** Let G = $\{[0, g] \mid g \in Z_7 \setminus \{0\}\}$ be an interval group under multiplication modulo 7. S = $\{[0, x] \mid x \in Q^+ \cup \{0\}\}$ be an interval semiring. SG = $\{\Sigma [0, x] [0, g] \mid x \in Q^+ \cup \{0\}, g \in Z_7 \setminus \{0\}\}$ is the interval group interval semiring. Clearly SG is of infinite order.

***Example 4.2.68:*** Let G = $\{[0, x] \mid x \in Z_{23} \setminus \{0\}\}$ be an interval group. S = $\{[0, a] \mid a \in Z^+ \cup \{0\}\}$ be the interval semiring. SG is the interval group interval semiring. SG is of infinite order. SG is a strict semiring. P = $\{\Sigma [0, a] [0, x] \mid a \in 3Z^+ \cup \{0\}, x \in Z_{23} \setminus \{0\}\} \subseteq$ SG is the interval group interval subsemiring which is an ideal of SG.

***Example 4.2.69:*** Let G = $\{[0, g] \mid g \in Z_{23} \setminus \{0\}\}$ be an interval group. S = $\{[0, a] \mid a \in R^+ \cup \{0\}\}$ be the interval semiring. SG = $\{\Sigma [0, a] [0, g] \mid a \in R^+ \cup \{0\}$ and $g \in Z_{23} \setminus \{0\}\}$ be the interval group interval semiring. SG is commutative. SG is a strict interval group interval semiring.

SG has no zero divisors. SG has idempotents.

Consider P = $\{\Sigma [0, a] [0, g] \mid a \in Q^+ \cup \{0\}$ and $g \in Z_{23} \setminus \{0\}\} \subseteq$ SG, P is an interval subsemiring and is not an ideal of SG. Infact SG has infinite number of subsemirings.

Consider T = $\{\langle \Sigma [0,1] [0, g] \rangle \mid g \in Z_{23} \setminus \{0\}\}$
= $\langle [0, 1] [0, \overline{1}] + [0, 1] [0, \overline{2}] + \ldots + [0, 1] [0, \overline{22}] \rangle \subseteq$ SG.
$\{\overline{1}, \overline{2}, \ldots, \overline{22}\} \in Z_{23} \setminus \{0\}$.

T is a interval subsemigroup as well as an ideal of SG.

Inview of this we have the following observations.

**THEOREM 4.2.17:** *Let $G = \{[0, g] \mid g \in Z_p \setminus \{0\}; p$ a prime$\}$ be an interval group. $S = \{[0, x] \mid x \in Z^+ \cup \{0\}$ or $R^+ \cup \{0\}$ or $Q^+ \cup \{0\}\}$ be the interval semiring. SG be the interval group interval semiring.*

*The following conditions are satisfied by SG.*
1. *SG is commutative, strict infinite interval semiring.*
2. *SG is a S-semiring.*
3. *SG has no zero divisors.*
4. *SG has idempotents if and only if S is built over $Q^+ \cup \{0\}$ or $R^+ \cup \{0\}$.*
5. *SG has infinite number of ideals if S is built over $Z^+ \cup \{0\}$.*
6. *SG has only one ideal if S is built over $R^+ \cup \{0\}$ or $Q^+ \cup \{0\}$.*

The proof is direct and hence is left as an exercise for the reader.

*Example 4.2.70:* Let $G = \{[0, a] \mid a \in Q^+\}$ be an interval group under multiplication. $S = \{[0, x] \mid x \in Z^+ \cup \{0\}\}$ the interval semiring.

SG be the interval group interval semiring. SG is of infinite order and commutative. SG has ideals given by $P = \{\Sigma[0, x] [0, a] \mid x \in 17Z^+ \cup \{0\}$ and $a \in Q^+\} \subseteq SG$; P is easily verified to be an ideal of SG. Infact SG has infinite number of ideals and subsemirings in it.

*Example 4.2.71:* Let $G = \{[0, a] \mid a \in R^+\}$ be an interval group. $S = \{[0, x] \mid x \in R^+ \cup \{0\}\}$ be an interval semiring. SG be the interval group interval semiring.

SG is of infinite order and is commutative. SG is a strict semiring. SG has no zero divisors or idempotents. SG has only interval subsemirings and no ideals. Take $P = \{[0, a] \mid a \in Q^+\} \subseteq G$ and $T = \{[0, x] \mid x \in Z^+ \cup \{0\}\} \subseteq S$. $TP = \{\Sigma [0, x] [0, a]$ where $a \in Q^+$ and $x \in Z^+ \cup \{0\}\} \subseteq SG$. TP is only an interval subsemiring and not an ideal.

The reader is left with the task of proving SG has no ideals.

We will now see some examples of finite interval group interval semirings.

***Example 4.2.72:*** Let $G = \{[0, a] \mid a \in Z_{13} \setminus \{0\}\}$ be an interval group. $T = \{[0, x] \mid x \in Z_{20}\}$ be an interval semiring. $TG = \{\Sigma [0, x] [0, a] \mid a \in Z_{13} \setminus \{0\}$ and $x \in Z_{20}\}$ be the interval group interval semiring. TG is of finite order. TG has zero divisors. TG is not a strict semiring but TG is commutative.

**THEOREM 4.2.18:** *Let $G = \{[0, a] \mid a \in Z_p \setminus \{0\}$ and $p$ a prime$\}$ be a finite interval group. $S = \{[0, x] \mid x \in Z_n\}$ be an interval semiring. $SG = \{\Sigma [0, x] [0, a] \mid x \in Z_n$ and $a \in Z_p \setminus \{0\}\}$ be the interval group interval semiring. The following are true.*
  1. *SG is a commutative finite semiring.*
  2. *SG is not strict.*
  3. *SG has ideals provided n is not a prime.*

The proof is direct and hence is left as an exercise to the reader. We can construct interval semirings which are non commutative.

***Example 4.2.73:*** Let $G = \{S_X$ the interval symmetric group where $X = \{[0, a_1], [0, a_2], [0, a_3]\}$ and $S = \{[0, x] \mid x \in Q^+ \cup \{0\}\}$ be an interval semiring.

SG is the interval group interval semiring of infinite order and SG is a non commutative semiring.

Thus using the interval symmetric group $S_X$, where $X = \{[0, a_1], \ldots, [0, a_n]\}$ we can get the class of non commutative interval semirings of infinite order. It is pertinent to mention here that $S(\langle X \rangle)$ will give the special interval semigroups and $S(X)$ the interval semigroups using all these we can get both infinite and finite classes of interval semigroup interval semirings.

**Chapter Five**

# INTERVAL NEUTROSOPHIC SEMIRINGS

In this chapter we define different types of neutrosophic interval semirings. We use neutrosophic intervals in the place of real intervals. The two types of neutrosophic intervals which we will be using are pure neutrosophic intervals and neutrosophic intervals.

Pure neutrosophic intervals will be of the form $\{[0, aI] \mid aI \in Z_nI \text{ or } Z^+I \cup \{0\} \text{ or } Q^+I \cup \{0\} \text{ or } R^+I \cup \{0\}\}$. Neutrosophic intervals would be of the form $\{[0, a + bI] \mid a, b \in Z_n \text{ or } Z^+ \cup \{0\} \text{ or } Q^+ \cup \{0\} \text{ or } R^+ \cup \{0\}\}$. Using these intervals we will just briefly define the algebraic structures and proceed onto define neutrosophic interval semirings.

**DEFINITION 5.1:** *Let $S = \{[0, aI] \mid a \in Z_n$ or $Z^+ \cup \{0\}$ or $Q^+ \cup \{0\}$ or $R^+ \cup \{0\}$; S under interval addition and multiplication is an interval semiring known as the pure neutrosophic interval semiring.*

We will give examples of them.

***Example 5.1:*** Let $S = \{[0, aI] \mid a \in Z_9, I^2 = I\}$. S is a finite pure neutrosophic interval semiring.

***Example 5.2:*** Let $S = \{[0, aI] \mid a \in Z^+I \cup \{0\}\}$ be an infinite pure neutrosophic interval semiring of infinite order. Both these interval semirings are commutative.
   If in the definition 5.1 we replace the interval [0, aI] by $\{[0, a+bI]$ ; $a, b \in Z_n$ or $Z^+ \cup \{0\}$ or $Q^+ \cup \{0\}$ or $R^+ \cup \{0\}\}$ is a neutrosophic interval semiring and is not a pure neutrosophic interval semiring.

We will illustrate this situation also by some examples.

***Example 5.3:*** Let $S = \{[0, a+bI] \mid a, b \in R^+ \cup \{0\}\}$ be the neutrosophic interval semiring. Now $P = \{[0, bI] \mid b \in R^+ \cup \{0\}\} \subseteq S$ as a neutrosophic interval subsemiring which is clearly pure neutrosophic.

***Example 5.4:*** Let $S = \{[0, a+bI] \mid a, b \in Z_{20}\}$ be a neutrosophic interval semiring. Clearly S is of finite order; where as the neutrosophic interval semiring given in example 5.3 is of infinite order. Both of them are commutative and both contain a subsemiring which is a pure neutrosophic interval subsemiring.

We can define substructures in them; as it is a matter of routine and so is left as an exercise to the reader. However we will illustrate this by some examples.

***Example 5.5:*** Let $S = \{[0, aI] \mid aI \in Z^+I \cup \{0\}\}$ be the interval pure neutrosophic semiring (also semifield). $P = \{[0, aI] \mid a \in 5Z^+ \cup \{0\}\} \subseteq S$ is an ideal as well as a subsemiring of S.

***Example 5.6:*** Let $S = \{[0, aI] \mid a \in Z_7\}$ be the pure neutrosophic interval semiring, S has no subsemirings. Hence has no ideals.

In view of this we have the following theorem.

**THEOREM 5.1:** *Let $S = \{[0, aI] \mid aI \in Z_p$; p a prime$\}$ be a pure neutrosophic interval semiring. S has no subsemiring and hence has no proper ideals.*

The proof is straight forward hence left as an exercise for the reader.

**THEOREM 5.2:** *Let $S = \{[0, aI] \mid aI \in Q^+I \cup \{0\}\}$ be a pure neutrosophic interval semiring. S has subsemirings but has no ideals.*

Proof is left as an exercise for the reader.

***Example 5.7:*** Let $S = \{[0, a + bI] \mid a, b \in Q^+I \cup \{0\}\}$ be a neutrosophic interval semiring. $P = \{[0, bI] \mid b \in Q^+ \cup \{0\}\} \subseteq S$ is a pure neutrosophic interval subsemiring of S. Also $T = \{[0, a] \mid a \in Q^+ \cup \{0\}\} \subseteq S$ is an interval subsemiring which is not a neutrosophic interval subsemiring of S.
　　This interval subsemiring will be defined as the subsemiring.

We will give more examples.

***Example 5.8:*** Let $S = \{[0, a + bI] \mid a, b \in Z^+ \cup \{0\}\}$ be a neutrosophic interval semiring. S has both pure neutrosophic interval subsemiring and a pseudo neutrosophic interval subsemiring.

The following theorems are of importance.

**THEOREM 5.3:** *Let $S = \{[0, aI] \mid aI \in Z_n$ or $Z^+ \cup \{0\}$ or $Q^+ \cup \{0\}$ or $R^+ \cup \{0\}\}$ be a pure neutrosophic interval semiring. S has no pseudo neutrosophic interval subsemiring.*

**THEOREM 5.4:** *Let $S = \{[0, a + bI] \mid a,b \in Z_n$ or $Z^+ \cup \{0\}$ or $Q^+ \cup \{0\}$ or $R^+ \cup \{0\}\}$ be an interval neutrosophic semiring. S has both pure neutrosophic interval subsemiring and pseudo neutrosophic interval subsemiring.*

**THEOREM 5.5**: *Let $S = \{[0, a+bI] \mid a,b \in Z^+ \cup \{0\}\}$ be a neutrosophic interval semiring. S has ideals.*

**THEOREM 5.6:** *Let $S = \{[0, \alpha I] \mid \alpha \in Z^+ \cup \{0\}\}$ be a pure neutrosophic interval semiring. S has ideals.*

Now we proceed onto define and give examples of other types of neutrosophic semirings.

**DEFINITION 5.2:** *Let*

$$S = \left\{ \sum_{i=0}^{\infty} [0, aI\,]x^i \,\middle|\, a \in Z_n \right.$$

*or $Z^+ \cup \{0\}$ or $Q^+ \cup \{0\}$ or $R^+ \cup \{0\}\}$; S under polynomial addition and multiplication is an interval semiring, known as the pure neutrosophic interval polynomial semiring.*

We will illustrate this by one or two examples.

*Example 5.9:* Let

$$S = \left\{ \sum_{i=0}^{\infty} [0, aI]x^i \,\middle|\, a \in Z_{12} \right\}$$

be an infinite pure neutrosophic polynomial interval semiring.

*Example 5.10:* Let

$$S = \left\{ \sum_{i=0}^{\infty} [0, aI]x^i \,\middle|\, a \in Z^+ \cup \{0\} \right\}$$

be an infinite pure neutrosophic polynomial interval semiring.

*Example 5.11:* Let

$$T = \left\{ \sum_{i=0}^{9} [0, aI]x^9 \,\middle|\, x^{10} = 1, a \in Q^+ \cup \{0\} \right\}$$

be an infinite pure neutrosophic polynomial interval semiring.

***Example 5.12:*** Let
$$S = \left\{ \sum_{i=0}^{5}[0,aI]x^5 \,\middle|\, x^6 = 1, \ a \in Z_{12} \right\}$$
be a pure neutrosophic polynomial interval semiring of finite order.

***Example 5.13:*** Let
$$S = \left\{ \sum_{i=0}^{\infty}[0,a+bI]x^i \,\middle|\, a,b \in Q^+ \cup \{0\} \right\}$$
be a neutrosophic interval polynomial semiring of infinite order.

***Example 5.14:*** Let
$$S = \left\{ \sum_{i=0}^{\infty}[0,a+bI]x^i \,\middle|\, a,b \in R^+ \cup \{0\} \right\}$$
be the neutrosophic interval polynomial semiring of infinite order.

***Example 5.15:*** Let
$$S = \left\{ \sum_{i=0}^{\infty}[0,a+bI]x^i \,\middle|\, a,b \in Z_{16} \right\}$$
be a neutrosophic interval polynomial semiring of infinite order.

***Example 5.16:*** Let
$$S = \left\{ \sum_{i=0}^{\infty}[0,a+bI]x^i \,\middle|\, a,b \in Z^+ \cup \{0\} \right\}$$
be a neutrosophic interval polynomial semiring of infinite order.

***Example 5.17:*** Let
$$S = \left\{ \sum_{i=0}^{7}[0,a+bI]x^i \,\middle|\, a,b \in Z_{19}; x^8 = 1 \right\}$$
be the neutrosophic interval polynomial semiring of finite order.

We will proceed onto give examples of neutrosophic polynomial interval subsemirings and ideals.

*Example 5.18:* Let
$$S = \left\{ \sum_{i=0}^{\infty} [0, aI] x^i \,\Big|\, a \in Z^+ \cup \{0\} \right\}$$
be a pure neutrosophic interval polynomial semiring.
Consider
$$P = \left\{ \sum_{i=0}^{\infty} [0, aI] x^i \,\Big|\, a \in 3Z^+ \cup \{0\} \right\} \subseteq S;$$
P is a subsemiring as well as ideal of S.
Take
$$T = \left\{ \sum_{i=0}^{\infty} [0, aI] x^{2i} \,\Big|\, a \in 7Z^+ \cup \{0\} \right\} \subseteq S;$$
T is a subsemiring and is not an ideal of S.

Thus S has infinite number of subsemirings which are not ideals and subsemirings which are ideals of S.

*Example 5.19:* Let
$$S = \left\{ \sum_{i=0}^{\infty} [0, aI] x^i \,\Big|\, aI \in R^+ I \cup \{0\} \right\}$$
be a pure neutrosophic interval polynomial semiring.
Consider
$$T = \left\{ \sum_{i=0}^{\infty} [0, aI] x^i \,\Big|\, a \in Q^+ \cup \{0\} \right\} \subseteq S$$
be a pure neutrosophic interval polynomial subsemiring.

Clearly T is not an ideal of S. S has infinite number of subsemirings.

*Example 5.20:* Let
$$S = \left\{ \sum_{i=0}^{\infty} [0, aI] x^i \,\Big|\, a \in Z_{12} \right\}$$
be a pure neutrosophic interval polynomial semiring.
Consider

$$P = \left\{ \sum_{i=0}^{\infty} [0, aI]x^i \,\middle|\, a \in \{0, 2, 4, 6, 8, 10\} \subseteq Z_{12} \right\} \subseteq S,$$

P is both a subsemiring as well as ideal of S.
Take

$$T = \left\{ \sum_{i=0}^{\infty} [0, aI]x^{2i} \,\middle|\, a \in \{0, 3, 6, 9\} \subseteq Z_{12} \right\} \subseteq S;$$

T is only a subsemiring and not an ideal of S.

Thus S has many subsemirings which are not ideals.

*Example 5.21:* Let

$$E = \left\{ \sum_{i=0}^{9} [0, aI]x^i \,\middle|\, x^{10} = 1, a \in Z_{70} \right\}$$

be a pure neutrosophic interval polynomial semiring.
Consider

$$F = \left\{ \sum_{i=0}^{9} [0, aI]x^i \,\middle|\, a \in \{0, 10, 20, ..., 60\} \subseteq Z_{70} \right\} \subseteq E;$$

F is a subsemiring as well as an ideal of E.

Having defined substructures in pure neutrosophic interval polynomial semirings we now proceed onto give examples of neutrosophic interval polynomial semirings and their substructures.

*Example 5.22:* Let

$$S = \left\{ \sum_{i=0}^{\infty} [0, aI + b]x^i \,\middle|\, a, b \in Q^+ \cup \{0\} \right\}$$

be a neutrosophic polynomial interval semiring.

$$T = \left\{ \sum_{i=0}^{\infty} [0, b]x^i \,\middle|\, b \in Q^+ \cup \{0\} \right\} \subseteq S$$

is a pseudo neutrosophic polynomial interval subsemiring of S.

$$W = \left\{ \sum_{i=0}^{\infty} [0, aI]x^i \,\middle|\, a \in Q^+ \cup \{0\} \right\} \subseteq S$$

is a pure neutrosophic interval polynomial subsemiring of S. Clearly W is also an ideal of S, however T is not an ideal of S. Infact S has several subrings and ideals.

*Example 5.23:* Let
$$S = \left\{ \sum_{i=0}^{9} [0, a + bI]x^i \,\bigg|\, a, b \in Z_{42}; x^{10} = 1 \right\}$$
be a neutrosophic interval polynomial semiring.
$$P = \left\{ \sum_{i=0}^{9} [0, bI]x^i \,\bigg|\, b \in Z_{42}; x^{10} = 1 \right\} \subseteq S$$
is a pure neutrosophic interval polynomial subsemiring as well as ideal of S.
$$T = \left\{ \sum_{i=0}^{9} [0, a + bI]x^i \,\bigg|\, a, b \in Z_{42}; x^{10} = 1 \right\} \subseteq S$$
is only pseudo neutrosophic interval polynomial subsemiring of S and is not an ideal of S.
$$W = \left\{ \sum_{i=0}^{9} [0, a + bI]x^i \,\bigg|\, x^{10} = 1, a, b \in \{0, 7, 14, 21, 28, 35\} \subseteq Z_{42} \right\}$$
$\subseteq S$ is both a subsemiring as well as an ideal of S.

Now we can as in case of interval matrix semirings define the notion of neutrosophic interval matrix semirings. Here we give only examples of them and their substructures.

*Example 5.24:* Let $S = \{([0, a_1I], \ldots, [0, a_9I])$ where $a_i \in Z_{40}$, $1 \leq i \leq 9\}$ be a pure neutrosophic row interval matrix semiring. $P = \{([0, aI], [0, bI], 0, 0, 0, 0, 0, [0, cI], [0, dI]) \mid a, b, c, d \in Z_{40}\} \subseteq S$ is a pure neutrosophic interval row matrix subsemiring as well as an ideal of S. Clearly S is finite order.

*Example 5.25:* Let $S = \{([0, aI], [0, bI], [0, cI]) \mid a, b, c \in Z^+ \cup \{0\}\}$ be a pure neutrosophic row interval matrix semiring. S has both ideals and subsemirings.

*Example 5.26:* Let $S = \{([0, a_1 + b_1I], [0, a_2 + b_2I], [0, a_3 + b_3I], [0, a_4 + b_4I] \mid a_i, b_i \in Z_7, 1 \leq i \leq 4\}$ be a neutrosophic row

interval matrix semiring. S has T = {([0, $a_1$], [0, $a_2$], [0, $a_3$], [0, $a_4$]} | $a_i \in Z_7$, $1 \le i \le 4$} $\subseteq$ S to be a pseudo neutrosophic row interval matrix subsemiring. W = {([0, $b_1$I], [0, $b_2$I], [0, $b_3$I], [0, $b_4$I]) | $b_i \in Z_7$, $1 \le i \le 4$} $\subseteq$ S is a pure neutrosophic row matrix interval subsemiring.

In view of this we have the following theorems the proof of which are direct.

**THEOREM 5.7:** *Let S = {([0, $a_1 + b_1$I], [0, $a_2 + b_2$I], ..., [0, $a_n + b_n$I]) | $a_i, b_i \in Z_p$, p a prime, $1 \le i \le n$} be a neutrosophic row matrix interval semiring. S has pseudo neutrosophic row interval matrix subsemiring, pure neutrosophic row interval matrix semiring and neutrosophic row matrix interval subsemiring.*

**THEOREM 5.8:** *Let S = {([0, $a_1$I], [0, $a_2$I], ..., [0, $a_n$I]) | $a_i \in Z_p$, p a prime, $1 \le i \le n$} be a pure neutrosophic interval matrix semiring. S has both pure neutrosophic interval matrix ideals and subsemirings.*

**THEOREM 5.9:** *Let S = {([0, aI], [0, aI], ..., [0, aI]) | $a \in Z_p$, p a prime} be a pure neutrosophic interval matrix semiring. S has no subsemirings.*

**COROLLARY 5.1:** *S given in theorem 5.9 has no ideals.*

*Example 5.27:* Let S = {([0, aI], [0, aI], [0, aI], [0, aI], [0, aI], [0, aI]) | $a \in Z_{11}$} be a pure neutrosophic interval matrix semiring. S has no subsemirings and ideals.

Now having seen pure neutrosophic and neutrosophic interval row matrix semirings we now proeed onto give examples of pure neutrosophic square matrix interval semrings.

*Example 5.28:* Let S = {All $3 \times 3$ pure neutrosophic interval matrices of the form [0, aI] where $a \in Z_{13}$} be a pure neutrosophic $3 \times 3$ interval matrix semiring.

Consider

$$A = \left\{ \begin{bmatrix} [0,aI] & [0,aI] & [0,aI] \\ [0,aI] & [0,aI] & [0,aI] \\ [0,aI] & [0,aI] & [0,aI] \end{bmatrix} \middle| a \in Z_{13} \right\} \subseteq S$$

is a pure neutrosophic $3 \times 3$ interval matrix subsemiring as well as ideal of S.

Take W = {All $3 \times 3$ pure neutrosophic upper triangular interval matrices with intervals of the form [0, aI] | a $\in Z_{13}$} $\subseteq$ S, is only a subsemiring and not an ideal of S.

*Example 5.29:* Let S = {$5 \times 5$ pure neutrosophic interval matrices with intervals of the form [0, aI] where a $\in Z^+ \cup \{0\}$} be a pure neutrosophic square matrix interval semiring. S has subsemirings and ideals given by P = {$5 \times 5$ pure neutrosophic intervals matrices with intervals of the form [0, aI] where a $\in 7Z^+ \cup \{0\}$} $\subseteq$ S; P is an ideal of S.

*Example 5.30:* Let

$$S = \left\{ \begin{bmatrix} [0,a] & [0,b] \\ [0,c] & [0,d] \end{bmatrix} \middle| \right.$$

$a = a_1 + a_2I$, $b = b_1 + b_2I$, $c = c_1 + c_2I$ and $d = d_1 + d_2I$, $a_i$, $c_i$, $b_i$, $d_i$ $\in Z_{12}$, $1 \leq i \leq 2$} be a neutrosophic interval matrix semiring. Clearly

$$P = \left\{ \begin{bmatrix} [0,a_1] & [0,b_1] \\ [0,c_1] & [0,d_1] \end{bmatrix} \middle| a_1, b_1, c_1, d_1 \in Z_{12} \right\} \subseteq S$$

is a pseudo neutrosophic interval matrix subsemiring of S which is not an ideal of S.

$$T = \left\{ \begin{bmatrix} [0,a_2I] & [0,b_2I] \\ [0,c_2I] & [0,d_2I] \end{bmatrix} \middle| a_2, b_2, c_2, d_2 \in Z_{12} \right\} \subseteq S$$

is a pure neutrosophic interval matrix subsemiring of S which is also an ideal of S.

Consider

$$W = \left\{ \begin{bmatrix} [0,a_1 + a_2I] & [0,b_1 + b_2I] \\ [0,c_1 + c_2I] & [0,d_1 + d_2I] \end{bmatrix} \middle| \right.$$

$d_i, a_i, c_i, b_i \in \{0, 2, 4, 6, 8, 10\} \subseteq Z_{12}, 1 \le i \le 2\} \subseteq S$; W is an ideal of S.

*Example 5.31:* Let S = {All 10 × 10 neutrosophic interval matrices with entries of the form [0, a + bI] where a, b ∈ $Z^+ \cup \{0\}$} be the neutrosophic matrix interval semiring. Consider W = {All 10 ×10 upper triangular interval neutrosophic matrices with intervals of the form [0, a + bI] where a, b ∈ $3Z^+ \cup \{0\}$} ⊆ S. W is only a neutrosophic interval matrix subsemiring of S and is not an ideal of S.

V = {All 10 × 10 neutrosophic interval matrices with intervals of the form [0, a+bI] where a, b ∈ $7Z^+ \cup \{0\}$} ⊆ S is an ideal of S.

We can derive all the properties associated with usual interval matrix semirings in case of pure neutrosophic and neutrosophic interval matrix semrings. We can also use either neutrosophic groupoids or neutrosophic groups or neutrosophic semigroups or neutrosophic loops and use usual neutrosophic interval semirings to built semirings both associative and non associative in an analogous way or use neutrosophic interval groupoids or groups or loops or semigroups over usual semirings to get both neutrosophic and pure neutrosophic interval semirings which are a matter of routine. We will give one or two examples of them.

*Example 5.32:* Let G = {$Z_{12}$, *, (7, 5)} be a groupoid. S = {[0, aI] | a ∈ $Z^+ \cup \{0\}$} be a pure neutrosophic semiring.
Consider
$$SG = \left\{ \sum_g [0, aI]g \,\middle|\, a \in Z^+ \cup \{0\} \right\}$$
and g ∈ G is the finite formal sums}; SG is an interval neutrosophic semiring which is non associative known as the groupoid interval neutrosophic semiring.

*Example 5.33:* Let G = {$Z_{19}$, *, (3, 6)} be a groupoid. S = {[0, a + bI] | a, b ∈ $Z_{45}$} be an interval neutrosophic semiring. SG is the groupoid interval neutrosophic semiring of the groupoid G

over the interval neutrosophic semiring S. Clearly SG is non associative and is of finite order.

***Example 5.34:*** Let $L = L_9(7)$ be a loop. $S = \{[0, aI] \mid a \in Q^+ \cup \{0\}\}$ be a pure neutrosophic interval semiring. SL be the loop interval neutrosophic semiring of the loop L over the pure neutrosophic interval semiring S. SL is non associative neutrosophic interval semiring of infinite order.

***Example 5.35:*** Let $L = L_{13}(7)$ be a loop of order 14. $S = \{[0, a + bI] \mid a, b \in Z_3\}$ be a neutrosophic interval semiring. SL be the loop neutrosophic interval semiring. SL is non associative but commutative, SL has non associative pseudo neutrosophic interval subsemiring of finite order given by

$$P = \left\{ \sum_g [0,a]g_i \,\middle|\, a \in Z_3 \text{ and } g_i \in L \right\} \subseteq SL$$

and pure neutrosophic interval subsemiring given by

$$T = \left\{ \sum_i [0,aI]g_i \,\middle|\, a \in Z_3 \text{ and } g_i \in L \right\}.$$

Further

$$V = \left\{ \sum [0_1, a+aI]g_i \,\middle|\, a \in Z_3 \text{ and } g_i \in L \right\} \subseteq S$$

is also a neutrosophic interval subsemiring of SL which is commutative.

***Example 5.36:*** Let $G = S_6$ be the symmetric group of degree 6. $S = \{[0, aI] \mid a \in Z^+ \cup \{0\}\}$ be the pure neutrosophic interval semiring.

$$SG = \left\{ \sum_g [0,aI]g \,;\right.$$

finite formal sums of the form with $a \in Z^+ \cup \{0\}$ and $g \in S_6\}$; SG is the group pure neutrosophic interval semiring of the group G over the semiring S. SG is associative but is a non commutative neutrosophic interval semiring of infinite order. SG has many subsemirings and ideals.

Take
$$T = \left\{ \sum_h [0, aI]h \mid a \in Z^+ \cup \{0\}; h \in A_6 \triangleleft S_6 \right\} \subseteq SG;$$
T is an ideal of SG.

***Example 5.37:*** Let $G = <g \mid g^5 = 1>$ be the cyclic group of degree 5. $S = \{[0, a + bI] \mid a, b \in Z_{10}\}$ be a neutrosophic interval semiring. SG the group neutrosophic interval semiring. SG is of finite order and commutative. SG has both pseudo neutrosophic interval subsemiring and pure neutrosophic interval subsemiring.

**THEOREM 5.10:** *Let $S_n$ be the symmetric group of degree n. $S = \{[0, a+bI] \mid a, b \in Q^+ \cup \{0\}\}$ be the neutrosophic interval semiring. $SS_n$ is the symmetric group neutrosophic interval semiring. $SS_n$ has ideals.*

*Proof:* We only give an hint to the proof of the theorem. $A_n \triangleleft S_n$; $A_n$ is the normal subgroup of $S_n$ $SA_n \subseteq SS_n$ and $SA_n$ is the neutrosophic interval subsemiring of $SS_n$ and is an ideal as $A_n$ is the normal subgroup of $S_n$.

It is interesting to note that in the above theorem $SS_n$ has pseudo neutrosophic interval subsemirings and pure neutrosophic interval subsemirings also.

Now having seen groupoid neutrosophic interval semirings, loop neutrosophic interval semirings and group neutrosophic interval semirings we now proceed onto give examples semigroup neutrosophic interval semirings.

***Example 5.38:*** Let $P = \{Z_{12}, \times\}$ be the semigroup under multiplication modulo 12. $S = \{[0, aI] \mid a \in Z_{20}\}$ be a pure neutrosophic interval semiring. SP the set of all finite formal sums of the form
$$\sum_{p_i} [0, aI]p_i ,$$
$a \in Z_{20}$ and $p_i \in P$ is the semigroup pure neutrosophic interval semiring. This semiring has subsemirings, ideals, zero divisors,

and idempotents. Further SP is a commutative semiring and is of finite order.

**Example 5.39:** Let P = {$Z_{25}$, under multiplication modulo 25} be a semigroup of order 25. S = {[0, a+bI] | a, b ∈ $Q^+$ ∪ {0}} be an interval neutrosophic semiring. SP the semigroup neutrosophic interval semiring is of infinite order, has subsemirings ideals and pseudo interval neutrosophic subsemirings and pure neutrosophic interval subsemirings.

**Example 5.40:** Let S(9) be the symmetric semigroup. S = {[0, aI] | a ∈ $Q^+$ ∪ {0}} be the pure neutrosophic interval semiring. SS (9) is the semigroup pure neutrosophic interval semiring of infinite order and is non commutative SS(9) has non trivial subsemirings and ideals in it.

**Example 5.41:** Let S (3) be the symmetric semigroup. S = {[0, aI] | a ∈ $Z_{14}$} be the pure interval neutrosophic semiring. SS(3) the semigroup interval neutrosophic semiring is of finite order, non commutative and has ideals and subsemirings SS(3) has also zero divisors.

**Example 5.42:** Let S (5) be the symmetric semigroup. S = {[0, a + bI] | a, b ∈ $R^+$ ∪ {0}} be the neutrosophic interval semiring. SS(5) is the semigroup neutrosophic interval semiring of infinite order non commutative but has subsemirings.

Now having seen these new structures we can still built neutrosophic interval semirings using neutrosophic interval semigroups, groupoids, groups and loops.

We will however illustrate this by examples.

**Example 5.43:** Let G = {[0, aI] | aI ∈ $Z_{11}I$ \ {0}} be a neutrosophic interval group under multiplication modulo 11. Take S = {$Z^+$ ∪ {0}} to be the semiring. SG = {Σ $a_i$ [0, gI] | gI ∈ $Z_{11}I$ \ {0}} is the group neutrosophic interval semiring of infinite order. SG is commutative and has subsemirings and ideals.

***Example 5.44:*** Let $G = \{[0, aI] \mid a \in Q^+\}$ be an interval neutrosophic group under multiplication. $S = \{R^+ \cup \{0\}\}$ be the semiring. SG the group neutrosophic interval semiring is of infinite order.

***Example 5.45:*** Let $P = \{[0, a+bI] \mid a, b \in Z^+ \cup \{0\}\}$ be an interval neutrosophic semigroup of infinite order under multiplication. $S = \{Q^+ \cup \{0\}\}$ be the semiring. SP is the semigroup interval neutrosophic semiring of the interval neutrosophic semigroup P over the semiring S.

***Example 5.46:*** Let $S = \{[0, aI] \mid a \in Z_{12}\}$ be an interval neutrosophic semigroup of finite order under multiplication. $F = \{Z^+ \cup \{0\}\}$ be the semiring. FS is the interval semigroup neutrosophic semiring. FS is commutative and is of infinite order has ideals and subsemirings.

***Example 5.47:*** Let $G = \{[0, aI] \mid a \in Z_{12}I, *, (3, 7)\}$ be an interval neutrosophic groupoid of finite order. $S = \{Z^+ \cup \{0\}\}$ be the semiring. SG is the groupoid interval neutrosophic semiring of the neutrosophic interval groupoid G over the semiring S. Clearly SG is a non commutative and non associative interval semiring of infinite order.

***Example 5.48:*** Let $S = \{[0, aI] \mid a \in Z^+ \cup \{0\}, *, (9, 12)\}$ be an interval neutrosophic groupoid of infinite order. $F = \{Q^+ \cup \{0\}\}$ be the semiring. FS is the neutrosophic interval groupoid semiring of infinite order which is non associative and contains finite formal sums of the form

$$\sum_i \alpha_i [0, sI]$$

where i runs over a finite number and $\alpha_i \in F$, $sI \in Z^+ \cup \{0\}$.

***Example 5.49:*** Let $G = \{[0, aI] \mid aI \in Z_7I, *, (2,3))$ be an interval neutrosophic groupoid. $F = \{Z_{12}\}$ be a finite semiring. FG is the interval neutrosophic groupoid semiring of finite order which is non associative and non commutative.

***Example 5.50:*** Let $L = \{L_nI(m)\} = \{$The interval loop which is neutrosophic having interval elements $[0, eI], [0, I], [0, 2I], \ldots, [0, nI]$, n-odd $n > 3$, with operation * such that $[0, aI] * [0, bI] = mbI + (m-1)aI \pmod{n}$ where $m \leq n$, $(m, n) = 1$ $aI, bI \in L$ and $(m-1, n) = 1\}$ be the interval neutrosophic loop of order $n + 1$. $S = \{Z^+ \cup \{0\}\}$ be the semiring. $SL = \{\Sigma s_i [0, aI]\}$ all finite formal sums where $s_i \in S$ and $[0, aI] \in L\}$ is the loop interval neutrosophic semiring of infinite order, non commutative if $m \neq n+1/2$ and non associative.

We can construct any number of such non associative semirings using this class of loops.

Now having seen several types of interval semirings the reader is left with the task of defining properties and results related with them.

Now we proceed onto give the applications of these new structures.

**Chapter Six**

# APPLICATION OF INTERVAL SEMIRINGS

These new structures will find applications in finite interval analysis. Since these structures do not in general involve negative numbers one can use them in real value problems and solutions. Since both associative and non associative interval semirings are defined; in situations were the operations are non associatives these structures can be used.

Semirings which happen to be a less strong than rings and which trivially includes the rings will find applications in place of rings. Likewise semifields of characteristic which contains the class of fields will certainly find applications wherever the field theory is applied. Further this interval semirings and interval semifields will also have lot of applications in stiffness matrices etc. The interlinking of groups, loops, semigroups and groupoids has not only make the structure interesting but stronger.

They can also find applications in performance evaluations and dynamic programming algorithms. For more applications of these interval semirings please refer Semirings and their applications by J.S. Golan (1999)] [3].

**Chapter Seven**

# SUGGESTED PROBLEMS

In this chapter we suggest around 118 problems which will make the interested reader get better grasp of the subject.

1. Obtain some interesting properties enjoyed by interval semirings built using $Z^+ \cup \{0\}$.

2. Is every interval semiring a S-interval semiring?

3. Give an example of an interval semiring which is not an S-interval semiring.

4. Does there exist an interval semiring for which every interval subsemiring is an S-interval subsemiring?

5. Does there exist a S-interval semiring for which no interval subsemiring is an S-interval subsemiring?

6. Obtain some interesting properties enjoyed by S-interval ideals of interval polynomial semiring constructed using $Z^+ \cup \{0\}$.

7. What are the special properties enjoyed by the interval semiring constructed using $Q^+ \cup \{0\}$?

8. Determine a method of solving interval polynomial equations built using $Q^+ \cup \{0\}$.

9. Let $S = \{[0, a] \mid a \in 3Z^+ \cup \{0\}\}$ be an interval semiring.

   a. Does S have ideals and interval subsemirings?
   b. Does S have interval subsemirings which are not ideals?
   c. Can S be a S-semiring?
   
   Justify your claims.

10. Let $P = \{[0, a] \mid a \in Q^+ \cup \{0\}\}$ be a semiring. Study the properties (a), (b) and (c) given in problem (9) for P.

11. Let $T = \{[0, a] \mid a \in R^+ \cup \{0\}\}$ be an interval semiring. Is T an interval semifield?

12. Obtain some interesting properties about interval ideals of an interval semiring.

13. How is an S-interval ideal different from an interval ideal in an interval semiring?

14. Obtain some interesting properties about S-dual interval ideals of an interval semiring.

15. Describe some nice properties enjoyed by S pseudo dual interval ideals of an interval semiring.

16. Define an interval semiring which has zero divisors but no S-anti zero divisors.

17. What are the properties enjoyed by S-zero divisors?

18. Can we say if an interval semiring has S-zero divisors then S has S-antizero divisors?

19. Can we say if an interval semiring has zero divisors then it has S-zero divisors?

20. Does there exist an interval polynomial semiring for which every interval polynomial subsemiring is a S-interval polynomial subsemiring?

21. Find for the interval group ring SG where S = {[0, a] / a ∈ $Z_7$} and G = $S_7$.
    a. What is the order of SG?
    b. Find ideals in SG.
    c. Find interval subsemirings in SG.
    d. Does SG have zero divisors?
    e. Can SG have units?

22. Let G be any group of finite order. S = {[0, a] / a ∈ $Z^+$ ∪ {0}} be an interval semiring, SG the group interval semiring of the group G over S, the interval semiring.
    a. Can SG have zero divisors?
    b. Can SG have S-zero divisors?
    c. Find ideals of SG and characterize them.
    d. Can SG have idempotents or S-idempotent? Justify.

23. Let G = $D_{27}$ = {a, b / $a^2$ = $b^7$ = 1, bab = a} be the dihedral group. S = [0, a] / a ∈ $Q^+$ ∪ {0}} be the interval semiring. Enumerate all the special properties enjoyed by SG, the group interval semiring.

24. Let G = {g / $g^{23}$ = 1} be a cyclic group of order 23. S = {[0, a] / a ∈ $R^+$ ∪ {0}} be an interval semiring. SG the group interval semiring of G over S. Can SG have group interval subsemirings? Can SG have ideals? Is SG a S-semiring? Justify your claim.

25. Let G = {g / $g^{28}$ = 1} be a cyclic group, S = {[0, a] / a ∈ $Z^+$ ∪ {0}} be an interval semiring. Let SG be the group interval semiring.
    a. Is SG a S-ring?
    b. Can SG have ideals?
    c. Does SG have interval subsemirings which are not ideals?
    d. Can SG have idempotents?
    e. Can SG have units?

  f. Can SG have S-anti zero divisors? Justify your answers.

26. Let $S_n$ be the symmetric group of degree n. $S = \{[0, a] / a \in Q^+ \cup \{0\}\}$ be an interval semiring. $SS_n$ be the symmetric group interval semiring of the symmetric group $S_n$ over the interval semiring S.
  a. Find subsemirings of $SS_n$
  b. Can $SS_n$ have idempotents?
  c. Can $SS_n$ have ideals?

27. Let G be a group of finite order say n, and $S_m$ be a symmetric group n < m. $S = \{[0, a] / a \in Q^+ \cup \{0\}\}$ be an interval semiring. SG and $SS_n$ be group interval semiring and symmetric group interval semiring respectively.
  a. Prove there exists an embedding of SG in $S_n$.
  b. $SS_n$ is always non commutative, prove.
  c. $SS_n$ have commutative interval subsemirings, prove!
  d. Find a homomorphism $\eta : SG \to SS_n$ so that $\eta$ is an injective homomorphism.

28. Determine interesting properties enjoyed by $SS_n$ whee $S = \{[0, a] / a \in R^+ \cup \{0\}\}$.

29. Distinguish between $SS_n$ and $TS_n$ where $S = \{[0, a] / a \in Q^+ \cup \{0\}\}$ and $T = \{[0, a] / a \in Z^+ \cup \{0\}\}$.

30. Let $S = \{[0, a] / a \in Z_{20}\}$ be an interval semiring. $G = \langle g / g^{10} = 1\rangle$ be a cyclic group. Find the properties enjoyed by SG. Can SG have zero divisors and S-zero divisors?

31. Let $S = \{[0, a] / a \in Z_{23}\}$ be an interval semiring. $G = \langle g / g^7 = 1\rangle$ be a cyclic group, SG be the group interval semiring. Can SG have zero divisors? Justify!

32. Find all properties enjoyed by $S = \{[0, a] / a \in Z_p$, p a prime$\}$, and $G = \langle g \mid g^q = 1$, q a prime$\rangle$ where S is an interval semiring and G is a group. SG the group interval semiring.

  a. Find all properties enjoyed by SG if p and q are prime.

b. If p and q are not primes but (p, q) = 1 study the additional properties satisfied by SG.
   c. If p and q are non primes with (p, q) = d, what are the properties enjoyed by SG?

   Which of the three group semirings SG given by (a), (b) and (c) are a rich algebraic structure? We say an algebraic structure is rich if it has several properties enjoyed by its elements or substructures.

33. Let G be a group of order n, n, a composite number. S = $\{[0, a] \mid a \in Z_m\}$, m a composite number such that (m, n) = d where d is again a composite number. Let SG be the group interval semiring of the group G over the semiring S. Enumerate all the properties enjoyed by SG.

34. Let S = $\{[0, a] / a \in Z^+ \cup \{0\}\}$ be an interval semiring. G = $\{(D_{2.3} \times S_4 \times P)$ where P is a cyclic group of order 12$\}$ be a group. SG is the group interval semiring of the group G over the interval semiring S.
   a. Is SG a S-semiring?
   b. Does SG contain ideals?
   c. Does SG contain subsemirings which are not ideals?
   d. Find some special properties enjoyed by SG.

35. Let S = $\{Z^+ \cup \{0\}\}$ be a semiring. G = $\{[0, a] / a \in Q^+\}$ be an interval group. SG be the interval group semiring. Determine some distinct properties enjoyed by SG.

36. Let S be the chain lattice $\{0 < a_1 < a_2 < \ldots < a_{12} < 1\}$ which is a semiring. G = $\{[0, a] / a \in Z_{19} \setminus \{0\}\}$ be an interval group under multiplication modulo 19. SG be the interval group semiring of the interval group G over the semiring S. What are the special and distinct properties satisfied by SG?

37. Let S = { 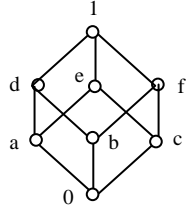

; 0, 1, a, b, c, d, e, f } be a finite semiring. G = {[0, a] / a ∈ $Z_{37} \setminus \{0\}$} be an interval group under multiplication modulo 37. SG be the interval group semiring.

   a. Find ideals in SG.
   b. Find zero divisors in SG
   c. Is SG a S-interval group semiring?
   d. Can SG have idempotents?
   e. Find subsemirings which are not ideals in SG.

38. Let S = { 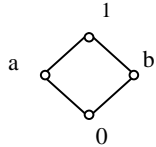

| 0, a, b, 1} be a semiring of order four. G = {[0, a] / a ∈ $Z_{41} \setminus \{0\}$} be an interval group. SG be the interval group ring of the interval group G over the semiring S.
   a. Find zero divisors in SG.
   b. Find the cardinality of SG.
   c. Find substructures of SG.

39. Let S (6) be the symmetric semigroup. F = {[0, a] / a ∈ $Z^+ \cup \{0\}$} be an interval semiring. Let FS (6) be the symmetric semigroup interval semiring of the symmetric semigroup S (6) over the interval semiring F. Determine the special properties enjoyed by FS (6). Prove FS (6) is a S-semigroup interval semiring. Does FS (6) have S-ideals? Does FS (6) have S-pseudo subsemirings?

40. Let S = {[0, a] / a ∈ $Q^+ \cup \{0\}$} be an interval semiring. P = {$Z_{20}$, ×} be a semigroup under multiplication modulo 20.

Find zero divisors, S-zero divisors, S-antizero divisor, idempotents if any in the semigroup interval semiring SP.

41. Let $S = \{[0, a] / a \in Z_{12}\}$ be an interval semiring. $P = \{Z_{24}, \times\}$ be a semigroup under multiplication modulo 24.

    a. Find the order of the semigroup interval semiring SP
    b. Is SP a S-semiring?
    c. Find ideals and subsemirings of SP.

42. Let $G = D_{2.7}$ be the dihedral group of order 14. $S = \{[0, a] / a \in Z_{21}\}$ be an interval semiring. $SD_{2.7} = SG$ be the group interval semiring of the group G over the interval semiring S.
    a. Find the order of SG.
    b. Can SG have ideals?
    c. Can SG have zero divisors?
    d. If SG has zero divisors are they S-zero divisors or S-antizero divisors?
    e. Is SG a S-semiring?
    f. Can SG have S-idempotents?
    g. Give reasons for SG to be only an interval semiring and not an interval semifield. Justify all your answers.

43. Let $G = \{[0, a] \mid a \in Z_{15}\}$ be an interval semigroup under multiplication modulo 15. $S = \{Z^+ \cup \{0\}\}$ be the semiring. SG the interval semigroup semiring of the semigroup G over the semiring S.
    a. Find zero divisors in SG.
    b. Can SG have S-zero divisors?
    c. Can SG have S-antizero divisors?
    d. Is SG a S-semiring?
    e. Can SG have S-ideals?
    f. Does SG have S-subsemiring?

44. Let $G = \{[0, a] \mid a \in Z_{19} \setminus \{0\}\}$ be an interval group. $S = \{Z^+ \cup \{0\}\}$ be the semiring. SG the interval group semiring of the interval group G over the semiring S. Find interesting properties of SG.

45. Let $S = \{[0, a] \mid a \in Z_{12}\}$ be an interval semigroup. $F = \{Q^+ \cup \{0\}\}$ be a semiring. FS the interval semigroup semiring of the interval semigroup S over the semiring F.
    a. Can FS have zero divisors?
    b. Can FS be a S-semiring?
    c. Can FS have ideal? Justify your claim.

46. Let $G = \{[0, a] \mid a \in Z_{11} \setminus \{0\}\}$ be an interval group. $S = \{Z^+ \cup \{0\}\}$ be a semiring. SG the interval group semiring.
    a. Is SG a S-semiring?
    b. Is SG a semifield?
    c. Obtain proper subsemirings of SG.

47. Let $G = \{[0, a] \mid a \in Z_{18}\}$ be an interval semigroup and $S = Z^+ \cup \{0\}$ be a semiring. SG the interval semigroup semiring of the interval semigroup G over the semiring S. Study the problems (a), (b) and (c) given in problem 46.

48. Let $S = \{[0, a] \mid a \in Z_{20}\}$ be an interval semiring. $G = \langle g \mid g^{12} = 1 \rangle$ be a group. SG the group interval semiring. Obtain some interesting properties enjoyed by SG.

49. Let $S = \{[0, a] \mid a \in Z_{20}\}$ be an interval semigroup. $F = Z^+ \cup \{0\}$ be a semiring. FS be the interval semigroup semiring of the interval semigroup S over the semiring S.
    a. Find zero divisors in FS.
    b. Can FS be a S-semiring?
    c. Does FS have S-ideals?
    d. Does FS have subsemirings which are not S-subsemirings?

50. Let $S = \{R^+ \cup \{0\}\}$ be a interval semiring. $G = \langle g \mid g^{24} = 1 \rangle$ be a group SG be the group interval semiring.
    a. Can SG have zero divisors?
    b. Can SG have idempotents?
    c. Is SG a S-semifield?
    d. Can SG have S-ideals?

51. Let $S = \{[0, a] / a \in Z^+ \cup \{0\}\}$ and $P = \{[0, x] / x \in Z_{15}\}$ be the interval semiring and interval semigroup respectively.

SP be the interval semigroup interval semiring. Discuss the special properties enjoyed by this new structure.

52. What are groupoid interval semirings? Give examples of finite and infinite groupoid interval semirings.

53. Let $G = \{[0, a] / a \in Z_{12}, *, (3, 7)\}$ be an interval groupoid. $S = \{Z^+ \cup \{0\}\}$ be a semiring. SG be the interval groupoid semiring. Find the special properties enjoyed by these structures.

54. Let $G = \{[0, a] / a \in Z_{11}, *, (3, 4)\}$ be an interval groupoid. $S = \{Q^+ \cup \{0\}\}$ be a semiring, SG be the interval groupoid semiring of G over S.
    a. Find zero divisors of SG if any.
    b. Find S-ideal of SG if any.
    c. Does SG contain S-subsemiring?
    d. Show SG is non associative.
    e. Does SG satisfy any special identity?

55. Let $G = \{[0, a], *, (3, 11), a \in Z^+ \cup \{0\}\}$ be an interval groupoid. $S = \{Z^+ \cup \{0\}\}$ be a semiring; SG is an interval groupoid semiring of infinite order which is a non associative semiring.
    a. SG has no zero divisors, prove.
    b. Can SG have units and idempotents?
    c. Is SG a S-semiring?
    d. Find subsemirings in SG.
    e. Can SG have S-subsemirings?
    f. Can SG have ideals and S-ideals?
    g. Does SG have S-subsemirings which are not ideals?

56. Let $G = \{[0, a] / a \in Q^+ \cup \{0\}, *, (3/7, 2/5)\}$ be an interval groupoid. $S = \{Z^+ \cup \{0\}\}$ be a semiring. SG the interval groupoid semiring of the interval groupoid G over the semiring S. Study all the questions 1 to 7 given in problem 55.

57. Let $G = \{[0, a] / a \in Z_7, *, (3, 2)\}$ be an interval groupoid. $S = Z^+ \cup \{0\}$ be a semiring. SG be the interval groupoid

semiring. Can SG be a S-semiring? Can SG have S-zero divisors?

58. Let G = {[0, a] / a ∈ $Z_{12}$, (3, 6), *} be an interval groupoid. S = $Z^+$ ∪ {0} be a semiring. SG be the interval groupoid semiring. Describe all the distinct properties enjoyed by SG.

59. Let S = {[0, a] / a ∈ $Z^+$ ∪ {0}} be an interval semiring. G = {a ∈ $Z_{12}$, *, (3, 7)} be a groupoid. SG the groupoid interval semiring. Find S-zero divisors and zero divisors if any in SG.

60. Let G = {$Z_{19}$, *, (2, 3)} be a groupoid. S = {[0, a] / a ∈ $Q^+$ ∪ {0}} be an interval semiring. SG be the groupoid interval semiring.
    a. Find S-subsemirings in SG.
    b. Can SG have idempotents?
    c. Does SG have S-ideals?
    d. Can SG have S-pseudo ideals?
    e. Can SG have S-anti zero divisors?

61. Let G = {$Z_{214}$, *, (0, 3)} be a groupoid. S = {[0, a] / a ∈ $Z_2$} be an interval semiring. SG the groupoid interval semiring of G over S.
    a. Find the order of SG.
    b. Is SG a S-ring?
    c. Can SG have an S-ideal?
    d. Can SG have S-subsemirings which are not S-ideals?
    e. Can SG have S-zero divisors?

62. Let G = ⟨$Z_6$, *, (0, 5)⟩ be the groupoid. S = {[0, a] / a ∈ $Z_{12}$} be the interval semiring. SG be the groupoid interval semiring. Find zero divisors in SG? Is every zero divisor a S-zero divisor in SG?

63. Let G = {$Z_{14}$, *, (0, 7)} be a groupoid. S = {[0, a] / a ∈ $Z^+$ ∪ {0}} be an interval semiring. SG be the groupoid interval semiring. Determine some properties enjoyed by SG.

64. Can there be a groupoid interval semiring SG of order 29?

65. Can there be a groupoid interval semiring SG of order p, p a prime? Justify your claim.

66. Can there be a groupoid interval semiring SG of order $2^{14}$?

67. Let G = {$Z_{12}$, *, (3, 6)} be a groupoid. S = {$Z^+ \cup \{0\}$} be the interval semiring. SG the groupoid interval semiring of G over S. Does SG satisfy any of the classical identities for groupoids?

68. Let G = {[0, a] | a ∈ $Z_{12}$, *, (3, 4); [0, a] * [0, b] = [0, 3a+4b (mod 12)]} be an interval groupoid. S = {$0 < a_1 < a_2 < a_3 < \ldots < a_{10} < 1$ } be a finite semiring. Study the interval groupoid semiring SG.

69. Let G = {[0, a] | a ∈ $Z_{15}$, (0, 7)} be an interval groupoid. S = $C_n$ be a chain lattice where n = 30. SG be the interval groupoid semiring.
    a. Can SG have S-subsemiring?
    b. Can SG have ideals?
    c. Can SG have S-ideals?
    d. Can SG have zero divisors?
    e. Obtain some interesting properties enjoyed by SG.
    f. Find the order of SG.
    g. Does the order of subsemiring divide the order of SG?

70. Let S = {[0, a] / a ∈ $Q^+ \cup \{0\}$} be an interval semiring. G = $L_{17}$ (5) be the loop of order 18. SG be the loop interval semiring of the loop G over the interval semiring S.
    a. Can SG have ideals?
    b. Is SG a S-semiring?
    c. Can SG have S-subsemirings which are not ideals of SG?
    d. Does SG contain zero divisors?
    e. Can SG have S-zero divisors?
    f. Prove SG has idempotents.

71. Let G = {$L_{11}$ (5)} be a loop of order 12. S = {[0, a] / a ∈ $Z_{12}$} be an interval semiring. Let SG be the loop interval semiring of G over the interval semiring S.

a. Find the order of SG.
b. Does SG have zero divisors?
c. Is every zero divisors of SG a S-zero divisor of SG.
d. Find S-ideals in SG (if any).
e. Can SG be S-semiring?
f. Can SG have S-subsemirings which are not S-ideals of SG?

72. Let $G = L_{15}(8)$ be a loop of order 16. $S = \{[0, a] / a \in Z_{16}\}$ be an interval semiring; SG be the loop interval semiring of the loop G over the interval semiring S.
    a. Find the order of SG.
    b. Is SG commutative?
    c. Does SG have S-ideals?
    d. Does the order of every subsemiring divide the order of SG?
    e. Can SG have subsemirings which are not ideals?
    f. Can SG have zero divisors?
    g. Does SG satisfy any classical identity for loops?

73. Let $G = \{[[0, a] / a \in L_5(3)\}$ be an interval loop. $S = Z^+ \cup \{0\}$ be a semiring. SG be the interval loop semiring of the interval loop G over the semiring S. Obtain some special properties enjoyed by SG.

74. Let $G = L_5(3)$ be a loop and $S = \{[0, a] \mid a \in Z^+ \cup \{0\}\}$ be an interval semiring SG be the loop interval semiring. Are semirings in problems (73) and (74) isomorphic?

75. Determine some interesting properties enjoyed by infinite interval loop semirings.

76. Derive some interesting properties enjoyed by loop interval semiring of infinite order.

77. What are the distinct properties satisfied by interval loop semirings?

78. Let $L = \{[0, x] \mid x \in L_{21}(17)\}$ be an interval loop. $S = \{[0, n] \mid n \in Z_{22}\}$ be an interval semiring. LS the interval loop interval semiring. Find the order of LS. Discuss the properties enjoyed by LS.

79. Enumerate some properties enjoyed by non associative interval semiring.

80. Construct non associative interval semirings other than those obtained by groupoid interval semirings, interval groupoid semirings, loop interval semirings and interval loop interval semirings.

81. Give an example of a groupoid interval semiring which is not a S-semiring.

82. Give an example of a loop interval semiring which is not a S-semiring.

83. Prove the class of interval semirings $SL_n(m)$; $S = \{[0, a]$ where a is in $Z^+ \cup \{0\}\}$ and $L_n(m)$ defined in [ ] are S-loop interval semirings.

84. Let $G = \{Z_{41}, *, (3, 12)\}$ be a groupoid. $S = Z^+ \cup \{0\}$ be a semiring. SG be the groupoid semiring. Consider $G_1 = \{[0, a] / a \in Z_{41}, *, (3, 12)\}$ be the interval groupoid. $SG_1$ be the interval groupoid semiring of $G_1$ over S.
    a. Is $SG \cong SG_1$?
    b. Is SG a S-semiring?
    c. Can $SG_1$ be a S-semiring?

85. Let $G = \{[0, aI] / a \in Z_{20}, *, (3, 7)\}$ be a neutrosophic interval groupoid. $S = \{Z^+ \cup \{0\}\}$ be the semiring. Study the properties of SG, the neutrosophic interval groupoid semiring.
    Suppose $H = \{Z_{20}, *, (3, 7)\}$ be a groupoid and $S = \{[0, aI] / a \in Z^+ \cup \{0\}\}$ be the neutrosophic interval semiring. Study the properties of SH the groupoid neutrosophic interval semiring.
    Compare SH and SG and the derive all the common properties enjoyed by them.

86. Study all the properties enjoyed by $S = \{[0, a + bI] / a, b \in Q^+ \cup \{0\}\}$ the neutrosophic interval semiring.

87. What are the probable applications of this new class of semirings?

88. Let $L = \{[0, aI] / a \in \{eI, I, 2I, \ldots, 15I\}, *, m = 8, [0, aI] * [0, bI] = [0, 8bI + 7aI \pmod{15}]\}$ be the pure neutrosophic interval loop. Derive all the properties associated with L.
    Suppose $T = \{[0, a + bI] / a, b \in \{e, 1, 2, \ldots, 15\}, *, m = 8\}$ be the neutrosophic interval loop. Compare the two loops L and T. What is the cardinality of T? Does T have interval subloops?

89. Let $G = \{[0, a+bI] / a, b \in Z_{32}\}$ be the neutrosophic interval semigroup. Discuss all the properties related with G.
    a. What is the cardinality of G?
    b. Does G have zero divisors?
    c. Can G have idempotents?
    d. Is G a S-semigroup?

90. Let $H = \{[0, a + bI] / a, b \in 31\}$ be the neutrosophic interval semiring. Study the properties (a) to (d) mentioned in problem 89 for H.

91. Let $P = \{[0, a + bI] / a, b \in Z^+ \cup \{0\}\}$ be a neutrosophic interval semigroup.
    a. Is P a S-semigroup?
    b. Can P have ideals?
    c. Can P have S-subsemigroups which are not S-ideals?
    d. Can P have zero divisors or idempotents?

92. Let $W = \{[0, a + bI] / a, b \in R^+ \cup \{0\}\}$ be a neutrosophic interval semigroup. Study properties (a) to (d) mentioned in problem 91 for P incase of W.

93. Let $F = \{[0, a + bI] \mid a, b \in Z^+ \cup \{0\}\}$ be a neutrosophic interval semiring.
    a. Is F a semifield?
    b. Is F a S-semiring?

c. Can F have S-ideals?

94. Let $K = \{[0, a + bI] \mid a, b \in R^+ \cup \{0\}\}$ be a neutrosophic interval semiring. Study properties (a), (b) and (c) mentioned in problem (93) for K.

95. Let $S = \{[0, a + bI] \mid a, b \in Z_{28}\}$ be a neutrosophic interval semigroup. $F = \{C_{10}; 0 < a_1 < \ldots < a_8 < 1\}$ be a semiring. FS the neutrosophic interval semigroup semiring.
    a. What is the cardinality of FS?
    b. Find S-ideals in any in FS.
    c. Can FS have S-zero divisors?
    d. Can FS have idempotents?
    e. Is FS a S-semiring?

96. Let $P = \{[0, aI] \mid a \in L_n(m), *\}$ be an interval neutrosophic loop. Study and characterize the properties enjoyed by P. Suppose $S = \{C_7, 0 < a_1 < \ldots < a_5 < 1\}$ be the semiring. SP the interval neutrosophic loop semiring.
    a. Prove order of SP is finite.
    b. Is SP a S-semiring?
    c. Is SP a semifield?
    d. Can SP have S-zero divisors?
    e. What are the special properties enjoyed by SP?

97. Let $L = \{[0, a + bI] \mid a, b \in L_{21}(11)\}$ be a neutrosophic interval loop.
    a. Find the order of L?
    b. Is L a S-loop?

98. Let $L = \{[0, a + bI] \mid a, b \in L_{23}(12)\}$ be a neutrosophic interval loop. $S = \{0 < a_1 < \ldots < a_{16} < 1\}$ be the semiring. SL the neutrosophic interval semiring of L over S.
    a. What is the order of SL?
    b. Is SL a S-semiring?
    c. Does SL have zero divisors and S-zero divisors?
    d. What are the S-ideals in SL? (if any).

99. Let $G = \{[0, aI] \mid a \in Z_{17} \setminus \{0\}\}$ be the pure neutrosophic interval group. $S = \{0 < a_1 < a_2 < \ldots < a_{10} < 1\}$ be the semiring; SG the group interval neutrosophic semiring.
    a. Find the order of SG?
    b. Is SG a S-semiring?
    c. Can SG have S-ideals?
    d. Does SG have S-zero divisors?
    e. Can SG have S-subsemirings which are not S-ideals?

100. Let $G = \{[0, aI + b] \mid b, a \in Z_{23} \setminus \{0\}\}$ be a neutrosophic interval group. $S = Z^+ \cup \{0\}$ be the semiring. SG be the group interval neutrosophic semiring. Study (b) to (e) given in problem (99) for this SG.

101. Let $S = \{[0, a + bI] \mid a, b \in Z_{240}\}$ be a neutrosophic interval semigroup. $K = \{0 < a_1 < \ldots < a_{100} < 1\}$ be a semiring. KS the neutrosophic interval semigroup semiring. Study the properties (a) to (e) given in problem 99.

102. Let $S = \{[0, a + bI] \mid a, b \in Z_{24}\}$ be a neutrosophic interval semigroup. $K = \{[0, a + bI] / a, b \in Z^+ \cup \{0\}\}$ be the interval neutrosophic semiring. KS the semigroup neutrosophic interval semiring. Obtain the properties enjoyed by KS.

103. Let $W = \{[0, a+bI] \mid a, b \in Z_{26}, *, (3, 11)\}$ be a neutrosophic interval groupoid. $S = \{0 < a_1 < a_2 < a_3 < 1\}$ be a semiring. SW be the neutrosophic interval groupoid semiring.
    a. Find the order of SW.
    b. Is SW a S-semiring?
    c. Does SW have S-zero divisors?
    d. Can SW have S-ideals?
    e. Is S a non associative S-division ring?

104. Let $G = \{[0, a +bI] / a, b \in Z_{20}, *, (10, 12)\}$ be a neutrosophic interval groupoid. $S = \{[0, a + bI] / a, b \in Q^+ \cup \{0\}\}$ be an neutrosophic interval semiring. SG be the neutrosophic interval semiring. Study the questions (b) to (e) mentioned in problem (103).

105. Let $G = \{[0, a + bI] \mid a, b \in Z_{46}\}$ be a neutrosophic interval semigroup. $S = \{[0, a + bI] \mid a, b \in Q^+ \cup \{0\}\}$ be an interval neutrosophic semiring. SG be the neutrosophic interval semigroup interval semiring. Enumerate all the properties enjoyed by SG.

106. Let $S = \left\{\sum_{i=0}^{\infty}[0, a + bI]x^i \mid a, b \in Z^+ \cup \{0\}\right\}$ be an interval neutrosophic polynomial semiring. Describe all the properties enjoyed by S.

107. Let $P = \left\{\sum_{i=0}^{\infty}[0, aI]x^i \mid aI \in Q^+ \cup \{0\}\right\}$ be an interval neutrosophic polynomial semiring.

    a. Can P have S-ideals?
    b. Is P a strict semiring?
    c. Can P be a semifield?
    Justify all your claims.

108. Let P = {all 10 × 10 neutrosophic interval matrices with intervals of the form $[0, a + bI] \mid a, b \in Z_7$} be a neutrosophic interval matrix semiring.
    a. Find the order of P.
    b. Is P a semifield?
    c. Find S-zero divisors in P.
    d. Find S-ideals if any in P.
    e. Can P have S-subsemirings which are not S-ideals? Justify your claim.

109. Let $S = \{([0, a_1I], [0, a_2I], \ldots, [0, a_9I]) \mid a \in Z_{45}\}$ be a neutrosophic interval row matrix semiring. Study properties (a) to (e) mentioned in problem 108.

110. Let T = {All 5 × 5 upper triangular neutrosophic interval matrices with intervals of the form $[0, a + bI]$ ; $a, b \in Z_{12}$} be a neutrosophic interval matrix semiring study properties (a) to (e) described in problem (108).

111. Let S = {([0, a + bI]), …, [0, d+eI] | a, b, d, e ∈ $Z_{24}$} be a 1 × 5 row neutrosophic interval matrix semiring. Study properties (a) to (e) mentioned in problem 108.

112. Let P = {All 6 × 6 pure neutrosophic interval matrices with intervals of the form [0, aI] where a ∈ $Q^+ \cup \{0\}$} be the pure neutrosophic interval matrix semiring.
    a. Is P a S-semiring?
    b. Can P be a semifield?
    c. Does P have S-zero divisors? Justify your claim!

113. Let W = $\left\{ \sum_{i=0}^{8} [0, a+bi]x^i \,\middle|\, x^9 = 1, a, b \in Z_5 \right\}$ be a neutrosophic interval polynomial semiring.
    a. Find the order of W.
    b. Can W be a S-semiring?
    c. Find S-ideals in W.
    d. Can W be a semifield?
    e. Find S-zero divisors in W.

114. Let P = $\left\{ \sum_{i=0}^{5} [0, a+bi]x^i \,\middle|\, x^6 = 1, a, b \in Z_6 \right\}$ be a neutrosophic interval polynomial semiring study properties (i) to (v) mentioned in problem 113. Is their any difference between W and P structure wise?

115. Obtain any interesting properties enjoyed by neutrosophic interval polynomial semirings built using $Z^+ \cup \{0\}$.

116. Study some special properties enjoyed by pure neutrosophic interval 5 × 5 matrix semiring built using $Z_n$, n finite and n not a prime.

117. Differentiate between neutrosophic interval semirings and pure neutrosophic interval semirings.

118. Let G = {All 4 × 4 neutrosophic interval matrices with intervals of the form [0, a + bI] | a, b ∈ $Z_{12}$} be a neutrosophic interval matrix semiring.

a. Find order of G.
b. Is G a S-semiring?
c. Find S-subsemirings which are not S-ideals of G.
d. Find S-zero divisors (if any).
e. Find zero divisors which are not S-zero divisors (if any).

# FURTHER READING

# INDEX

**W**

# ABOUT THE AUTHORS

**Dr.W.B.Vasantha Kandasamy** is an Associate Professor in the Department of Mathematics, Indian Institute of Technology Madras, Chennai. In the past decade she has guided 13 Ph.D. scholars in the different fields of non-associative algebras, algebraic coding theory, transportation theory, fuzzy groups, and applications of fuzzy theory of the problems faced in chemical industries and cement industries. She has to her credit 646 research papers. She has guided over 68 M.Sc. and M.Tech. projects. She has worked in collaboration projects with the Indian Space Research Organization and with the Tamil Nadu State AIDS Control Society. She is presently working on a research project funded by the Board of Research in Nuclear Sciences, Government of India. This is her 53$^{rd}$ book.

On India's 60th Independence Day, Dr.Vasantha was conferred the Kalpana Chawla Award for Courage and Daring Enterprise by the State Government of Tamil Nadu in recognition of her sustained fight for social justice in the Indian Institute of Technology (IIT) Madras and for her contribution to mathematics. The award, instituted in the memory of Indian-American astronaut Kalpana Chawla who died aboard Space Shuttle Columbia, carried a cash prize of five lakh rupees (the highest prize-money for any Indian award) and a gold medal.
She can be contacted at vasanthakandasamy@gmail.com
Web Site: http://mat.iitm.ac.in/home/wbv/public_html/
or http://www.vasantha.in

**Dr. Florentin Smarandache** is a Professor of Mathematics at the University of New Mexico in USA. He published over 75 books and 150 articles and notes in mathematics, physics, philosophy, psychology, rebus, literature.

In mathematics his research is in number theory, non-Euclidean geometry, synthetic geometry, algebraic structures, statistics, neutrosophic logic and set (generalizations of fuzzy logic and set respectively), neutrosophic probability (generalization of classical and imprecise probability). Also, small contributions to nuclear and particle physics, information fusion, neutrosophy (a generalization of dialectics), law of sensations and stimuli, etc. He can be contacted at smarand@unm.edu